\documentclass[a4paper,reqno,12pt]{amsart}

\usepackage{fullpage}

\usepackage[plainpages=false,colorlinks,linkcolor=black,citecolor=black,urlcolor=black,breaklinks]{hyperref}

\usepackage{amsmath,amsthm,amssymb}
\usepackage{mathrsfs}
\usepackage[shortlabels]{enumitem}
\usepackage{graphicx}
\usepackage[font=small]{caption}
\usepackage{xcolor}
\usepackage{url}

\newcommand{\N}{\mathbb{N}}
\newcommand{\R}{{\mathbb{R}}}

\newcommand{\C}{{\mathbb{C}}}
\newcommand{\Z}{{\mathbb{Z}}}
\newcommand{\D}{{\mathbb{D}}}
\newcommand{\dd}{{{\rm d}}}
\newcommand{\ii}{{\rm i}}

\newcommand{\ov}{\overline}
\newcommand\wt{\widetilde}

	\newcommand{\la}{\lambda}
\newcommand{\La}{\Lambda}
\newcommand{\eps}{\varepsilon}

\newcommand{\spp}{\sigma_{\rm p}}

\newcommand{\Dom}{{\operatorname{Dom}}}
\newcommand{\Ker}{{\operatorname{Ker}}}

\newcommand{\Rank}{{\operatorname{rank}}}
\renewcommand{\Re}{\operatorname{Re}}
\renewcommand{\Im}{\operatorname{Im}}

\newcommand{\dist}{\operatorname{dist}}
\newcommand{\sgn}{\operatorname{sgn}}
\newcommand{\supp}{\operatorname{supp}}

\newcommand{\Num}{\operatorname{Num}}

\newcommand{\BigO}{\mathcal{O}}

\newcommand{\Dt}{-\frac{\dd^2}{\dd x^2}}
\newcommand{\Dtp}{\frac{\dd^2}{\dd x^2}}
\newcommand{\Do}{\frac{\dd}{\dd x}}

\theoremstyle{plain}

\newtheorem{theorem}{Theorem}[section]
\newtheorem{lemma}[theorem]{Lemma}

\newtheorem{proposition}[theorem]{Proposition}
\newtheorem{corollary}[theorem]{Corollary}

\theoremstyle{definition}

\newtheorem{remark}[theorem]{Remark}
\newtheorem{asm-sec}[theorem]{Assumption}

\newcommand\cD{\mathcal D}

\newcommand\cH{\mathcal H}

\newcommand\cJ{\mathcal J}

\newcommand\cP{\mathcal P}

\newcommand\cS{\mathcal S}
\newcommand\cT{\mathcal T}

\newcommand\frs{\mathfrak s}

\newcommand\frU{\mathfrak U}

\numberwithin{equation}{section}
\numberwithin{figure}{section}

\usepackage[
backend=biber
,style=alphabetic
,firstinits=true
,maxnames=4 
,sorting=nyvt
,isbn=false
]{biblatex}

\renewbibmacro{in:}{} 

\DeclareFieldFormat[article,inbook,incollection,inproceedings,patent,thesis,unpublished,misc]{title}{#1}
\DeclareFieldFormat{pages}{#1}

\AtEveryBibitem{
	\ifentrytype{article}{
		\clearfield{number}
	}{}
}

\begin{document}
\title[Spectral projections of a complex anharmonic oscillator]{Spectral projections of an anharmonic oscillator with complex polynomial potential}

\author{Boris Mityagin}

\address[Boris Mityagin]{Department of Mathematics, The Ohio State University, 231 West 18th Ave, Columbus, OH 43210, USA}
\email{mityagin.1@osu.edu, boris.mityagin@gmail.com}

\author{Petr Siegl}

\address[Petr Siegl]{Institute of Applied Mathematics, Graz University of Technology, Steyrergasse 30, 8010 Graz, Austria}
\email{siegl@tugraz.at}

%\thanks{The author ...}



\date{\today}

\begin{abstract}
For a broad class of polynomial potentials $V$, with an important and instructive representative being $V(x) = x^{2a} + \ii x^b$, $x \in \R$, $a, b \in \N$, we show that the system of spectral projections $\{P_n\}_n$ of an anharmonic operator $L = - (\dd/ \dd x)^2 + V(x)$ does not generate a (Riesz) basis 
in $L^2(\R)$ if $a - 1 < b < 2a$. Moreover, for $\sigma = [b - (a - 1)]/(1 + a)$ and $\gamma > 0$ small enough, 
$\limsup_n \|P_n\|/ \exp(\gamma n^\sigma) = \infty$. Proofs are based on two groups of results which are of great interest
on their own: (a) relationship between behavior (growth) of the norms of projections $\|P_n\|$ and of the resolvent $\|(z - L)^{-1}\|$  outside of the spectrum $\sigma(L)$; (b) partial fraction decompositions of special
meromorphic functions $1/F$ where $F(w) = \prod_{k=1}^\infty \left( 1 + \frac{w}{a_k} \right)$, $a_{k+1} \geq a_k>0$, $k \in \N$, and the generalization of the first resolvent identity.
\end{abstract}

\maketitle

%\newpage 
\tableofcontents

%\newpage

\section{Introduction}
\label{sec:intro}

A special case of the differential operator $L$ in $L^2(\R)$ analyzed in this paper reads
\begin{equation}\label{L.intro}
L = \Dt + x^{2a} + \ii x^b, \qquad a,b \in \N,  \ b<2a,
\end{equation}
cf.~\eqref{L.def}.
It is an unbounded m-accretive operator with compact resolvent in a separable Hilbert space. K.~C.~Shin showed that the eigenvalues $\{\la_n\}_{n \in \N}$ of $L$ are simple for all $n>N_0 > 0$ and they satisfy as $n \to \infty$
\begin{equation}\label{la.L.intro}
	\begin{aligned}
		\la_n &= (d n)^\kappa(1+o(1)), 
		\\
		\la_{n+1}-\la_n & = \kappa d  (dn)^{\kappa-1} \left(1 + o(1)\right),   \quad \kappa = \frac{2a}{a+1} \in [1,2),
	\end{aligned}
\end{equation}
where $d>0$ is as in \eqref{d.L.def}, see \cite{Shin-2010-35} or Proposition~\ref{prop:Shin} below; cf.~also~\cite{Caliceti-2006-39,Caliceti-2014-19}.

The operator $L$ has a structure $L=A+B$ where $A = - (\dd/ \dd x)^2 + x^{2a}$ is self-adjoint and $B = \ii x^b$ is its relatively bounded perturbation. The system of orthogonal spectral projections of $A$ is complete, however, since $B$ is non-symmetric, the orthogonality of spectral projections of $L$ is lost. Nonetheless, the system of spectral projections $S_0 \cup \{P_n\}_{n>N_0}$ of $L$ is complete (in fact, even for $b \geq 2a$); here 
\begin{equation}\label{Pn.L.infro}
	P_n = \frac{1}{2 \pi \ii} \int_{|z-\la_n|=\delta_{\la_n}} (z-L)^{-1} \, \dd z, \quad \delta_\la = \frac12 \dist(\la,\sigma(L)\setminus\{\la\}), \quad n > N_0,
\end{equation}
and $S_0$ is the projection with $\Rank (S_0) = N_0$ related to the remaining $N_0$ eigenvalues, i.e.~, with $R_0>0$ such that $|\la_{N_0}|<R_0<|\la_{N_0+1}|$ 
\begin{equation}\label{S0.L.intro}
	S_0 = \frac{1}{2 \pi \ii} \int_{|z|=R_0} (z-L)^{-1} \, \dd z.
\end{equation}
It can be expected that the basis property of the spectral projections of $L$ depends on the strength of the perturbation and the goal is to find the critical point between the basis and no basis property. 

It is known that if $b < a-1$, then the system of spectral projections of $L$ has the Riesz property, see \cite[Thm.~6.4]{Mityagin-2019-139}. Consequently, the resolvent of $L$ satisfies (with suitable constants $M,R>0$)
\begin{equation}
\|(z-L)^{-1}\| \leq \frac{M}{\dist(z,\sigma(L))}, \quad z \in \rho(L), \quad |z| \geq R.
\end{equation}
For a collection of recent results on the Riesz basis property for abstract and differential operators (Dirac, Hill, Schr\"odinger), see \cite{Adduci-2012-10,Adduci-2012-73,Djakov-2011-351,Djakov-2010-283,Djakov-2012-263,Mityagin-2016-106,Mityagin-2019-139,Motovilov-2017-8,Motovilov-2019-53,Shkalikov-2010-269,Shkalikov-2016-71,Wyss-2010-258}.

The resolvent behavior of $L$ changes drastically if $b>a-1$. In detail, if $b$ is odd, for any fixed $\beta \in \R$ and any $N \in \N$, there exists $C>0$ such that for all $\alpha + \ii \beta \in \rho(L)$ satisfying $\alpha \geq 1$, we have 
\begin{equation}\label{L.res.lb.intro}
\|(\alpha + \ii \beta-L)^{-1}\| \geq C \alpha^N,
\end{equation}
see~\cite[Thm.~3.7, Ex.~3.8]{Krejcirik-2019-276}; the case of even $b$ is slightly more complicated, but qualitatively similar. In fact in both cases, the resolvent growth is exponential; in detail for $b$ odd, there exists $\eta>0$ such that 
\begin{equation}\label{L.res.lb.exp.intro}
	\|(\alpha + \ii \beta-L)^{-1}\| \geq \exp(\eta \alpha^\tau), \quad \tau = \frac{b-(a-1)}{2a},
\end{equation}
see Theorem~\ref{thm:L.res.lb}. The latter is a generalization \cite[Thm.~7]{Krejcirik-2015-56} dealing with the case $a=b=1$ and following up on the semi-classical methods, see~e.g.~\cite[Chap.~2]{Dimassi-1999}.
Due to the known eigenvalue asymptotics \eqref{la.L.intro}, both \eqref{L.res.lb.intro} and \eqref{L.res.lb.exp.intro} can be used to show that the system of spectral projections of $L$ does not have the Riesz property. 

However, it turns out that for $b>a-1$ the spectral projections of $L$ do not have even a basis property. More precisely, as our main result for $L$, we show that for a small enough $\gamma>0$ we have
\begin{equation}\label{Pn.L.limsup.intro}
	\limsup_{n \to \infty} \frac{\|P_n\|}{\exp (\gamma n^\sigma)} = \infty,
\end{equation}
where
\begin{equation}\label{sigma.L.intro}
	\sigma = \frac{b-(a-1)}{1+a} \in (0,1).
\end{equation}
In fact, like in the special case $a=b=1$ studied in \cite{Mityagin-2017-272}, \emph{we expect that \eqref{Pn.L.limsup.intro} holds with $\lim$ instead of $\limsup$ and also that $\|P_n\|$ actually grow with the exponential order $\sigma$,} see Remark~\ref{rem:L.shifted.HO}. \emph{We also expect that the transition in the basis properties at $b=a-1$ holds for more general than polynomial potentials,} see Subsection~\ref{ssec:no.basis} where also the Abel basis property for $b=a-1$ is discussed.

Our method to show \eqref{Pn.L.limsup.intro} is of independent interest and applies to abstract m-accretive operators $T$ with compact resolvent and with complete system of spectral projections (the perturbative structure $A+B$ is not needed). Our strategy is to propagate an exponential bound of order $\varrho>0$ on the spectral projections to an exponential bound on the resolvent with the same order. For differential operators with quadratic symbols, such behavior was observed earlier by J.~Viola, cf.~\cite[2013-14-ICMS]{nsa-portal-op}.

We explain the strategy in the case inspired by $L$ where $\varrho < 1/2$ and the eigenvalues are eventually simple and satisfy $\la_n = (d n)^\kappa(1+o(1))$ with $d,\kappa>0$ as $n \to \infty$.  In this case, Theorem~\ref{thm:T} shows that if the spectral projections $P_n$ of $T$ satisfy 
\begin{equation}\label{Pn.bd.intro}
\|P_n\| \leq \exp ( \gamma |\la_n|^\varrho (1+o(1))), \quad n \to \infty,
\end{equation}
for $\gamma>0$ and $\varrho <1/2$, then for any $\delta>0$ there exists $C>0$ such that for $z \in \rho(T)$ with $\Re z \geq 0$
\begin{equation}\label{T.res.intro}
\|(z-T)^{-1}\| \leq \frac{C} {\dist(z,\sigma(T))} \exp \big((1+\delta) \gamma \cos (\varrho \arg z) |z|^\varrho \big), \quad z \to \infty.
\end{equation}
Thus the exponential order $\varrho$ is propagated from the projection bounds to the resolvent; moreover, the information on the type $\gamma$ is transported as well. Theorem~\ref{thm:L} states this result for $L$ and the claim \eqref{Pn.L.limsup.intro} on spectral projections follows by combining \eqref{Pn.bd.intro} and \eqref{T.res.intro} with \eqref{la.L.intro} and \eqref{L.res.lb.exp.intro}, see Theorem~\ref{thm:L.Pn} for details. We note that the assumptions of Theorem~\ref{thm:T}, which is our second main result, are more flexible than the case described above. Only $|\la_n| \geq (d n)^\kappa$, $n \in \N$, is required and $\varrho \geq 1/2$ is included under additional assumptions on the resolvent growth on certain rays in $\C$.

The proof of Theorem~\ref{thm:T} relies on a remarkable new identity for the resolvent. For $T$ with simple eigenvalues it reads 
\begin{equation}\label{res.id.intro}
B_z(T) = \frac{1}{F(z)}(z-T)^{-1} + \sum_{n=1}^\infty \frac{1}{(z+a_n) F'(-a_n)} (a_n+T)^{-1}, \quad z \in \rho(T),
\end{equation}
where 
\begin{equation}\label{B.F.intro}
B_z(w) = \frac{1}{F(w)(z-w)}, \quad F(w) = \prod_{k=1}^{\infty} \left(1 + \frac{w}{a_k}\right), \quad w \in \C \setminus \{z\}, 
\end{equation}
and
\begin{equation}\label{ak.intro}
a_k = \left( \frac{k}{\nu} \right)^{\frac 1 \varrho}, \quad \nu>0, \ \varrho < 1/2, \ k \in \N.
\end{equation}
The key step is to realize that \eqref{res.id.intro} is an operator version of the partial fraction decomposition of the meromorphic function $B_z(\cdot)$. 

The identity \eqref{res.id.intro} indeed allows for propagating the exponential bound on the spectral projections to the resolvent. Namely, the entire function $F$, studied in particular by E.~C.~Titchmarsh in \cite{Titchmarsh-1927-26}, obeys for any $\theta \in (-\pi,\pi)$ the asymptotics 
\begin{equation}\label{F.asym.intro}
\log|F(r e^{\ii \theta})| = \pi \nu r^\varrho \frac{\cos (\theta \varrho)}{\sin (\pi \varrho)}(1+o(1)), \qquad r \to + \infty.
\end{equation}
By taking $\varrho$ as in \eqref{Pn.bd.intro} and selecting a sufficiently large $\nu>0$ in \eqref{ak.intro}, we obtain a gauge function $F$ satisfying
\begin{equation}
\log |F(\la_n)| = \pi \nu |\la_n|^\varrho \frac{\cos (\varrho \arg \la_n )}{\sin (\pi \varrho)}(1+o(1)), \quad n \to \infty.
\end{equation}
This $F$ regularizes the growth of the spectral projections so that $B_z$ can be represented as a convergent series 
\begin{equation}
B_z(T) = \sum_{j=1}^\infty \frac1{z-\la_j} \frac1{F(\la_j)} P_j.
\end{equation}
Expressing the resolvent from \eqref{res.id.intro}, we obtain the desired bound on $\|(z-T)^{-1}\|$ using \eqref{F.asym.intro}, the boundedness of $B_z(T)$ and the convergence of the series on the r.h.s.~\eqref{res.id.intro}. The latter is justified by the convergence of the partial fraction decomposition of $1/F$
\begin{equation}\label{F.pfd.intro}
\frac{1}{F(w)} = \sum_{n=1}^{\infty} \frac{1}{F'(-a_n)} \frac{1}{w+a_n}, \quad w \in \C \setminus \{-a_n\}_{n \in \N},
\end{equation}
which holds if $\varrho <1/2$ and in which case
\begin{equation}
|F'(-a_n)| = \exp \big( \pi n \cot (\pi \varrho) + \BigO(\log n)\big), \quad n \to \infty,
\end{equation}
see Subsection~\ref{ssec:prod.pfd} and Appendix~\ref{app:pfc}. 

The identity \eqref{res.id.intro}--\eqref{B.F.intro} is inspired by a ``polynomial`` version 
\begin{equation}\label{Davies.ident.intro}
	(z-L)^{-1} \prod_{k=1}^m {(k+L)^{-1}} = \left( \prod_{k=1}^m \frac{1}{k+z} \right) (z-L)^{-1} + \sum_{k=1}^m \frac{c_k}{z+k} (k+L)^{-1},
\end{equation}
for certain $c_k \in \C$, $1\leq k \leq m < \infty$, written by E.~B.~Davies in \cite[p.~434]{Davies-2000-32}. He was motivated by the even imaginary oscillators in $L^2(\R)$
\begin{equation}\label{H.even.intro}
H = \Dt + \ii |x|^b, \quad b>0.
\end{equation}
Employing \eqref{Davies.ident.intro} and lower resolvent bounds of the type \eqref{L.res.lb.intro} found earlier in \cite{Davies-1999-200}, he showed that the spectral projections $P_n$ of $H$ are not ``tame'', i.e., for any $N>0$, 
\begin{equation}
\limsup_{n \to \infty} \frac{\|P_n\|}{n^N} = \infty.
\end{equation}
A numerical evidence indicated that $\|P_n\|$ grow exponentially which was for $b \in 2\N$ confirmed by a precise direct analysis by R.~Henry, namely
\begin{equation}
	c(b) \equiv	\lim_{n \to \infty} \frac{\log \|P_n\|}{n} >0,
\end{equation}
see \cite{Henry-2012-350,Henry-2014-4} where a more precise asymptotics is given. 

For $H$ in \eqref{H.even.intro}, Theorem~\ref{thm:T} and the existing lower bounds on $\|(z-H)^{-1}\|$ yield that for $b>2/3$ and a sufficiently small $\gamma>0$
\begin{equation}\label{H.even.limsup}
\limsup_{n \to \infty} \frac{\|P_n\|}{\exp(\gamma n)} = \infty,
\end{equation}
see Theorem~\ref{thm:H.even.Pn}. \emph{We expect that \eqref{H.even.limsup} holds with $\lim$ instead of $\limsup$} as in the special cases with $b \in 2\N$. The restriction $b>2/3$ appears natural in the sense that for $b=2/3$ the completeness of $\{P_n\}$ is non-trivial and it was established only recently by S.~Tumanov in \cite{Tumanov-2021-109,Tumanov-2022-319} for operators in $L^2(\R_+)$ with Dirichlet boundary condition at $0$, see Subsection~\ref{ssec:b23} for more details.

In Subsection~\ref{ssec:H.odd}, we show how Theorem~\ref{thm:T} can be employed also for odd imaginary oscillators
\begin{equation}\label{H.odd.intro}
H = \Dt + \ii x^{2b+1}, \quad b \in \N,
\end{equation}
with the prominent and extensively studied cubic case, see e.g.~\cite{Bender-1998-80,Dorey-2001-34,Shin-2002-229,Siegl-2012-86,Henry-2014-15,Giordanelli-2015-16}. This leads to the result analogous to \eqref{H.even.limsup} for a range of $b \in \N$, see Theorem~\ref{thm:H.odd.Pn}; the special case $b=1$ was analyzed directly in \cite{Henry-2014-15}, see Remark~\ref{rem:H.odd}. 

As a final application, we employ Theorem~\ref{thm:T} for conjugated oscillators in $L^2(\R)$
\begin{equation}\label{H.conj.intro}
	H = e^v \left(\Dt + |x|^b \right) e^{-v} = \left(-\ii \Do + \ii v'\right)^2 + |x|^b 
\end{equation}
where $b \geq 2$ and $v \in C^2(\R)$ is a non-negative odd function, increasing on $(0,\infty)$ and satisfying for some $ s \in (0,1)$
\begin{equation}\label{v.conj.intro}
	v(x) = \frac12 x^{\frac{2+b}{2} s}, \quad x>1.
\end{equation}
For these operators $H$, asymptotic behavior of the spectral projections norms was found in \cite{Mityagin-2017-272} and Theorem~\ref{thm:T} yields new resolvent norm bounds for $H$, see Subsection~\ref{ssec:MSV.op} for details.

This paper is organized as follows. In Section~\ref{sec:prelim} we recall useful results from operator theory, Schr\"odinger operators with complex potentials and infinite products and their partial fraction decompositions; a more detailed presentation of the latter is included in Appendix \ref{app:pfc}. In Section \ref{sec:L} we analyze $L$ from \eqref{L.intro} or \eqref{L.def}. In Subsection~\ref{ssec:L.basic} we summarize spectral properties of $L$, in particular the resolvent estimate \eqref{L.res.lb.even} is stated in Theorem~\ref{thm:L.res.lb}; a proof is given in Appendix~\ref{app:res.L.lb}. In Subsections \ref{ssec:Pn.to.res} -- \ref{ssec:L.res.bound} we prove a suitable version of the key identity \eqref{res.id.intro} for the resolvent of $L$ and obtain the resolvent bound in Theorem~\ref{thm:L} as well as the result on the spectral projections in Theorem~\ref{thm:L.Pn}. The consequences for the basis properties are discussed in Subsection~\ref{ssec:basis}. In Section~\ref{sec:T} we consider an abstract operator $T$ and establish \eqref{res.id.intro} and the related result on the resolvent norm (Theorem~\ref{thm:T}). The applications to imaginary and conjugated oscillators are in Subsections~\ref{ssec:H.even} -- \ref{ssec:MSV.op}. Finally, in Section~\ref{sec:misc} we collect various remarks on our results, their sharpness, critical cases and conjectures.

\section{Preliminaries}
\label{sec:prelim}

We summarize selected results from operator theory, analysis of Schr\"odinger operators and entire functions in a specialized form suitable for this paper; more general results can be found in the included references.

\subsection{Some facts from operator theory}
\label{ssec:ot}

Let $\cH$ be a separable Hilbert space and let $T$ be a densely defined operator with compact resolvent and $0 \in \rho(T)$. The spectrum of $T$ consists of isolated eigenvalues of finite algebraic multiplicity. Notice that the spectrum of such $T$ can be empty or a finite set, however, these cases are not in our focus here and we further assume that there exist infinitely many eigenvalues $\la_n$, $n \in \N$. We use the convention that the eigenvalues are repeated according to their algebraic multiplicities and ordered with non-decreasing moduli, i.e., $|\la_{n+1}|\geq|\la_n|$, $n \in \N$. 

The singular values of $T^{-1}$ are denoted by $\frs_n \equiv \frs_n(T^{-1})$, $n \in \N$, with the convention that they are repeated according to the multiplicity and ordered in a non-increasing way, i.e., $\frs_{n+1} \leq \frs_n$, $n \in \N$, see e.g.~\cite[Chap.~II]{Gohberg-1969}. $T^{-1}$ is in the Schatten classes $\cS_q$, $q \in (0,\infty)$, if
\begin{equation}
	\sum_{n=1}^\infty \frs_n(T^{-1})^q < \infty. 
\end{equation}
We recall that if there exists $\kappa >0$ such that 
\begin{equation}
\frs_n = \BigO(n^{-\kappa} ), \quad n \to \infty, 
\end{equation}
then
\begin{equation}\label{lak.sk.ineq}
 |\la_n^{-1}| = \BigO(n^{-\kappa}), \quad n \to \infty,
\end{equation}
see~\cite[Cor.~II.3.2]{Gohberg-1969}.

The Riesz projections of $T$ are denoted by 
\begin{equation}\label{Q.la.T.def}
Q_\la = \frac{1}{2 \pi \ii} \int_{|z-\la|=\delta_\la} (z-T)^{-1} \, \dd z, 
\quad
\delta_\la := \frac12 \dist (\la, \sigma(T)\setminus \{\la\}),
\quad
 \la \in \sigma(T),
\end{equation}
Notice that if $\Rank (Q_\la) = 1$, i.e.~ the eigenvalue $\la$ is simple, then $Q_\la = \langle \cdot, f_\la^* \rangle f_\la$ where $f_\la \in \Ker(T-\la)$ and $f_\la^* \in \Ker(T^*-\ov \la)$ are normalized so that  $\langle f_\la, f_\la^* \rangle =1$.

The numerical range of $T$ reads
\begin{equation}
\Num(T)	= \{  \langle Tf,f \rangle \, : \, f \in \Dom(T), \|f\|=1\}.
\end{equation}
For $T$ with compact resolvent, we have $\sigma(T) = \spp(T) \subset \Num(T)$ and 
\begin{equation}\label{res.Num.T}
\|(z-T)^{-1}\| \leq \frac{1}{\dist(z,\Num(T))}, \qquad z \in \left(\ov{\Num(T)} \right)^c;
\end{equation}
see e.g.~\cite[Chap.~V.3]{Kato-1966} or \cite[Chap.~26.1]{Markus-1988} for details. 

We recall the well-known result on the completeness of Riesz projections, see \cite[Cor.~XI.9.31]{DS2}. 

\begin{proposition}\label{prop:compl}
	Let $T$ be a densely defined closed operator in a separable Hilbert space $\cH$. Assume that $\rho(T) \neq \emptyset$ and for some $z_0 \in \rho(T)$
	\begin{equation}
		(z_0-T)^{-1} \in \cS_q
	\end{equation} 
	with $q \in (0,\infty)$. Let $\gamma_1, \dots, \gamma_j$ be non-overlapping differentiable arcs having a limit direction at infinity and suppose that no adjacent pair of arcs forms an angle as great as $\pi/q$ at infinity. Suppose that, for some $K \in \Z$, the resolvent of $T$ satisfies 
	\begin{equation}\label{res.pol.gr}
		\|(z-T)^{-1}\| = \BigO(|z|^K)
	\end{equation}
	along each $\gamma_i$, $i=1,\dots, j$, as $z \to \infty$. Then the system of Riesz projections $Q_{\lambda}$, $\la \in \sigma(T)$, is complete in $\cH$. 
\end{proposition}

Let the numerical range of $T$ be contained in a sector
\begin{equation}\label{Num.T.sector}
	\Num(T) \subset \{ z \in \C \, : \, |\arg z| \leq \vartheta \}, 
\end{equation}
with $\vartheta \leq \pi/2$. Then by \eqref{res.Num.T}, for any $\varphi \in (-\pi,-\vartheta) \cup (\vartheta,\pi]$
\begin{equation}
	\|(r e^{\ii \varphi} - T)^{-1}\| = \BigO(r^{-1}), \quad r \to + \infty.
\end{equation}
Thus if $T^{-1} \in \cS_q(\cH)$ for $0 < q < \pi/(2 \vartheta)$, Proposition~\ref{prop:compl} yields that the spectrum of $T$ comprises infinitely many eigenvalues and the system of Riesz projections $Q_\la$, $\la \in \sigma(T)$, defined in \eqref{Q.la.T.def} is complete in $\cH$.

Let the numerical range of $T$ be contained in a ``parabolic'' region
\begin{equation}\label{T.Num.par}
	\begin{aligned}
		\Num (T) \subset B_R(0)
				\cup \{ z \in \C \, : \, \Re z \geq 1, |\Im z| \leq h(\Re z)\}
	\end{aligned}	
\end{equation}
where $h : \R_+ \to \R_+$ is continuous and satisfies 
\begin{equation}\label{T.Num.h}
	\lim_{x \to + \infty} \frac{h(x)}{x} = 0.
\end{equation}	 
Then for any $\vartheta>0$, there exists $\alpha > 0$ such that 
\begin{equation}
	\Num (T+\alpha) \subset \{  z \in \C \, : \, |\arg z| \leq \vartheta\} 
\end{equation}
and also there exists $R_\vartheta>0$ such that 
\begin{equation}
\Num (T) \subset B_{R_\vartheta} \cup \{  z \in \C \, : \, |\arg z| \leq \vartheta\} 
\end{equation}
Thus by \eqref{res.Num.T}, for any $\varphi \in (-\pi,0) \cup (0,\pi]$,
\begin{equation}\label{T.Num.par.res}
\|(r e^{\ii \varphi} - T)^{-1}\| = \BigO(r^{-1}), \quad r \to + \infty.
\end{equation}
It follows from~Proposition~\ref{prop:compl} that if $T^{-1} \in \cS_q(\cH)$ for $q>0$, then the spectrum of $T$ comprises infinitely many eigenvalues and the system of Riesz projections $Q_\la$, $\la \in \sigma(T)$, is complete in $\cH$ .

\subsection{One-dimensional Schr\"odinger operators with complex potentials}
\label{ssec:1D.Schr}

We summarize the special case of \cite[Thm.~2.5, Prop.~2.6]{Krejcirik-2017-221}; cf.~also \cite{Almog-2015-40,Boegli-2017-42,Semoradova-2022-54, Tumanov-2024-215}. 
\begin{proposition}\label{prop:1D.Schr}
Let $V: \R \to \C$ be such that $V = V_0 + V_1$ where $V_0 \in C^1(\R)$ satisfies $\Re V_0 \geq 0$,  
\begin{equation}\label{V.unbdd}
	|V_{0}(x)| \to \infty \quad \text{and} \quad |V_0'(x)| = o(|V_0(x)|^\frac 32), \quad |x| \to \infty,
\end{equation}
and $V_1 \in L^\infty_{\rm loc}(\R)$ satisfies 
\begin{equation}\label{V1.oV0}
V_1(x) = o(|V_0(x)|), \quad |x| \to \infty.
\end{equation}
Let $H$ be a Schr\"odinger operator in $L^2(\R)$  
\begin{equation}
H =\Dt + V, \quad \Dom(H) = \{ f \in H^2(\R) \, : \, Vf \in L^2(\R) \}.
\end{equation}
Then $\rho(H) \neq \emptyset$, the resolvent of $H$ is compact and for every $\eps \in (0,1)$, there exists $C_\eps>0$ such that
\begin{equation}\label{H.graph.sep}
	\| H f \|^2 \geq (1-\eps) (\|f''\|^2 + \|V_0 f\|^2 ) - C_\eps \|f\|^2, \quad f \in \Dom(H).
\end{equation}
\end{proposition}
We remark that under the assumption \eqref{V1.oV0} on $V$, the condition  $Vf \in L^2(\R)$ in the characterization of $\Dom(H)$ is can be changed to $V_0f \in L^2(\R)$.

In particular, let $W: \R \to \R$ be such that $W = W_0 + W_1$ with $W_0 \in C^1(\R;\R)$, $W_0 \geq 0$,  $W_1 \in L^\infty_{\rm loc}(\R;\R)$ and
\begin{equation}
W_0(x) \to + \infty, \quad |W_0'(x)| = o(|W_0(x)|^\frac 32) \quad \text{and} \quad	W_1(x) = o(|W_0(x)|)
\end{equation}
as $|x| \to \infty$. Then, as a consequence of Proposition~\ref{prop:1D.Schr}, the operator
\begin{equation}
G =\Dt + W, \quad \Dom(G) = \{ f \in H^2(\R) \, : \, W f \in L^2(\R) \} 
\end{equation}
is self-adjoint with compact resolvent and the graph norm separation 
\begin{equation}\label{G.graph.sep}
	\| G f \|^2 \geq (1-\eps) (\|f''\|^2 + \|W_0 f\|^2 ) - C_\eps \|f\|^2, \quad f \in \Dom(G)
\end{equation}
holds. Hence, if $W_0$ and $|V_0|$ are of the same size at infinity, namely,
\begin{equation}\label{V0.W0.inf}
	\lim_{|x|\to \infty} \frac{|V_0(x)|}{W_0(x)} =1,
\end{equation} 
\eqref{H.graph.sep}, \eqref{G.graph.sep} and \eqref{V0.W0.inf} allow for comparing the singular values of the resolvents of $H$ and $G$, namely, 
there exist $\alpha, c_1, c_2 >0$ such that
\begin{equation}\label{sn.GH}
c_1 \frs_n((G+\alpha)^{-1}) 	\leq \frs_n((H+\alpha)^{-1}) \leq c_2 \frs_n((G+\alpha)^{-1}), \quad n \in \N;
\end{equation}
cf.~also \cite[Thm.~1.1, 1.2]{Almog-2015-40} and \cite[Chap.~XII]{EE}. 

The important case for our analysis is the power-potential, namely the self-adjoint anharmonic oscillator
\begin{equation}\label{G0.def}
	G_0 = \Dt + |x|^l, \qquad \Dom(G_0) = H^2(\R) \cap \Dom(|x|^l), \quad l \in (0,\infty).
\end{equation}
All eigenvalues of $G_0$ are simple, positive and we denote them by $\mu_n$, $n \in \N$; we order $\mu_n$ in an increasing way. It is known that they obey the asymptotics as $n \to \infty$
\begin{equation}\label{mu.n.asym}
	\mu_n = (d n)^\kappa(1+o(1)), 
	\quad 
	\mu_{n+1} - \mu_n = \kappa d (dn)^{\kappa-1}(1+o(1)),
\end{equation}
where
\begin{equation}\label{kap.G_0.even}
	d = \frac{\sqrt \pi \ \Gamma(\frac 32 + \frac 1 l)}{\Gamma(1 + \frac 1 l)}, \qquad 	\kappa = \frac{2l}{l+2},
\end{equation}
see \cite[Chap.~VII]{Titchmarsh-1962-book1}, \cite[Thm.~XIII.82]{Reed4} or \cite[Sec. ~6]{Mityagin-2019-139}.

\subsection{Infinite products and partial fraction decompositions}
\label{ssec:prod.pfd}

We summarize some known facts on infinite products and related partial fraction decompositions. 

Let $F: \C \to \C$ be an infinite product 
\begin{equation}\label{F.def}
	F(w) = \prod_{k=1}^\infty \left( 1 + \frac{w}{a_k} \right), \qquad a_{k+1} \geq a_k>0, \ k \in \N,
\end{equation}
where 
\begin{equation}\label{n.r.asm}
\exists \nu >0 \ \exists \varrho \in (0,1): \ \	n(r) := \# \{ k \, : \, a_k \leq r\} = \nu r^\varrho (1+o(1)), \quad r \to + \infty.
\end{equation}
The following properties of $F$ are known, see Appendix~\ref{app:pfc} for details, references and proofs.
\begin{proposition}\label{prop:F.Titch}
Let $F$ be as in \eqref{F.def} with \eqref{n.r.asm}. Then the infinite product \eqref{F.def} is absolutely convergent, $F$ is an entire function and for any $\theta \in (-\pi,\pi)$
\begin{equation}\label{F.asym.Titch}
\log|F(r e^{\ii \theta})| = \pi \nu r^\varrho \frac{\cos (\theta \varrho)}{\sin (\pi \varrho)}(1+o(1)), \qquad r \to + \infty;
\end{equation}
moreover, the remainder estimate is uniform for $\theta \in [-\pi + \eps, \pi - \eps]$ with any fixed $\eps>0$.
\end{proposition}

An essential ingredient in our proofs below is the partial fraction decomposition of $1/F$. We formulate claims when zeros $\{-a_k\}_{k \in \N}$ are determined by
\begin{equation}\label{an.special}
a_k := \left( \frac{k}{\nu} \right)^{\frac 1 \varrho}, \qquad k \in \N,
\end{equation}
where $\nu>0$ and $\varrho \in (0,1)$. (Notice that then \eqref{n.r.asm} is satisfied.) 

\begin{proposition}\label{prop:F.pfc.new}
Let $F$ be as in \eqref{F.def} where $a_k$, $k \in \N$, satisfy \eqref{an.special} with $\nu>0$ and $\varrho \in (0,1)$.
Then 
\begin{equation}\label{F'.an.asym}
|F'(-a_k)| = \exp \big( \pi k \cot (\pi \varrho) + \BigO(\log k)\big), \qquad k \to \infty.
\end{equation}
If in addition $\varrho \in (0,1/2)$, then 
\begin{equation}\label{pfc.F.b>2.new}
\frac{1}{F(w)} = \sum_{n=1}^\infty \frac{1}{F'(-a_n)} \frac{1}{w+a_n},
\end{equation}
where the series converges uniformly on compacts subsets of $\C \setminus \{-a_k\}_{k=1}^\infty$.

Hence, for all $w \neq z$ and $w \neq -a_n$, $n \in \N$, 
\begin{equation}\label{pf.decomp}
	\frac{1}{(z-w) F(w)} = \frac{1}{F(z)(z-w)} + \sum_{n=1}^\infty \frac{1}{(z+a_n) F'(-a_n)(w+a_n)}.
\end{equation}
\end{proposition}
For $\varrho =1/2$, the terms $|F'(-a_k)|$ are not exponentially growing and the partial fraction decomposition like \eqref{pfc.F.b>2.new} is not available. Nonetheless, Euler's product formula for $\sin$ leads to the decomposition \eqref{pfc.F.12} below.

\begin{proposition}\label{prop:F.pfc.new.12}
Let $F$ be as in \eqref{F.def} where $a_k$, $k \in \N$, satisfy \eqref{an.special} with $\nu>0$ and $\varrho = 1/2$, $k \in \N$.
Then 
\begin{equation}\label{pfc.F.12}
	\frac{1}{F(w)} 
	= 1 + 2 \sum_{n=1}^\infty (-1)^n \frac{w}{w+a_n}, 
\end{equation}
where the series converges uniformly on compacts subsets of $\C \setminus \{-a_k\}_{k=1}^\infty$.

Hence, for all $w \neq z$ and $w \neq -a_n$, $n \in \N$, 
\begin{equation}\label{pf.decomp.12}
	\frac{1}{(z-w) F(w)} = \frac{1}{F(z)(z-w)} +  2 \sum_{n=1}^\infty (-1)^{n+1}\frac{a_n}{(z+a_n)(w+a_n)}.
\end{equation}
\end{proposition}

\begin{remark}\label{rem:ak.modification}
The results above are stated for $a_k$, $k \in \N$, satisfying \eqref{an.special}. In fact, they are based on the analysis of 
\begin{equation}
\Phi(w) = \prod_{k=1}^\infty \left(1 + \frac{w}{k^\frac1 \varrho}\right),	
\end{equation}
the special case with $\nu=1$. The claims for a general $\nu>0$ can be explained by rescaling
\begin{equation}
F(w) = \Phi(t w), \quad t = \nu^\frac1\varrho;
\end{equation}
cf.~Appendix~\ref{ssec:App.A.3}.

One can also perturb $a_k$ in \eqref{an.special} by a bounded sequence, i.e.,~consider
\begin{equation}
a_k := \left( \frac{k}{\nu} \right)^{\frac 1 \varrho} + b_k, \qquad k \in \N, \qquad b_k = \BigO(1), \quad k\to \infty.
\end{equation}
It is showed in \cite[Lemma~6]{Sedletskii-1983-34} that the asymptotics \eqref{F'.an.asym} remains valid and hence also \eqref{pfc.F.b>2.new} and \eqref{pf.decomp} hold. 
\end{remark}

\section{Spectral projections of complex perturbations of anharmonic oscillators}
\label{sec:L}

The first main object of our analysis is an operator in $L^2(\R)$
\begin{equation}\label{L.def}
L = \Dt + x^{2a} + V_1(x), \quad \Dom(L) = \{ f \in H^2(\R) \, : \, x^{2a} f(x) \in L^2(\R) \}
\end{equation}
with $a \in \N$ and a polynomial potential $V_1$ with
\begin{equation}\label{V.def}
\deg (\Re V_1) \leq 2a-1, \quad \deg (\Im V_1) = b, \quad 0 \leq a-1 < b <2a, 
\end{equation}
and 
\begin{equation}\label{V.deg}
\Im V_1(x) = c_b x^b + \dots, \quad  c_b >0.
\end{equation}
We will further assume that $\Re V_1 \geq 1$, which can be always achieved by a shift, i.e.,~by considering $L+\alpha$ instead of $L$ with a sufficiently large $\alpha>0$. (Clearly, adding of the scalar multiple of the identity operator preserves all spectral properties of $L$.)

\subsection{Spectral properties of $L$}
\label{ssec:L.basic}

Basic facts about $L$ are summarized in the following.
\begin{proposition}\label{prop:L.basic}
Let $L$ be as in \eqref{L.def} -- \eqref{V.deg}. Then
	\begin{enumerate}[\upshape (i)]
		\item for any $\eps>0$, there exists $R_\eps > 0$ such that the numerical range of $L$ satisfies 
		\begin{equation}\label{Num.L.odd}
		\begin{aligned}				
		\Num (L) & \subset \{ z \in \C \, : \, \Re z \geq 1, |z| \leq R_\eps \} 
		\\
		& \qquad \cup \left\{   z \in \C  \,:\, \Re z \geq 1, |\Im z| \leq (c_b+\eps) (\Re z)^{\frac{b}{2a}} \right\} \quad \text{if $b$ is odd}
		\end{aligned}
		\end{equation}
		and 
		\begin{equation}\label{Num.L.even}
		\begin{aligned}				
		\Num (L) & \subset \{ z \in \C \, : \, \Re z \geq 1, |z| \leq R_\eps \} 
		\\
		& \qquad \cup \left\{   z \in \C  \,:\, \Re z \geq 1, 0 \leq \Im z \leq (c_b+\eps) (\Re z)^{\frac{b}{2a}} \right\} \quad \text{if $b$ is even};
		\end{aligned}
		\end{equation}
		\item $L$ is m-accretive with compact resolvent, moreover, 
		\begin{equation}\label{sk.L}
			\frs_k(L^{-1}) = \BigO(k^{- \kappa}), \quad k \to \infty,
		\end{equation}
		where
		\begin{equation}\label{kap.L.def}
		\kappa=	\frac{2a}{a+1};
		\end{equation}
%		thus in particular $L^{-1} \in \cS_p$ for any $p> (a+1)/(2a)$;
		\item \label{prop:L.basic.complete}
		the spectrum of $L$ comprises infinitely many eigenvalues and the system of Riesz projections $Q_\la$, $\la \in \sigma(L)$, as in \eqref{Q.la.T.def} with $T=L$, is complete in $L^2(\R)$.
	\end{enumerate}
\end{proposition}
\begin{proof}[Sketch of proof]
\begin{enumerate}[\upshape (i), wide]
\item The claims \eqref{Num.L.odd}, \eqref{Num.L.even} follow by straightforward estimates since the potential $V(x) = x^{2a} + V_1(x)$, $x \in \R$, satisfies $\Re V(x) = x^{2a} + \BigO(x^{2a-1})$ and $\Im V(x)= c_b x^b + \BigO(x^{b-1})$ as $|x| \to \infty$. 
\item The conditions on $V$ in Subsection~\ref{ssec:1D.Schr} with $V_0(x) = x^{2a}$ are satisfied, so the resolvent set of $L$ is non-empty and the resolvent is compact (the accretivity of $L$ follows from the enclosure of $\Num(L)$ above). The behavior of the singular values in \eqref{sk.L} and \eqref{kap.L.def} is obtained by \eqref{sn.GH} where $H=L$ and $G=G_0$ (with $l=2a$) are used, see \eqref{G0.def} and  \eqref{mu.n.asym}.
\item The numerical range of $L$ is contained in a ``parabolic'' region on the r.h.s.~of \eqref{T.Num.par} with $h(t)=(c_b+\eps)t^{\frac b{2a}}$, $t>0$. Moreover, from \eqref{sk.L}, we have $L^{-1} \in \cS_q$ for any $q > 1/\kappa$, so the claims on eigenvalues and completeness follow by Proposition~\ref{prop:compl} and remarks below it.
\qedhere
\end{enumerate}
\end{proof}

Starting from the results in \cite{Sibuya-1975}, K.~C.~Shin proved the following more precise spectral results for $L$, which we formulate in the form convenient for us; cf.~also \cite{Shin-2002-229,Shin-2005-38}.
\begin{proposition}[{\cite[Thm.~1.1, 1.3, Cor.~2.1, 2.2]{Shin-2010-6}}]\label{prop:Shin}
	Let $L$ be as in  \eqref{L.def} -- \eqref{V.deg}, let $\kappa$ be as in \eqref{kap.L.def} and let $\{\la_n\}_{n \in \N}$ be the eigenvalues of $L$ (repeated with the algebraic multiplicity and ordered with a non-decreasing modulus). Then 
	\begin{enumerate}[\upshape (i)]
		\item there exists $N_0 \in \N$ such that for all $n \geq N_0$ the eigenvalues $\la_n$
%%		
%		\begin{equation}
%			\la_n = t_n + \ii s_n
%		\end{equation}
%%
		are simple and  
		\begin{equation}
			|\la_n| < |\la_{n+1}|;
		\end{equation}
		\item the real part of $\la_n$ satisfies
		\begin{equation}\label{tn.L.asym}
			\Re \la_n = (d n)^\kappa(1 + \BigO(n^{-\frac{1}{1+a}})), \quad n \to \infty,
		\end{equation}
		where
		\begin{equation}\label{d.L.def}
		d = \frac{\pi}{B\left(\tfrac12,1+ \tfrac1{2a}\right)} >0;
		\end{equation}
	\item the ``gaps'' satisfy
	\begin{equation}\label{la.L.gaps}
		\la_{n+1} - \la_n = \kappa d (dn)^{\kappa-1} \left(1 + \BigO(n^{-\frac{1}{1+a}})\right), \quad n \to \infty;
	\end{equation}
	\item the imaginary part of $\la_n$ satisfies
	\begin{equation}
		\Im \la_n = (d' n)^{\kappa'} \left(1 + \BigO(n^{-\frac{1}{1+a}})\right), \quad n \to \infty,
	\end{equation}	
	where
	\begin{equation}
	d' \neq 0, \quad \kappa > \kappa' = \kappa - \frac{j}{1+a} \quad \text{for some $j$ with  } 1 \leq j \leq 2a,
	\end{equation}
	or
	\begin{equation}
	\Im \la_n = \BigO(n^{-\frac{1}{1+a}}), \quad n \to \infty;
	\end{equation}
	
	\item if $V_0$ is even and $V_1$ is odd, then the eigenvalues $\{\la_n\}_{n\in \N}$ are real with possibly finitely many exceptions.
	\end{enumerate}
\end{proposition}

%\begin{remark}\label{rem:L.shift}
We remark that the case $a=1$ and $b=1$ is more explicit, namely, for 
\begin{equation}\label{L.shift.HO.def}
	L = \Dt + x^2 + \alpha_1 x + \alpha_0, \quad \alpha_1, \alpha_0 \in \C, 
\end{equation}
all eigenvalues are simple and satisfy
\begin{equation}\label{la.shifted.HO}
	\la_n = 2n - 1 + \alpha_0 - \frac{\alpha_1^2}{4}, \qquad n \in \N;
\end{equation}
cf.~\cite{Mityagin-2017-272,Krejcirik-2015-56,Mityagin-2021-22}.  
%\end{remark}
%

The resolvent norm of $L$ satisfies the following lower bounds. 
\begin{theorem}\label{thm:L.res.lb}
Let $L$ be as in \eqref{L.def} -- \eqref{V.deg} and let
\begin{equation}\label{rho.L.lb}
	\tau: = \frac{b+1}{2a} - \frac 12 = \frac{b-(a-1)}{2a}.
\end{equation}
(Notice that due to \eqref{V.def}, we have $\tau \in (0,1/2)$.)
Then the following holds.
\begin{enumerate}[\upshape (i)]
	\item If $b$ is odd, then for every $\eps \in (0,c_b)$  there exist $\eta , R >0$ such that for all $z \in \rho(L)$ satisfying
	\begin{equation}
		\Re z >R \quad \text{and} \quad |\Im z| \leq (c_b-\eps) (\Re z)^\frac{b}{2a},
	\end{equation}
	we have 
	\begin{equation}\label{L.res.lb.odd}
		\|(z-L)^{-1}\| \geq \exp(\eta (\Re z)^{\tau}).
	\end{equation}

	\item If $b$ is even, then for every $\eps \in (0,c_b)$ and $\omega \in (0, \omega_0)$ with
	\begin{equation}\label{var.kap.def}
		\omega_0 := \frac{b}{2a} - \frac{b}{2(b+1)} =\frac b2 \left(\frac1a - \frac1{b+1}\right)= \frac{b}{b+1} \tau, 
	\end{equation}	
	 there exist $\eta, R >0$ such that for all $z \in \rho(L)$ satisfying 
	\begin{equation}\label{pseudo.region.even}
		\Re z >R, \quad \Im z \geq (\Re z)^{\frac{b}{2(b+1)} + \omega} \quad \text{and} \quad \Im z \leq (c_b-\eps) (\Re z)^\frac{b}{2a},
	\end{equation}
	we have
	\begin{equation}\label{L.res.lb.even}
		\|(z-L)^{-1}\| \geq \exp\left(\eta \frac{(\Im z)^{\frac{b+1}b}}{(\Re z)^\frac12} \right) \geq 
		\exp\left(\eta (\Re z)^{\omega \frac{b+1}b} \right).
	\end{equation}
\end{enumerate}
\end{theorem}

We remark that $\omega_0 >0$ since by \eqref{V.def}
\begin{equation}
	\frac{1}{a} > \frac{1}{b+1}.
\end{equation}
\begin{proof}
The proof is given in Appendix~\ref{app:res.L.lb} and it is based on construction of pseudomodes, see in particular \eqref{pseudo.est} and \eqref{L.res.lb.proof}.
\end{proof}

\subsection{From spectral projections to the resolvent}
\label{ssec:Pn.to.res}

The goal of this section is to relate the upper bounds on the norms of spectral projections and on the resolvent norm. Relying on Proposition~\ref{prop:Shin}, we fix $R_0>0$ such that $|\la_{N_0}|<R_0<|\la_{N_0+1}|$ and introduce
\begin{equation}\label{S0.L.def}
	S_0 = \frac{1}{2 \pi \ii} \int_{|z|=R_0} (z-L)^{-1} \, \dd z,
\end{equation}
and for $n > N_0$ 
\begin{equation}\label{Pn.L.def}
	P_n =  Q_{\la_n} = \frac{1}{2 \pi \ii} \int_{|z-\la_n|=\delta_{\la_n}} (z-L)^{-1} \, \dd z,
\end{equation}
where 
\begin{equation}\label{del.n.L.def}
\delta_{\la} = \frac12 \dist (\la, \sigma(L)\setminus \{\la\}).
% = \frac 12 \min_{j \neq n} \{|\la_j-\la_n|\}. 	
\end{equation}
Notice that 
\begin{equation}
	\Rank (S_0) =N_0, \quad \Rank (P_n) =1, \quad n > N_0,
\end{equation}
and that the system of projections $\{S_0\} \cup \{P_n\}_{n>N_0}$ is disjoint and complete; see Proposition~\ref{prop:L.basic} and \eqref{Pn.L.def}. It follows that the linear space
\begin{equation}\label{DL.def}
	\cD_L := \left\{ h \in L^2(\R) \,: \, \exists N(h) > N_0 \ \forall k > N(h) \  P_k h = 0 \right\} 
\end{equation}
is dense in $L^2(\R)$. Moreover, the resolvent of $L$ acts on $\cD_L$ as
\begin{equation}\label{L.res.DL}
(z-L)^{-1} h = (z-L S_0)^{-1} S_0 h + \sum_{j = N_0+1}^{N(h)} \frac{1}{z-\la_j} P_j h, \qquad h \in \cD_L, \ z \in \rho(L).
\end{equation}

Under proper assumptions or estimates, \eqref{L.res.DL} could lead to identities for all $h \in L^2(\R)$ or to the operator norm inequalities. The latter is in particular important for us when, being inspired by \cite{Davies-2000-32}, we use and develop a resolvent identity 
\begin{equation}\label{Davies.ident}
(z-L)^{-1} \prod_{k=1}^m {(k+L)^{-1}} = \left( \prod_{k=1}^m \frac{1}{k+z} \right) (z-L)^{-1} + \sum_{k=1}^m \frac{c_k}{z+k} (k+L)^{-1},
\end{equation}
for certain constants $c_k \in \C$, $1 \leq k \leq m < \infty$, see~\cite[p.~434]{Davies-2000-32}.

We treat \eqref{Davies.ident} as an operator version of the partial fraction decomposition
\begin{equation}\label{Davies.pfc}
\begin{aligned}
\frac{1}{z-w} \frac{1}{\Phi(w)} &= \frac{1}{z-w} \frac{1}{\Phi(z)} + \sum_{k=1}^m \frac{c_k}{z+a_k} \frac{1}{a_k+w}, 
\\	
c_k & = \frac{1}{\Phi'(-a_k)}, \qquad 1 \leq k \leq m < \infty,
\end{aligned}
\end{equation}
where
\begin{equation}\label{F.polyn}
	\Phi(w) := \prod_{k=1}^m (w+a_k)
\end{equation}
is a polynomial of $\deg \Phi=m$ and with simple negative zeros $\{ -a_k\}_{k=1}^m$; in the case of \eqref{Davies.ident} we have $a_k = k$, $1 \leq k \leq m$.
This point of view gives us the tools to deal with ``wild'' systems of projections $\{P_n\}_{n \in \N}$, i.e.,~with the case when the norms of $\|P_n\|$ could grow faster than any polynomial; cf.~\cite{Davies-2000-32} for the results on the ``tame'' system, i.e., when the growth of the norms $\|P_n\|$ is polynomial.

\subsection{Regularization of the resolvent expansion}
\label{ssec:L.reg.res}

To be able to regularize exponentially growing spectral projections norms, we aim at establishing a version of \eqref{Davies.ident} for infinite products instead of the polynomial in \eqref{F.polyn} above. We select a gauge function $F: \C \to \C$ in the form \eqref{F.def} where $a_k$, $k \in \N$, satisfy \eqref{an.special} with $\nu>0$ and $\varrho \in (0,1/2)$. 

As the first step in the regularization (like \eqref{Davies.ident}), we introduce an operator $B_z$ and estimate its norm.

\begin{lemma}\label{lem:B.L.bound} 
Let $L$ be as in \eqref{L.def} -- \eqref{V.deg} and $\kappa$, $d$ as in \eqref{kap.L.def} and \eqref{d.L.def}. Let $P_n$, $n >N_0$, be as in \eqref{Pn.L.def} and suppose that there exist $M, \gamma, \sigma >0$ such that
\begin{equation}\label{Pn.L.asm}
 \|P_n\| \leq M \exp \left(\gamma (dn)^\sigma\right), \qquad n > N_0,
\end{equation}
and also such that (with $\kappa$ from \eqref{kap.L.def})
\begin{equation}\label{rho.L.def}
 \varrho := \frac{\sigma}{\kappa} < \frac 12.
\end{equation}
For $\delta>0$, let $F$ be as in \eqref{F.def} where $a_k$, $k \in \N$, satisfy \eqref{an.special} with $\varrho$ from \eqref{rho.L.def} and 
\begin{equation}\label{nu.L.set}
\nu = \gamma \frac{\sin\left(\pi \varrho \right)}{\pi} (1+\delta).
\end{equation}	
For each $z \in \rho(L)$, consider the densely defined operator (recall \eqref{DL.def})
	\begin{equation}\label{B.L.def}
		\begin{aligned}
			B_z h := 
			\sum_{j=N_0+1}^{N(h)} \frac{1}{z-\la_j} \frac{1}{F(\la_j)} P_{j} h,  
			\qquad
			h \in \Dom(B_z) := \cD_L.
		\end{aligned}
	\end{equation}
Then $B_z$ admits a unique bounded extension on $L^2(\R)$ and it satisfies
\begin{equation}\label{B.L.norm}
\|B_z \| \leq \frac{M_1}{\dist(z, \{\la_n\}_{n > N_0} )}, \quad z \in \rho(L),
\end{equation}
where 
\begin{equation}\label{M1.L.def}
	M_1:=	\sum_{j=N_0+1}^\infty \frac{\|P_j\|}{|F(\la_j)|} < \infty.
\end{equation}
\end{lemma}
\begin{proof}
It follows from Proposition~\ref{prop:Shin} and \eqref{Num.L.odd} or \eqref{Num.L.even} that for an arbitrarily small $\vartheta > 0$ all eigenvalues $\la_j$ with a sufficiently large $j$ lie in a sector around the positive real axis with a half-angle $\vartheta$. Thus, using also \eqref{tn.L.asym}, there exists  $j_\delta \in \N$ such that for all $j > j_\delta$
\begin{equation}\label{la.j.est}
|\la_j| \geq \frac{1}{\left(1 + \frac \delta 2\right)^\frac{1}{3 \varrho}} (d j)^\kappa, \qquad \cos(\varrho \arg \la_j) \geq  \frac{1}{\left(1 + \frac \delta 2\right)^\frac{1}{3}}.
\end{equation}
Moreover, by Proposition~\ref{prop:F.Titch}, $j_\delta$ can be selected such that for all $j > j_\delta$
\begin{equation}
|F(\la_j)| \geq \exp \left(
\frac{1}{\left(1+ \frac \delta 2\right)^\frac13} \pi \nu |\la_j|^\varrho \frac{\cos\left(\varrho \arg \la_j  \right)}{\sin ( \pi \varrho )}
\right).
\end{equation}
Next, by \eqref{la.j.est}, \eqref{nu.L.set} and \eqref{rho.L.def}, we obtain for all $j > j_\delta$
\begin{equation}
\begin{aligned}
|F(\la_j)| 	& \geq \exp \left(
	\frac{1}{1+ \frac \delta 2} \nu (dj)^{\kappa\varrho} \frac{\pi}{\sin ( \pi \varrho )}
	\right)
	 = 
	\exp \left(
	\frac{1+\delta}{1+ \frac \delta 2} \gamma (dj)^{\sigma} 
	\right).	
	\end{aligned}
\end{equation}
Employing the assumption \eqref{Pn.L.asm}, we have for all $j > j_\delta$ that
\begin{equation}\label{P.F.L.est}
	\frac{\|P_{j}\|}{|F(\la_j)|} \leq M
	\exp \left(-  \frac{\delta}{2+\delta} \gamma (dj)^\sigma 
	\right),
\end{equation}
thus the series in \eqref{M1.L.def} converges. Returning to \eqref{B.L.def}, we obtain that
\begin{equation}
\|B_z h\| \leq \frac{M_1}{\dist(z, \{\la_n\}_{n > N_0} )} \|h\|, \quad z \in \rho(L).
\end{equation}
Since $\Dom(B_z)=\cD_L$ is dense in $L^2(\R)$, the operator $B_z$ indeed admits a unique bounded extension and it satisfies \eqref{B.L.norm}.
\end{proof}

\subsection{Resolvent estimate for $L$ and its corollaries}
\label{ssec:L.res.bound}
In the next step we employ the partial fraction decomposition from Proposition~\ref{prop:F.pfc.new} and obtain an upper bound on the resolvent norm of $L$.

\begin{theorem}\label{thm:L}
Let $L$ be as in \eqref{L.def} -- \eqref{V.deg} and $\kappa$, $d$ as in \eqref{kap.L.def} and \eqref{d.L.def}.
Let $P_n$, $n >N_0$, be as in \eqref{Pn.L.def} and assume that \eqref{Pn.L.asm} and \eqref{rho.L.def} are satisfied.
Then for any $\delta>0$, there exist $C \equiv C(\delta)>0$ and $r \equiv r(\delta)>0$ such that for all $z \in \rho(L)$ satisfying $|z| \geq r$ and $|\arg z| \leq  \pi/2$
	\begin{equation}\label{res.norm.L}
		\begin{aligned}
			\|(z-L)^{-1} \| 
			& \leq 
			 \frac{C}{\dist(z, \sigma(L))} 
			 \exp \left( 
			 (1+\delta) \gamma \cos (\varrho \arg z) |z|^\varrho  
			 \right).  
		\end{aligned}
	\end{equation}
\end{theorem}

\begin{proof}
Let $F$ be as in \eqref{F.def}, where $a_k$, $k \in \N$, satisfy \eqref{an.special} with $\varrho$ from \eqref{rho.L.def} and 
\begin{equation}\label{nu.L.wt}
\nu = \gamma \frac{\sin\left(\pi \varrho \right)}{\pi} (1+\wt \delta)
\end{equation}	
for $\wt \delta>0$. (Recall that then and \eqref{n.r.asm} holds with these $\nu$ and $\varrho$.) We start with \eqref{B.L.def} and use the partial fraction decomposition \eqref{pf.decomp} and the estimate of the norm of $B_z$ in Lemma~\ref{lem:B.L.bound}. 
Namely, we have from \eqref{pf.decomp} that for all $h \in \cD_L$ and $z \in \rho(L)$,
	\begin{equation}\label{Bh.L.PFD}
		\begin{aligned}
			B_z h & = 	\sum_{j=N_0+1}^{N(h)} \frac{1}{(z-\la_j)F(\la_j)} P_j h 
			\\
			& =  \frac1{F(z)}\sum_{j=N_0+1}^{N(h)} 
			\frac{1}{z-\la_j} P_j h 
%			\\  & \qquad 
			+ \sum_{n=1}^\infty \frac{1}{(z+a_n) F'(-a_n)}  \sum_{j=N_0+1}^{N(h)}  \frac{1}{\la_j+a_n} 
			P_j h
			\\	 
			& = 
			\frac1{F(z)} (z-L)^{-1} (h-S_0h)
%			\\
%			& \qquad 
			+  \sum_{n=1}^\infty \frac{1}{(z+a_n) F'(-a_n)}  (a_n+L)^{-1}(h-S_0h).
		\end{aligned}
	\end{equation}
	Since $L$ is m-accretive, by \eqref{res.Num.T}, we get
	\begin{equation}\label{res.L.accr}
	\|(a_n+L)^{-1}\| \leq \frac{1}{a_n},
	\end{equation}
	thus \eqref{an.special} and \eqref{F'.an.asym} yield
	\begin{equation}\label{M2.L.def}
		M_2:=\sum_{n=1}^\infty \frac{\|(a_n+L)^{-1}\|}{|F'(-a_n)|} < \infty. 
	\end{equation}
	Since $\sigma(L S_0) \subset B_{R_0}(0)$, we also have
	\begin{equation}\label{res.S0.est.L}
		\| (z-L S_0)^{-1} S_0\| \leq \frac{\|S_0\|}{|z| - \|L S_0\|} \leq \frac{2 \|S_0\|}{|z|} , \quad |z| \geq  2 \|LS_0\|. 
	\end{equation}
	
	Rearranging \eqref{Bh.L.PFD}, we get
	\begin{equation}\label{B.res.rearr}
		\begin{aligned}
		\frac{1}{F(z)}	(z-L)^{-1} h &= 
			B_z h 
			- \sum_{n=1}^\infty \frac{1}{(z+a_n) F'(-a_n)}  (a_n+L)^{-1}(h-S_0h)
			\\
			& \quad 
			+ \frac{1}{F(z)} (z-L S_0)^{-1} S_0 h.
		\end{aligned}
	\end{equation}
	Hence, since $\cD_L$ is dense in $L^2(\R)$, by \eqref{B.L.norm}, \eqref{M1.L.def}, \eqref{M2.L.def} and \eqref{res.S0.est.L}, we obtain for all $z \in \rho(L)$ satisfying $|z| \geq 2 \|LS_0\|$ and  $|\arg z| \leq \pi/2$ that
	\begin{equation}\label{L.res.F.est1}
		\begin{aligned}
			\|(z-L)^{-1} \| & \leq |F(z)| \left(
			\frac{M_1}{\dist(z, \{\la_n\}_{n>N_0})}
			+
			\frac{M_2(1+\|S_0\|)}{\dist(z, \{-a_n\}_{n \in \N})}
			\right)
			\\
			& \quad 
			+ \|(z-L S_0)^{-1} S_0\|
			\\
			& \leq 
			|F(z)| \left(
			\frac{M_1}{\dist(z, \sigma(L))}
			+
			\frac{M_2(1+\|S_0\|)}{|z|}
			\right)
			+ \frac{2\|S_0\|}{|z|}.
		\end{aligned}
	\end{equation}
It follows from Proposition~\ref{prop:F.Titch} with \eqref{an.special}, \eqref{rho.L.def} and \eqref{nu.L.set} that there exists $r' \equiv r'(\wt \delta)>0$ such that for all $z \in \C$ satisfying $|z|> r'$ and  $|\arg z| \leq \pi/2$
\begin{equation}
|F(z)|  
\leq \exp \left( \pi \nu |z|^\varrho \frac{\cos (\varrho \arg z)}{\sin (\pi \varrho)} (1+\wt\delta)\right)
= 
\exp \left( (1+\wt \delta)^2 \gamma \cos (\varrho \arg z) |z|^\varrho  \right).
\end{equation}
Hence, returning to \eqref{L.res.F.est1}, we obtain that there exist $\wt r \equiv \wt r(\wt \delta) \geq r'$ and $\wt C \equiv \wt C (\wt \delta) >0$ such that for all $z \in \rho(L)$ satisfying  $|z|> \wt r$ and  $|\arg z| \leq \pi/2$
\begin{equation}\label{L.res.F.est2}
\|(z-L)^{-1} \| \leq \frac{\wt C}{\dist(z, \sigma(L))} 
\exp \left( 
(1+\wt \delta)^2 \gamma \cos (\varrho \arg z) |z|^\varrho  
\right); 
\end{equation}
notice that the constants $M_1$ and $M_2$ depend on the choice of $\wt \delta$.

Finally, the claim \eqref{res.norm.L} follows by setting $\delta = (1+\wt \delta)^2-1$, $C := \wt C$ and $r := \wt r$. 
\end{proof}

Our main result on $L$ is the following corollary on the norms of spectral projections. It relies on the lower resolvent norm estimate in Theorem~\ref{thm:L.res.lb}, asymptotics of eigenvalues of $L$ in Proposition~\ref{prop:Shin} and Theorem~\ref{thm:L}.

\begin{theorem}\label{thm:L.Pn}
Let $L$ be as in \eqref{L.def} -- \eqref{V.deg}, $P_n$ as in \eqref{Pn.L.def}, $\kappa$  as in \eqref{kap.L.def} and $\tau$ as in \eqref{rho.L.lb}. Set 
\begin{equation}\label{sigma.L.def}
	\sigma := \kappa \tau = \frac{b-(a-1)}{1+a}.
\end{equation}
(Notice that $\sigma \in (0,1)$.)
Then there exists $\gamma>0$ such that
\begin{equation}\label{Pn.L.limsup}
	\limsup_{n \to \infty} \frac{\|P_n\|}{\exp (\gamma n^\sigma)} = \infty.
\end{equation}
In particular, the eigensystem of $L$ does not contain a basis.
\end{theorem}
\begin{proof}
We argue by contradiction. Suppose that the spectral projections satisfy 
\begin{equation}\label{Pn.L.asm.contr}
	\|P_n\| \leq M \exp \left(\wt \gamma (dn)^\sigma\right), \qquad n > N_0,
\end{equation}
with some $M>0$, $d>0$ as in \eqref{d.L.def}, $\sigma$ as in \eqref{sigma.L.def} and $\wt \gamma = \gamma d^{-\sigma}$.

Consider the points
\begin{equation}\label{zk.def}
z_k = t_k + \ii s_k, \quad k > N_0,
\end{equation}
where 
\begin{equation}\label{sk.tk.def}
t_k := \frac{\Re \la_k+\Re \la_{k+1}}{2}, \quad k > N_0, 
\quad 
s_k = 
\begin{cases}
0 & \text{if $b$ is odd},
\\
\frac12 c_b  t_k^\frac{b}{2a} &  \text{if $b$ is even},
\end{cases}
\end{equation}
where $N_0$ is as in Proposition~\ref{prop:Shin}.

Theorem~\ref{thm:L.res.lb} yields that there exists $\eta'>0$ such that for all $k > N_1$ with a sufficiently large $N_1>N_0$
\begin{equation}\label{L.zk.lb}
	\|(z_k-L)^{-1} \| 
	\geq \exp( \eta' t_k^\tau). 
	%\geq C_1 \exp(\eta' k^\sigma),
\end{equation}
Indeed, if $b$ is odd, then $\eta'=\eta$ from \eqref{L.res.lb.odd}. If $b$ is even, then we use \eqref{L.res.lb.even} with $\eps = c_b/2$ and obtain
\begin{equation}
\eta \frac{(\Im z_k)^\frac{b+1}b}{(\Re z_k)^\frac12} 
= \eta \left(\frac{c_b}2\right)^\frac{b+1}{b} t_k^{\frac{b+1}{2a}-\frac12}
=: \eta' t_k^\tau.
\end{equation}

On the other hand, by Theorem \ref{thm:L} with $\varrho = \sigma/\kappa = \tau \in (0,1/2)$, for any $\delta>0$ there exists $C>0$ such that for all $k > N_2$ with a sufficiently large $N_2>N_1$
\begin{equation}\label{res.norm.L.1}
	\|(z_k-L)^{-1} \| 
	\leq 
	\frac{C}{\dist(z_k, \sigma(L))} \exp \left( 
	(1+\delta) \wt \gamma \cos(\varrho \arg z_k) |z_k|^\tau  
	\right). 
\end{equation}
Notice that $\arg z_k \to 0$ as $k \to \infty$. Moreover, it follows from \eqref{la.L.gaps} that $\Re \la_k$ satisfy 
\begin{equation}
	\Re \la_{k+1}-\Re \la_k = \kappa d^\kappa  k^\frac{a-1}{a+1}(1+\BigO(k^{-\frac{1}{a+1}})), \quad k \to \infty,
\end{equation}
thus the sequence $\{\Re \la_k\}_{k \in \N}$ is eventually increasing and with gaps that are eventually uniformly positive (in fact, growing to infinity for $a>1$). Hence $\dist(z_k, \sigma(L))$ is eventually uniformly positive. In summary, for any $\delta>0$ there exist $N_3>N_2$ and $C'>0$ such that for all $k>N_3$
\begin{equation}\label{res.norm.L.2}
	\begin{aligned}
		\|(z_k-L)^{-1} \| 
		& \leq 
		C' \exp \left( 
		(1+ 2\delta) \wt \gamma t_k^\tau 
		\right).
	\end{aligned}
\end{equation}
Hence, we conclude that if $\wt \gamma < \eta'$, i.e., $\gamma < d^\sigma \eta'$, then \eqref{L.zk.lb} and \eqref{res.norm.L.2} are in contradiction. 

Finally, if the eigensystem of $L$ has a basis property, then the spectral projections $P_n$, $n > N_0$, are uniformly bounded. However, this is not satisfied due to \eqref{Pn.L.limsup}.
\end{proof}

\begin{remark}
In the proof of Theorem~\ref{thm:L.Pn} we used an upper bound for the resolvent of $L$ in \eqref{res.norm.L.2}, obtained from Theorem~\ref{thm:L}, for a special sequence $z_k$, $k >N_0$, in \eqref{zk.def}, \eqref{sk.tk.def}, which leads to a contradiction with the lower bound on the resolvent in \eqref{L.zk.lb} from Theorem~\ref{thm:L.res.lb}. Nonetheless, an appropriate adjustment of the upper bound \eqref{res.norm.L.2} remains valid for all $z$ with $\dist(z,\sigma(L)) \geq \eps$ for any $\eps>0$. In fact, one can also consider $z$ such that $\dist(z,\sigma(L)) \geq |z|^{-m}$ with some $m>0$ as $\Re z \to \infty$ which results in an additional term $ m \log|z| = o((\Re z)^\tau)$ in the exponential in \eqref{res.norm.L.2}. In either case, a contradiction with the lower bound \eqref{L.zk.lb} can be achieved.
\end{remark}

\begin{remark}\label{rem:L.shifted.HO}
In the special case \eqref{L.shift.HO.def} with $a=b=1$, it is known that 
\begin{equation}\label{Pn.shifted.HO}
	\|P_n\| = \frac{1}{2^\frac34} \frac{1}{\sqrt{|\pi \Im \alpha_1 |}} \frac{\exp \left(\sqrt 2 |\Im \alpha_1| n^\frac12 \right)}{n^\frac14}  
	\left(
	1+ \BigO(n^{-\frac12})
	\right), \quad n \to \infty,
\end{equation}
see \cite[Thm.~2.6]{Mityagin-2017-272}. Notice that for $a=b=1$, we have $\sigma=1/2$ in \eqref{Pn.L.limsup}, thus the exponential order $\sigma$ in Theorem~\ref{thm:L.Pn} is sharp. Moreover, \eqref{Pn.shifted.HO} and Theorem~\ref{thm:L} show that (recall \eqref{la.shifted.HO})
\begin{equation}\label{shifted.HO.res.ub}
	\|(z - L)^{-1}\| \leq  \frac{C}{\dist(z, \sigma(L))} \exp \left( 
	(1+\delta) |\Im \alpha_1| |z|^\frac12(1+o(1))  
	\right)  
\end{equation}
for all $z \in \rho(L)$ with  $|\arg z| \leq \pi/2$ and a sufficiently large $|z|$. Thus also the exponential order $\varrho=1/2$ in \eqref{L.res.lb.odd} is sharp in this case.

\emph{In fact, we expect that \eqref{Pn.L.limsup} holds with $\lim$, i.e.~not only for a subsequence, and also that the exponential order $\sigma$ in \eqref{Pn.L.limsup} is sharp in all cases in Theorem~\ref{thm:L.Pn}.}
\end{remark}

\subsection{Absence of basis property}
\label{ssec:basis}

Theorem~\ref{thm:L.Pn}  states in particular that the eigensystem of $L$ does not contain a basis as the spectral projections $P_n$, $n > N_0$, are not uniformly bounded. In fact, if $P_n$, $n > N_0$, were uniformly bounded, a much more restrictive upper bound on $\|(z-L)^{-1}\|$ than \eqref{res.norm.L} would hold, see \eqref{res.norm.L.basis}, \eqref{res.L.basis.impr} below. Thus only a mild growth of $\|(z-L)^{-1}\|$ suffices to exclude the basis properties of $L$. Resolvent estimates for the operators in $L^2(\R)$
\begin{equation}
\Dt + |x|^{2a} + \ii W(x)
\end{equation}
with $a>0$ and a range of perturbations $W$, resulting in various rates of the resolvent norm, can be found in \cite[Thm.~3.7, Ex.~3.8, Thm.~5.2, Ex.~5.3]{Krejcirik-2019-276}; for some details see Subsection~\ref{ssec:no.basis} below.
\begin{proposition}\label{prop:L:no.basis}
Let $L$ be as in \eqref{L.def} -- \eqref{V.deg} with $a>1$. Let $P_n$, $n >N_0$, be as in \eqref{Pn.L.def} and suppose that there exists $M >0$ such that
\begin{equation}\label{Pn.L.basis}
	\|P_n\| \leq M , \qquad n > N_0.
\end{equation}
Then there exist $C>0$ and $r>0$ such that for all $z \in \rho(L)$ satisfying $|z| \geq r$ and $|\arg z| \leq  \pi/2$
\begin{equation}\label{res.norm.L.basis}
\|(z-L)^{-1} \| 
 \leq 
\frac{C|z|}{\dist(z, \sigma(L))}.  
\end{equation}
\end{proposition}
\begin{proof}
We regularize the resolvent of $L$ like in \eqref{Davies.ident} using a very simple gauge function $\Phi(w) = w+1$. In fact, this amounts to the first resolvent identity 
\begin{equation}\label{res.id.L}
(z-L)^{-1} - (-1-L)^{-1} = (z+1)(z-L)^{-1} (1+L)^{-1}, \quad z \in \rho(L).
\end{equation}
For all $h \in \cD_L$, cf.~\eqref{B.L.def}, we have (recall \eqref{L.res.DL})
\begin{equation}
(z-L)^{-1} (1+L)^{-1} h = (1+L)^{-1} (z-L S_0)^{-1} S_0 h + \sum_{j=N_0+1}^{N(h)} \frac{1}{(1+\la_j)(z-\la_j)} P_j h.
\end{equation}
Hence, employing \eqref{res.S0.est.L}, \eqref{res.L.accr} with $a_n=1$ and \eqref{Pn.L.basis}, we obtain for $|z| \geq 2 \|L S_0\|$
\begin{equation}\label{Rz.R1.L.basis}
\begin{aligned}
\|(z-L)^{-1} (1+L)^{-1} h\| 1
& \leq 
\frac{2 \|S_0\|}{|z|} \| h\| + \sum_{j=N_0+1}^{\infty} \frac{M}{(1+\Re \la_j)|z-\la_j|}\|h\|
\\
& \leq 
\frac{2 \|S_0\|}{|z|} \| h\| + \frac{M}{\dist(z,\sigma(L))} \sum_{j=N_0+1}^{\infty} \frac{1}{1+\Re \la_j}\|h\|.
\end{aligned}
\end{equation}
Since $a>1$ and hence $\kappa>1$, see~\eqref{kap.L.def}, the asymptotic of $\Re \la_j$ in \eqref{tn.L.asym} yields
\begin{equation}
M':=\sum_{j=N_0+1}^{\infty} \frac{1}{1+\Re \la_j} < \infty.	
\end{equation}
Thus for all $|z| \geq 2 \|L S_0\|+1$,
\begin{equation}
\|(z-L)^{-1} (1+L)^{-1}\| 
\leq 
\frac{1}{\dist(z,\sigma(L))} 
\left(
\frac{2 \|S_0\|}{|z|} + M M'
\right)
 .
\end{equation}
Returning to \eqref{res.id.L}, we get that for all $|z| \geq 2 \|L S_0\|$
\begin{equation}
\begin{aligned}
\|(z-L)^{-1}\| 
&\leq 
\|(1+L)^{-1}\| + |z+1|\|(z-L)^{-1} (1+L)^{-1}\| 
\\
& \leq 1 + \frac{|z+1|}{\dist(z,\sigma(L))} 
\left(
\frac{2 \|S_0\|}{|z|} + M M'
\right),
\end{aligned}
\end{equation}
thus the claim \eqref{res.norm.L.basis} follows.
\end{proof}

\begin{remark}
Let the spectral projections $P_n$, $n>N_0$, satisfy a polynomial bound 
\begin{equation}
\|P_n\| \leq M n^\sigma, \quad n > N_0,	
\end{equation}
where $M,\sigma>0$. Then there exist $C>0$ and $r>0$ such that  for all $z \in \rho(L)$ satisfying $|z| \geq r$ and $|\arg z| \leq  \pi/2$
\begin{equation}\label{res.norm.L.pol}
\|(z-L)^{-1} \| 
\leq 
\frac{C|z|^m}{\dist(z, \sigma(L))}, 
\end{equation}
where
\begin{equation}
	m= \left \lfloor \frac{\sigma+1}{\kappa}  \right \rfloor +1.
\end{equation}
The proof is analogous to the one of Proposition~\ref{prop:L:no.basis}, using the polynomial gauge \eqref{F.polyn} with this $m$ and $a_k = k$. The bound \eqref{res.norm.L.pol} can be viewed as a more quantitative version of \cite[Thm.~3]{Davies-2000-32}.
\end{remark}

\begin{remark}\label{rem:L.basis.precise}
	In fact, a more precise estimate of the series in \eqref{Rz.R1.L.basis}, and hence an improvement of of \eqref{res.norm.L.basis}, is possible. Namely, for any $\eps>0$, there exists $C'>0$ such for all $z \in \rho(L)$ satisfying $|z| \geq e$,  $|\arg z| \leq \pi/2$ and $\dist(z,\sigma(L)) \geq \eps$, we have
	\begin{equation}
		\sum_{j=N_0+1}^{\infty} \frac{1}{(1+\Re \la_j)|z-\la_j|} \leq  C'
		\begin{cases}
			\displaystyle \frac{\log|z|}{|z|} & \text{if $a= 1$},
			\\
			\frac{1}{|z|} & \text{if $a > 1$}.
		\end{cases}
	\end{equation}
	Such an estimate can be proved like in \cite[Lemma~4.3]{Mityagin-2019-139} employing the asymptotics of $\Re \la_j$ from \eqref{tn.L.asym} and \eqref{la.L.gaps}. Hence for all $z \in \rho(L)$ satisfying $|z| \geq r$,  $|\arg z| \leq \pi/2$ and $\dist(z,\sigma(L)) \geq \eps$, the estimate  \eqref{res.norm.L.basis} can be improved to
	\begin{equation}\label{res.L.basis.impr}
		\|(z-L)^{-1} \| 
		\leq 
		\wt C
		\begin{cases}
			\log|z| & \text{if $a=1$},
			\\	
			1 & \text{if $a>1$},
		\end{cases}
	\end{equation}
	with $\wt C\equiv \wt C(\eps) >0$.
\end{remark}

\section{Faster resolvent growth}
\label{sec:T}

The proof of Theorem~\ref{thm:L} indicates that the same claim can be obtained for abstract operators $T$ which share properties of $L$ (indeed, this is Theorem~\ref{thm:T} with $p=1$ below, see Remark~\ref{rem:L.T}). Nonetheless, the strategy is not applicable if the exponential order of the resolvent norm is larger than $1/2$, e.g.~for the imaginary anharmonic oscillators, see \eqref{H.even.def} and \eqref{H.res.lb.even} below. The limitation originates in the partial fraction decomposition in Proposition~\ref{prop:F.pfc.new} and can be overcome under additional assumptions on the location of spectrum and the resolvent growth. 

The starting point is the decomposition \eqref{pf.decomp} or \eqref{pf.decomp.12}. For any $p \in \N$, a substitution and a straightforward manipulation yield in the setting of Proposition~\ref{prop:F.pfc.new} that
\begin{equation}
	\frac{1}{(z^p-w^p) F(w^p)}  = \frac{1}{F(z^p)(z^p-w^p)} + \sum_{n=1}^\infty \frac{1}{(z^p+a_n) F'(-a_n)(w^p+a_n)},	
\end{equation}
thus
\begin{equation}\label{pf.decomp.p}
	\frac{1}{(z-w) F(w^p)}  = \frac{1}{F(z^p)(z-w)} + \sum_{n=1}^\infty \frac{\sum_{k=0}^{p-1} w^k z^{p-1-k}}{(z^p+a_n) F'(-a_n)(w^p+a_n)}.	
\end{equation}
Since $|F'(-a_n)|$ grow exponentially, see \eqref{F'.an.asym}, the additional terms in the sum in \eqref{pf.decomp.p} can be compensated if additionally the resolvent norm of $T$ does not grow too fast on the special rays in $\C$. For instance for $p=2$, the main term to control reads
\begin{equation}
	\begin{aligned}
		\|T (a_n+T^2)^{-1}\| & \leq  \| T (\ii a_n^\frac12+T)^{-1}\| \| (-\ii a_n^\frac12+T)^{-1}\|
		\\
		& 
		\leq (1+a_n^\frac12 \|(\ii a_n^\frac12+T)^{-1}\|) \| (-\ii a_n^\frac12+T)^{-1}\|,
	\end{aligned}
\end{equation}
so the resolvent norm of $T$ on the imaginary axis, at least, on the sequences $\{\ii a_n^{1/2}\}$, $\{-\ii a_n^{1/2}\}$, should grow sufficiently slower comparing to $|F'(-a_n)|$.

Similar steps in the case $\varrho = 1/2$ and \eqref{pf.decomp.12} give the identity
\begin{equation}\label{pf.decomp.p.12}
	\frac{1}{(z-w) F(w^p)}  = \frac{1}{F(z^p)(z-w)} + 2 \sum_{n=1}^\infty (-1)^{n+1} \frac{a_n}{z^p + a_n} \frac{\sum_{k=0}^{p-1} w^k z^{p-1-k}}{w^p+a_n}.	
\end{equation}

\subsection{From projections to the resolvent for an abstract $T$}

\begin{asm-sec}\label{asm:T}
	Let $T$ be a closed operator with compact resolvent in a separable Hilbert space $\cH$ and denote the eigenvalues of $T$ by $\la_n$, $n \in \N$ (repeated with the algebraic multiplicity and ordered such that $|\la_{n+1}|\geq|\la_n|$). Let $p \in \N$ and assume that	
	\begin{enumerate}[\upshape (i), wide]
		\item \label{asm:T.rho.sect} the spectrum of $T$ is contained in a union of a disc and a sector (with some $R>0$ and $\epsilon_0 \in (0,1]$)
		\begin{equation}\label{T.sp.encl}
			\sigma(T) \subset B_{R}(0) \cup \left\{ z \in \C \, : \, |\arg z| \leq (1 - \epsilon_0) \frac{\pi}{p} 
			\right\}; 
		\end{equation}		
		\item \label{asm:T.res.est} the resolvent norm is polynomially bounded on the following rays in $\C$ 
		\begin{equation}\label{T.res.pol.asm}
			\|(r e^{\ii \frac{\pi}{p}(1+2k)} - T)^{-1}\| \leq K r^N , \quad r \geq R+1 , \quad k = 0, 1, \dots, p-1,
		\end{equation}
		with some $K \geq 1$ and $N \in [-1,\infty)$;
		\item \label{asm:T.EV.simple} there exists $R_0 > R > 0$ such that the geometric and algebraic multiplicities coincide for all $\la \in \sigma(T)$ satisfying $|\la| > R_0$;

		\item\label{asm:T.compl} the system of spectral projections $Q_\la$, $\la \in \sigma(T)$, is complete in $\cH$;
		
		\item  the eigenvalues $\la_n$  satisfy
		\begin{equation}\label{la.n.asm}
			|\la_n| \geq c n^\kappa, \qquad n \in \N,
		\end{equation}
		with some $c,\kappa>0$. 
	\end{enumerate}
\end{asm-sec}

\begin{remark}\label{rem:L.T}
The operator $L$ in \eqref{L.def} satisfies Assumption~\ref{asm:T} with $p=1$. By \eqref{Num.L.odd}, \eqref{Num.L.even}, for any $\epsilon_0 \in (0,1)$, there exists $R>0$ such that the spectrum of $L$ satisfies \eqref{T.sp.encl} with $p=1$; moreover, the resolvent estimate \eqref{T.res.pol.asm} holds with $N=-1$, see Subsection~\ref{ssec:ot}, in particular around \eqref{T.Num.par}. The eigenvalues $\la_n$ of $L$ are simple for $n>N_0$ by Proposition~\ref{prop:Shin} and for $n >N_0$ the spectral projections $P_n \equiv Q_{\la_n}$ have $\Rank$ 1. The completeness of the spectral projections follows from Proposition~\ref{prop:L.basic} \ref{prop:L.basic.complete}. Finally, the eigenvalues asymptotics in Proposition~\ref{prop:Shin} justify \eqref{la.n.asm}. In fact, to verify that \eqref{la.n.asm} holds, it suffices to use the behavior of singular values \eqref{sk.L}, see \eqref{lak.sk.ineq}.
\end{remark}

\begin{theorem}\label{thm:T}
Let Assumption~\ref{asm:T} be satisfied. Assume that 
	\begin{equation}\label{Qla.asm}
		\|Q_\la\| \leq \exp \left(\gamma |\la|^{p\varrho}(1+o(1))\right), \quad \la \in \sigma(T), \ |\la| \to + \infty,
	\end{equation}
	for some $\gamma >0$ and $\varrho>0$. If either $\varrho < 1/2 $ or 
	\begin{equation}\label{asm.N.12}
		\varrho = \frac12 \quad \text{and} \quad  N < -1 + \frac 1 {2p},	
	\end{equation}	
	then for any $\epsilon_1 \in (0,1]$ and any $\delta>0$, there exist $C \equiv  C(\epsilon_1,\delta)>0$ and $r \equiv r(\epsilon_1,\delta)>0$ such that for all $z \in \rho(T)$ satisfying $|z|> r$ and $|\arg z| \leq (1 - \epsilon_1) \pi/p$,
	\begin{equation}\label{res.norm.T.p}
		\begin{aligned}
			\|(z-T)^{-1} \| 
			& \leq 
			\frac{C}{\dist(z, \sigma(T))} \exp  \left( 
			(1+\delta) \gamma \frac{ \cos (p \varrho \arg z)}{\cos \left( 
				(1-\epsilon_0) \pi \varrho \right)} |z|^{p\varrho}
			\right).
		\end{aligned}
	\end{equation}
\end{theorem}

\begin{proof}
We follow the steps in the proof of Theorem~\ref{thm:L}. 

Consider $F$ as in \eqref{F.def} with $a_k$, $k \in \N$, satisfying \eqref{an.special}, i.e.,
\begin{equation}\label{an.special.again}
	a_n := \left( \frac{n}{\nu} \right)^{\frac 1 \varrho}, \qquad n \in \N,
\end{equation}
 where $\varrho \leq 1/2$ is as in \eqref{Qla.asm} and
\begin{equation}\label{nu.p}
	\nu :=  \gamma \frac{\sin\left(\pi \varrho \right)}{\pi} \frac{1}{\cos \left( 
		(1-\epsilon_0) \pi \varrho 
		\right)}(1+ \delta).
\end{equation}

Notice that it suffices to prove the claim \eqref{res.norm.T.p} for all $\delta \in (0,\delta_*)$ with some $\delta_* \in (0,1)$ (since then \eqref{res.norm.T.p} holds for any $\delta \geq \delta_*$ as well). We select $\delta_*>0$ so small that $F(\la_n^p) \neq 0$ for all $n \in \N$ and all $\delta \in (0,\delta_*)$. This is possible as there are only finitely many $\la_n^p <0$, hence only for finitely many $\delta \in (0,1)$ there could exist $n \in \N$ such that $\la_n^p \in \{-a_m\}_{m \in \N}$.

For each $z \in \rho(T)$, consider the operator 
\begin{equation}\label{B.def.p}
	B_z h := \sum_{\substack{\la \in \sigma(T) \\ R_0 < |\la| < R(h) }} \frac{1}{z-\la} \frac{1}{F(\la^p)} Q_\la h,  
\end{equation}
defined on a dense subspace of $\cH$
\begin{equation}
\begin{aligned}
\cD_T :=  \big\{ h \in \cH \,: \,  \text{there exists } \ & R(h) > R_0  \  \text{such that}\ Q_\la h = 0  
	\\ 
	& \quad \text{for all} \ \la \in \sigma(T) \ \text{with } \ |\la|>R(h) \big\},	
\end{aligned}
\end{equation}
cf.~\eqref{DL.def}. Hence by Proposition~\ref{prop:F.Titch} and the assumption \eqref{T.sp.encl} on the location of the spectrum of $T$, we can select $r_\delta > R_0$ so large that we have
\begin{equation}\label{F.laj.p.geq}
	\begin{aligned}
		|F(\la^p)| &\geq 
		\exp \left(
		\frac{1}{1+ \frac \delta 4} \pi \nu |\la^p|^\varrho \frac{\cos\left( \arg(\la^p) \varrho \right)}{\sin ( \pi \varrho )}
		\right)
		\\
		&\geq 
		\exp \left(
		\frac{1}{1+ \frac \delta 4} \pi \nu |\la|^{p\varrho} \frac{\cos\left((1-\epsilon_0) \pi \varrho \right)}{\sin ( \pi \varrho )}
		\right)
		\\
		&=
		\exp \left(
		\frac{4+4 \delta}{4+ \delta} \gamma |\la|^{p\varrho} 
		\right),  \quad \la \in \sigma(T), \ |\la| >  r_\delta;
	\end{aligned}
\end{equation}
the choice of $\nu$ in \eqref{nu.p} is used in the last step. The hypothesis~\eqref{Qla.asm} yields that $r_\delta$ can be in addition selected so large that 
\begin{equation}
	\|Q_\la\| \leq \exp \left(\gamma |\la|^{p\varrho} \frac{4+2\delta}{4+\delta} \right), \quad \la \in \sigma(T), \ |\la| > r_\delta.
\end{equation}
Hence
\begin{equation}\label{B.norm.p}
	\|B_z h\| \leq \frac{1}{\dist(z,\sigma(T))} \sum_{\substack{\la \in \sigma(T)\\ R_0 < |\la| \leq  R(h) }} \frac{\|Q_\la\|}{|F(\la^p)|}  \|h\|
	\leq 
	\frac{M_1}{\dist(z,\sigma(T))} \|h\|
\end{equation}
where
\begin{equation}\label{M1.def.p}
	M_1= 
	\sum_{\substack{\la \in \sigma(T)\\ R_0 < |\la|  }} \frac{\|Q_\la\|}{|F(\la^p)|} 
	\leq 
	\sum_{\substack{\la \in \sigma(T)\\ R_0 < |\la|  }} \exp \left(-\gamma \frac{2 \delta}{4+\delta} |\la|^{p \varrho} \right).
\end{equation}
The condition \eqref{la.n.asm} guarantees that the right sum in \eqref{M1.def.p} is finite.
%We have $|\la_n| \geq dn^\kappa$, $n \in \N$, by the assumption \eqref{la.n.asm}. Let $\mu_m$, $m\in \N$, be eigenvalues of $T$ enumerated \emph{without} the repetition according to the algebraic multiplicity and ordered such that $|\mu_{m+1}| \geq |\mu_m|$. If $|\mu_m| = |\la_n|$ for some $n \in \N$, then $n \geq m$ and so
%%
%\begin{equation}
%	|\mu_m| = |\la_n| \geq d n^\kappa \geq d m^\kappa, \quad m \in \N.
%\end{equation}
%%
%Hence we get
%%
%\begin{equation}
%	M_1 \leq 
%	\sum_{m=1}^\infty \exp \left(-\gamma \frac{2 \delta}{4+\delta} |\mu_m|^{p \varrho} \right)
%	\leq 
%	\sum_{m=1}^\infty \exp \left(-\gamma \frac{2 \delta}{4+\delta} |d m^\kappa|^{p \varrho} \right) < \infty.
%\end{equation}

In the second step, we employ and analyze the expansion \eqref{pf.decomp.p} if $\varrho <1/2$ and \eqref{pf.decomp.12} if $\varrho =1/2$. To this end, we first prepare the estimates \eqref{Tk.res.p} and \eqref{zp.ap.ineq}.

We have
\begin{equation}
T (T - \zeta)^{-1} = 1 + \zeta (T - \zeta)^{-1}, \quad \zeta \in \rho(T),	
\end{equation}
and for any $v > 0$ and $k = 0, \dots, p-1$
\begin{equation}
\begin{aligned}
v^p + T^p & = \prod_{j=0}^{p-1} \left(T - v e^{\ii \frac \pi p (1+2j)} \right),
\\
T^k (v^p + T^p)^{-1}  
 &=\left(
\prod_{j=0}^{k-1} T \left(T - v e^{\ii \frac \pi p (1+2j)} \right)^{-1}
\right)
\left( 
\prod_{j=k}^{p-1} \left(T - v e^{\ii \frac \pi p (1+2j)} \right)^{-1} 
\right).
\end{aligned}
\end{equation}
Thus using the condition~\eqref{T.res.pol.asm}, for any $k = 0, \dots, p-1$ and $v \geq R +1$ 
\begin{equation}\label{Tk.res.p}
\begin{aligned}
 \|T^k(v^p+T^p)^{-1}\| 
& \leq 
\left( \prod_{j=0}^{k-1} \left(1+ v \|(v e^{\ii \frac{\pi}{p}(1+2j)} - T)^{-1} \| \right)	 \right)
\prod_{j=k}^{p-1} \|(v e^{\ii \frac{\pi}{p}(1+2j)} - T  )^{-1} \|	
\\
& 
\leq \left(\prod_{j=0}^{k-1}(1+ K v^{1 +  N }) \right) \prod_{j=k}^{p-1} K v^{N}
\leq  (2K)^p v^{p N + k}.
\end{aligned}
\end{equation}

Next, we verify that for all $z$ satisfying $|\arg z| \leq (1 - \epsilon_1) \pi / p $ and all $v>0$,
\begin{equation}\label{zp.ap.ineq}
|z^p + v^p| \geq 2^{1-p} \left( \sin \frac {\pi \epsilon_1} 2 \right)  (|z| + v)^p.
\end{equation}
Indeed, let $\psi \equiv \pi (1 - \epsilon_1)$ and notice first that
\begin{equation}
|p \arg z | \leq \psi, \quad \cos (p \arg z ) \geq \cos \psi, \quad \cos \frac \psi 2 = \sin \frac{\pi \epsilon_1}2.
\end{equation}
Then \eqref{zp.ap.ineq} can be justified as follows
\begin{equation}
\begin{aligned}
|z^p + v^p|^2 & = \left(|z|^p \cos (p \varphi) + v^p \right)^2 + \left(|z|^p \sin (p \varphi)\right)^2	
\\
& = |z|^{2p} + 2 |z|^p v^p \cos(p \arg z) + v^{2p} 
\geq 
|z|^{2p} + 2 |z|^p v^p \cos \psi + v^{2p}
\\
& = (1-\cos \psi)(|z|^{2p} + v^{2p}) +  (|z|^p+v^p)^2 \cos \psi
\\
& \geq \left(
\frac12(1- \cos \psi) + \cos \psi
\right)(|z|^p+v^p)^2
= \frac12(1+ \cos \psi)(|z|^p+v^p)^2
\\
& = \left(\cos \frac \psi 2\right)^2 (|z|^p+v^p)^2 = \left(\sin \frac{\pi \epsilon_1}2 \right)^2 (|z|^p+v^p)^2 \\
&\geq 2^{2(1-p)} \left(\sin \frac{\pi \epsilon_1}2 \right)^2 (|z|+v)^{2p}.
\end{aligned}
\end{equation}

Further we treat the two cases separately.

\noindent $\bullet$ $\varrho <1/2:$
Similarly as in \eqref{Bh.L.PFD}, we obtain that for all $h \in \cD_T$ and $z \in \rho(T)$ satisfying $|\arg z| \leq (1 - \epsilon_1) \pi/p$,
\begin{equation}\label{Bh.PFD.p}
\begin{aligned}
B_z h & = 	
\frac1{F(z^p)} (z-T)^{-1} (h-S_0h)
\\
& \quad+  \sum_{k=0}^{p-1} \sum_{n=1}^\infty \frac{z^{p-1-k}}{z^p+a_n}    \frac{1}{F'(-a_n)}  T^k(a_n+T^p)^{-1}(h-S_0h).
\end{aligned}
\end{equation}
Notice that $-a_n \in \rho(T^p)$, $n \in \N$, since $a_n$ were selected above such that $F(\la_n^p) \neq 0$, $n \in \N$.
Recalling that $a_n>0$ and $|F'(-a_n)|$ grows exponentially, see \eqref{F'.an.asym}, using \eqref{Tk.res.p} and \eqref{zp.ap.ineq}, we obtain that for $z$ satisfying $|\arg z| \leq (1 - \epsilon_1) \pi/p$
\begin{equation}\label{RHS.sum.p}
\begin{aligned}
& \left\| 
\sum_{k=0}^{p-1} \sum_{n=1}^\infty \frac{z^{p-1-k}}{z^p+a_n}    \frac{1}{F'(-a_n)}  T^k(a_n+T^p)^{-1}(I-S_0)
\right\|
\\
& \quad \leq \frac {\|(I-S_0)\|} {2^{1-p}\sin \frac {\pi \epsilon_1} 2}  	
\sum_{k=0}^{p-1} \sum_{n=1}^\infty \frac{|z|^{p-1-k}}{(|z|+a_n^\frac1p)^p}    \frac{\|T^k(a_n+T^p)^{-1} \|}{|F'(-a_n)|}  
\\
& \quad \leq \frac {\|(I-S_0)\|} {2^{1-p} \sin \frac {\pi \epsilon_1} 2}  
\left(
\sum_{n=1}^\infty    \frac{\|T^k(a_n+T^p)^{-1} \|}{|F'(-a_n)|}
\right)  
\left(\sum_{k=0}^{p-1} \frac{|z|^{p-1-k}}{|z|^p} \right) 
= \BigO(|z|^{-1}), 
\end{aligned}
\end{equation}
as $z \to \infty$. Rearranging \eqref{Bh.PFD.p} similarly as in \eqref{B.res.rearr}, using the density of $\cD_T$ in $\cH$ and the estimates \eqref{B.norm.p}, \eqref{RHS.sum.p} and an analogue of \eqref{res.S0.est.L}, we get that for $z$ satisfying $|\arg z| \leq (1 - \epsilon_1) \pi/p$
\begin{equation}\label{T.res.p.bound<12}
\begin{aligned}
\|(z-T)^{-1} \| & \leq |F(z^p)| \left(
\frac{M_1}{\dist(z, \sigma(T))}
+
\BigO(|z|^{-1})
\right) + \BigO(|z|^{-1}), \quad z \to \infty.
\end{aligned}
\end{equation}

\noindent 
$\bullet$ $\varrho = 1/2$: 	
We use \eqref{pf.decomp.p.12} instead of \eqref{pf.decomp.p}. We rewrite \eqref{pf.decomp.p.12} in the form
\begin{equation}
	\begin{aligned}
		\frac{1}{F(z^p)(z-w)}  = \frac{1}{F(w^p)(z-w)} 
		&- 2 \sum_{n=1}^\infty (-1)^n \frac{z^p}{z^p+a_n} \frac{\sum_{k=0}^{p-1} w^k z^{p-1-k}}{w^p+a_n} 
		\\ 
		& 
		+ 2 \sum_{n=1}^\infty (-1)^n \frac{\sum_{k=0}^{p-1} w^k z^{p-1-k}}{w^p+a_n}.
	\end{aligned}
\end{equation}
Like in \eqref{Bh.L.PFD} we obtain for all $h \in \cD_T$ and $z \in \rho(T)$ that
\begin{equation}\label{B.res.id.12}
	\begin{aligned}
		\frac1{F(z^p)} (z-T)^{-1} (h-S_0h) =  B_z h 
		& -2 \sum_{k=0}^{p-1} \sum_{n=1}^\infty  (-1)^n \frac{z^{2p-1-k}}{z^p+a_n}  T^k(a_n+T^p)^{-1} h
		\\
		& +2 \sum_{k=0}^{p-1} z^{p-1-k} \sum_{n=1}^\infty (-1)^n   T^k(a_n+T^p)^{-1}h.
	\end{aligned}
\end{equation}
The first sum on the r.h.s.~of \eqref{B.res.id.12} can be estimated in a straightforward way using \eqref{Tk.res.p} and \eqref{zp.ap.ineq}. Namely, for $z$ satisfying $|\arg z| \leq (1 - \epsilon_1) \pi/p$, we have for all $a_n \geq (R+1)^p$
\begin{equation}
\begin{aligned}
\sum_{k=0}^{p-1} \frac{|z|^{2p-1-k}}{|z^p+a_n|}  \|T^k(a_n+T^p)^{-1}\| 
	&\leq
	(2K)^p \frac { 2^{p-1}} { \sin \frac {\pi \epsilon_1}2}  \sum_{k=0}^{p-1} \frac{|z|^{2p-1-k}}{(|z|+a_n^\frac1p)^{p-k+k}} a_n^{N + \frac kp}
\\
&	\leq p (2K)^p
	\frac {2^{p-1}} {  \sin \frac {\pi \epsilon_1} 2} |z|^{p-1} a_n^{N}.
\end{aligned}
\end{equation}
Using the condition on $N$ in \eqref{asm.N.12}, we have 
\begin{equation}
	N  < -1 + \frac{1}{2p} \leq - \frac12, 
\end{equation}
thus, recalling \eqref{an.special.again} with $\varrho=1/2$, we obtain that for $z$ satisfying $|\arg z| \leq (1 - \epsilon_1) \pi/p$
\begin{equation}
	\left\| \sum_{k=0}^{p-1} \sum_{n=1}^\infty  (-1)^n \frac{z^{2p-1-k}}{z^p+a_n}  T^k(a_n+T^p)^{-1}  \right\|
	= \BigO(|z|^{p-1}), \quad z \to \infty.
\end{equation}

We rewrite the second sum on the r.h.s.~of \eqref{B.res.id.12} as 
\begin{equation}
	\begin{aligned}
		& \sum_{k=0}^{p-1} z^{p-1-k} \sum_{n=1}^\infty T^k \left( (a_{2n}+T^p)^{-1} - (a_{2n-1}+T^p)^{-1}\right) h
		\\
		& \quad =
		\sum_{k=0}^{p-1} z^{p-1-k} \sum_{n=1}^\infty T^k \left(  (a_{2n-1} - a_{2n}) (a_{2n}+T^p)^{-1}(a_{2n-1}+T^p)^{-1}\right) h.
	\end{aligned}
\end{equation}
Using \eqref{Tk.res.p}, we obtain that for all $k = 0, \dots, p-1$ and all $a_n \geq (R+1)^p$
\begin{equation}
\begin{aligned}
& \|(a_{2n-1} - a_{2n}) T^k (a_{2n}+T^p)^{-1}(a_{2n-1}+T^p)^{-1}\| 
\\
& \quad \leq 
|a_{2n-1} - a_{2n}| \|T^k (a_{2n}+T^p)^{-1}\| \|(a_{2n-1}+T^p)^{-1}\| 
\\
& \quad \leq |a_{2n-1} - a_{2n}| (2 K)^p a_{2n}^{N+\frac kp} K a_{2n-1}^N
= \BigO(a_n^{2N + \frac32 - \frac1p}), \quad n \to \infty,
\end{aligned}
\end{equation}
where we use that $|a_{2n-1} - a_{2n}| = \BigO(a_n^{1/2})$ as $n\to \infty$ in the last step, see \eqref{an.special.again} with $\varrho=1/2$.
From the condition on $N$ in \eqref{asm.N.12}, we have 
\begin{equation}
	2N + \frac32 - \frac1p < -\frac12, 
\end{equation}
thus  we obtain that for $z$ satisfying $|\arg z| \leq (1 - \epsilon_1) \pi/p$
\begin{equation}
	\left\|\sum_{n=1}^\infty (-1)^n \sum_{k=0}^{p-1} z^{p-1-k} T^k(a_n+T^p)^{-1} \right\| =  \BigO(|z|^{p-1}), \quad z \to \infty.
\end{equation}

Analogous arguments as in the justification of \eqref{T.res.p.bound<12} yield that for $z$ satisfying $|\arg z| \leq (1 - \epsilon_1) \pi/p$
\begin{equation}\label{T.res.p.bound=12}
	\begin{aligned}
		\|(z-T)^{-1} \| & \leq |F(z^p)| \left(
		\frac{M_1}{\dist(z, \sigma(T) )}
		+
		\BigO(|z|^{p-1})
		\right) + \BigO(|z|^{-1}), \quad z \to \infty.
	\end{aligned}
\end{equation}

Finally, the claim \eqref{res.norm.T.p} follows from \eqref{T.res.p.bound<12} or \eqref{T.res.p.bound=12},  \eqref{F.asym.Titch} and the choice of $\nu$ in \eqref{nu.p}.  (Notice that  $\BigO(|z|^{p-1})$ in \eqref{T.res.p.bound=12} can be absorbed in the remainder $o(1)$ in \eqref{res.norm.T.p}. For some further details see the part of the proof of Theorem~\ref{thm:L} below \eqref{L.res.F.est1}.)
\end{proof}

\begin{remark}\label{rem:theta}
Suppose that, in addition to the assumptions of Theorem~\ref{thm:T}, 
\begin{equation}\label{asm.thet0}
	\lim_{n \to \infty} \arg \la_n = \vartheta_0 \in (-\tfrac \pi p, \tfrac \pi p);
\end{equation}
for an example of such operators see Subsection~\ref{ssec:H.even} below. Then for any $\epsilon_1 \in (0,1]$ and any $\delta>0$, there exist $C \equiv C(\epsilon_1,\delta)>0$ and $r\equiv r(\epsilon_1,\delta)>0$ such that for all $z \in \rho(T)$ satisfying $|z|> r$ and $|\arg z| \leq (1 - \epsilon_1) \pi/p$,
\begin{equation}\label{res.norm.T.p.theta}
	\begin{aligned}
		\|(z-T)^{-1} \| 
		& \leq 
		\frac{C}{\dist(z, \sigma(T))} \exp  \left( 
		(1+\delta) \gamma \frac{ \cos (p \varrho \arg z)}{\cos \left( 
			p \varrho \vartheta_0 \right)} |z|^{p\varrho}
		\right).
	\end{aligned}
\end{equation}
Indeed,  instead of $\cos\left( \arg(\la^p) \varrho \right) \geq \cos\left((1-\epsilon_0) \pi \varrho \right)$ in \eqref{F.laj.p.geq}, it suffices to use that
\begin{equation}
\cos\left( \arg(\la^p) \varrho \right) = \cos(p \varrho \vartheta_0)(1+o(1)), \quad \la \in \sigma(T), \ \la \to \infty,
\end{equation}
and adjust \eqref{nu.p} as well as the second estimate in \eqref{F.laj.p.geq} accordingly. We omit the details.
\end{remark}

We formulate a corollary of Theorem~\ref{thm:T} for operators $T$ having numerical range contained in a ``parabolic'' region, see \eqref{T.Num.par} in Subsection~\ref{ssec:ot}. For an example of such an operator see $L$ in \eqref{L.def} and Proposition~\ref{prop:L.basic}.

\begin{corollary}\label{cor:T.par}
	Let $T$ be a closed operator with compact resolvent in a separable Hilbert space $\cH$ and let the numerical range of $T$ satisfy \eqref{T.Num.par}, \eqref{T.Num.h}.
	Let the singular values $\frs_k$ of $T^{-1}$ satisfy 
		\begin{equation}\label{cor.T.par.sk}
			\frs_k = \BigO(k^{-\frac1 \kappa}), \quad k \to \infty,
		\end{equation}
	for some $\kappa>0$ and let there exist $R_0 > 0$ such that the geometric and algebraic multiplicities coincide for all $\la \in \sigma(T)$ satisfying $|\la| > R_0$.
	Assume that
	\begin{equation}\label{Qla.asm.cor}
		\|Q_\la\| \leq \exp \left(\gamma |\la|^{p\varrho}(1+o(1))\right), \quad \la \in \sigma(T),  \ |\la| \to  +\infty,
	\end{equation}
	for some $\gamma>0$, $p \in \N$ and $\varrho \in (0,1/2]$. 
	Then for any $\epsilon_1 \in (0,1]$ and any $\delta>0$, there exist $C \equiv  C(\epsilon_1,\delta)>0$ and $r \equiv r(\epsilon_1,\delta)>0$ such that for all $z \in \rho(T)$ satisfying $|z|> r$ and $|\arg z| \leq (1 - \epsilon_1) \pi/p$, 
	\begin{equation}\label{res.norm.T.p.cor}
		\begin{aligned}
			\|(z-T)^{-1} \| 
			& \leq 
			\frac{C}{\dist(z, \sigma(T))} \exp  \left( 
			(1+\delta) \gamma  \cos (p \varrho \arg z) |z|^{p\varrho}
			\right).
		\end{aligned}
	\end{equation}
\end{corollary}
\begin{proof}
We verify the conditions in Assumption~\ref{asm:T}. The spectrum of $T$ is contained in $\Num(T)$, thus \eqref{T.sp.encl} holds for any $p \in \N$ and any $\epsilon_0>0$ (with a sufficiently large $R>0$). In fact, in this case, we have even \eqref{asm.thet0} with $\vartheta_0 =0$. 
Also \eqref{la.n.asm} holds due to \eqref{cor.T.par.sk}, see Subsection~\ref{ssec:ot}.

We have the resolvent bounds \eqref{T.Num.par.res}, thus \eqref{T.res.pol.asm} is satisfied with $N=-1$ (notice also that if $\varrho =1/2$, then the restriction on $N$ in \eqref{asm.N.12} is satisfied). Since $T^{-1} \in \cS_{\frac 1\kappa + \eps}$ for any $\eps>0$, see \eqref{cor.T.par.sk}, the completeness of spectral projections $Q_\la$, $\la \in \sigma(T)$, follows by Proposition~\ref{prop:compl}, see Subsection~\ref{ssec:ot}. 

The equality of the algebraic and geometric multiplicities of eigenvalues is assumed. 

Hence the claim follows by Theorem~\ref{thm:T}, being applied with $p \in \N$ and $\varrho \in (0,1/2]$ as in \eqref{Qla.asm.cor}, and Remark~\ref{rem:theta} with $\vartheta_0=0$.
\end{proof}

Like in Theorem~\ref{thm:L.Pn}, a lower bound on the resolvent norm can be propagated to spectral projections.

\begin{corollary}\label{cor:T.Pn}
	Let $T$ satisfy Assumption~\ref{asm:T} for $p \in \N$ and let there exist $\epsilon_1 \in (0,1]$ and a continuous curve 
	\begin{equation}
		\Gamma: \R_+ \to \left\{z \in \C \, : \, \Re z \geq 1, |\arg z| \leq (1-\epsilon_1)\frac{\pi}{p} \right\}	
	\end{equation}
	such that $|\Gamma(t)| \to + \infty$ as $t \to + \infty$ and for some $\eta,r>0$ and $\tau \leq  p/2$, 
	\begin{equation}\label{res.lb.limsup}
		\|(z-T)^{-1}\| \geq \exp(\eta |z|^\tau), \quad z \in \Im(\Gamma) \cap \rho(T), \ |z| > r. 
	\end{equation}
	If $\tau=p/2$, assume in addition that $N$ in \eqref{T.res.pol.asm} satisfies $N<-1+1/(2p)$.
	
	 Then for every $\gamma \in (0,\gamma_0)$, where (with $\epsilon_0$ from \eqref{T.sp.encl})
	\begin{equation}
		\gamma_0 = \eta \cos \left( (1-\epsilon_0) \frac \tau p  \pi \right),
	\end{equation}
	we have
	\begin{equation}\label{cor.T.Pn}
		\limsup_{n \to \infty} \frac{\|Q_{\la_n}\|}{\exp(\gamma |\la_n|^\tau)} = \infty. 
	\end{equation}
	
	If, in addition, 
	\begin{equation}\label{cor.T.thet0}
	\lim_{n \to \infty} \arg \la_n = \vartheta_0 \in (-\tfrac \pi p, \tfrac \pi p), 
	\end{equation}
	then \eqref{cor.T.Pn} holds for all $\gamma \in (0,\gamma_0)$ with 
	\begin{equation}\label{gam.0.theta}
	\gamma_0 = \eta \cos \left( \tau \vartheta_0 \right).
	\end{equation}
\end{corollary}
\begin{proof}
	We argue by contradiction. Suppose that there exists $r>0$ such that for all $\la \in \sigma(T)$ with $|\la|>r$
	\begin{equation}
		\|Q_\la\| \leq \exp(\gamma |\la|^\tau) = \exp(\gamma |\la|^{p \varrho}),
	\end{equation}
	where $\varrho = \tau/p \leq 1/2$.
	
	Due to the assumption \eqref{la.n.asm} and \cite[Lemma 3]{Agranovich-1995-28}, there exists $c>0$ and a strictly increasing sequence $\{j_k\}_{k \in \N} \subset \N$ such that
	
	\begin{equation}\label{Agr.gaps.0}
		|\la_{j_k+1}| - |\la_{j_{k}}| \geq c |\la_{j_k+1}|^\frac{\kappa-1}{\kappa}.
	\end{equation}
	We select points $z_k$ in $\Im(\Gamma)$ such that
	\begin{equation}
		|z_k| = \frac{|\la_{j_k+1}| + |\la_{j_{k}}|}{2}, \quad k \in \N.
	\end{equation}
	Since
	\begin{equation}
		|z_k| \leq  |\la_{j_{k+1}}|  \leq 2 |z_k|,
	\end{equation}
	we get
	\begin{equation}
		\dist(z_k,\sigma(T)) \geq \frac c2 |\la_{j_k+1}|^\frac{\kappa-1}{\kappa} \geq \frac c2 \min\{1,2^\frac{\kappa-1}{\kappa} \}  |z_k|^\frac{\kappa-1}{\kappa} =: c'_\kappa |z_k|^\frac{\kappa-1}{\kappa}.
	\end{equation}
	Applying Theorem~\ref{thm:T}, we obtain
	\begin{equation}\label{zk.res.ub}
		\begin{aligned}
			\|(z_k-T)^{-1}\| &\leq 
			\frac{C}{\dist(z_k, \sigma(T))} \exp  \left( 
			(1+\delta) \gamma \frac{ \cos (p \varrho \arg z_k)}{\cos \left( 
				(1-\epsilon_0) \pi \varrho \right)} |z_k|^{p\varrho}
			\right)
			\\
			& \leq 
			\exp  \left( 
			(1+\delta) \gamma \frac{ 1}{\cos \left( 
				(1-\epsilon_0) \pi \varrho \right)} |z_k|^\tau - \frac{\kappa-1}{\kappa} \log |z_k| + \log \frac{C}{c'_\kappa}
			\right)
			\\
			& = 
			\exp  \left( 
			(1+\delta) \gamma \frac{ 1}{\cos \left( 
				(1-\epsilon_0) \pi \frac{\tau}{p} \right)} |z_k|^\tau(1+o(1))
			\right), \quad k \to \infty.
		\end{aligned}
	\end{equation}
	However, this leads to a contradiction with \eqref{res.lb.limsup} for $z=z_k$ if
	\begin{equation}
		\gamma < \eta \cos \left( (1-\epsilon_0) \pi \frac{\tau}{p} \right).
	\end{equation}
	
	Finally, if additionally \eqref{cor.T.thet0} is satisfied, then \eqref{res.norm.T.p.theta} can be used and the resulting estimate of $\|(z_k-T)^{-1}\|$, analogous to \eqref{zk.res.ub}, leads to a contradiction with \eqref{res.lb.limsup} for all $\gamma<\gamma_0$ with $\gamma_0$ in \eqref{gam.0.theta}.
\end{proof}

\subsection{Imaginary even anharmonic oscillators}
\label{ssec:H.even}

Let $b\in (0,\infty)$ and consider the imaginary even oscillators
\begin{equation}\label{H.even.def}
	H = \Dt + \ii |x|^b, \quad \Dom(H) = \{ f \in H^2(\R) \, : \, |x|^b f(x) \in L^2(\R) \},
\end{equation}
studied in particular in \cite{Davies-1999-200,Davies-2000-32,Henry-2012-350,Henry-2014-4}. We summarize known properties of $H$.
\begin{proposition}\label{prop:H.even}
Let $H$ be as in \eqref{H.even.def} with $b \in (0,\infty)$. Then
\begin{enumerate}[\upshape (i)]
	\item the numerical range of $H$ satisfies
	\begin{equation}\label{Num.H.even}
	\Num (H) \subset \{ z \in \C \, : \, \Re z \geq 0, \Im z \geq 0\}; 
	\end{equation}
	\item $H$ is m-accretive with compact resolvent, moreover
	\begin{equation}\label{sk.H.even}
		\frs_k((1+H)^{-1}) = \BigO(k^{- \kappa}), \quad k \to \infty,
	\end{equation}
	where
	\begin{equation}\label{kap.H.even}
		\kappa = \frac{2b}{b+2};
	\end{equation}
	
	\item the spectrum of $H$ comprises infinitely many simple eigenvalues $\{\la_n\}_{n\in \N}$ of the form
	\begin{equation}\label{la.n.H.even}
	\la_n = e^{\ii \frac{\pi}{b+2}} \mu_n, \quad n \in \N,
	\end{equation}
	where $\{\mu_n\}_{n \in \N}$ are the eigenvalues of the self-adjoint anharmonic oscillator $G_0$ in \eqref{G0.def} with $l=b$; 
	\item if $b>2/3$, then the system of one-dimensional Riesz projections $P_n \equiv Q_{\la_n}$, $n \in \N$, related to the eigenvalues $\la_n$ of $H$, see	\eqref{Q.la.T.def} with $T=H$, is complete in $L^2(\R)$;
	\item we have as $r \to + \infty$
	\begin{equation}\label{H.even.res.iR}
	\begin{aligned}
	\|(\ii r - H)^{-1}\|  = \BigO \left(r^{-\frac23 \left(1-\frac1b\right)} \right), 
\qquad
	\|(r - H)^{-1}\|  = \BigO \left(r^{-\frac{1}{2} \frac b{b+1}} \right);
	\end{aligned}
	\end{equation}
	\item\label{prop.H.even.lb} for any $\eps \in (0,\pi/4)$, there exist $\eta, R >0$ such that for all $z \in \rho(H)$ with 
	\begin{equation}
		|z| >R \quad \text{and} \quad \arg z \in \left[\eps, \tfrac{\pi}{2}- \eps \right],
	\end{equation}
	we have (with $\kappa$ as in \eqref{kap.H.even})
	\begin{equation}\label{H.res.lb.even}
		\|(z-H)^{-1}\| \geq \exp(\eta |z|^{\frac 1 \kappa}).
	\end{equation}
\end{enumerate}
\end{proposition}
\begin{proof}
\begin{enumerate}[\upshape (i), wide]
	\item The enclosure of the numerical range \eqref{Num.H.even} follows by a straightforward estimate.
	\item For m-accretivity and compactness of resolvent see Subsection~\ref{ssec:1D.Schr}. The estimate \eqref{sk.H.even} follows by \eqref{sn.GH} in Subsection~\ref{ssec:1D.Schr} by comparing $H$ with $G_0$ for $l=b$.
	\item The claim is based on \eqref{mu.n.asym} and complex scaling, see e.g. \cite{Davies-2000-32}, \cite[Chap.~XIII.10]{Reed4}.
	\item By \eqref{sk.H.even}, the resolvent of $H$ belongs to $\cS_{\frac 1\kappa + \eps}$ for any $\eps>0$. Notice that for $b>2/3$
	\begin{equation}\label{kap.b.23}
	\frac{1}{\kappa} = \frac12 + \frac 1b < 2.
	\end{equation}
	By the enclosure of the numerical range \eqref{Num.H.even} and the resolvent bound \eqref{res.Num.T}, we obtain that for every $\omega \in [\pi/2, 3\pi/2]$
	\begin{equation}
	\|(r e^{\ii \omega}-(H + e^{\ii \frac \pi 4}) )^{-1} \| = \BigO(1), \qquad r \to + \infty.
	\end{equation}
	Thus the claim on the completeness follows from Proposition~\ref{prop:compl}.
	\item See \cite[Thm.~3.2, Thm.~4.2, Sec.~5.3, 5.4 and 7.1]{Arnal-2023-284}; cf.~also \cite{Boulton-2002-47,Pravda-Starov-2006-73,BordeauxMontrieux-2013}.
	\item The proof follows the arguments in Appendix~\ref{app:res.L.lb} for $L$ from \eqref{L.def} which are for $H$ technically slightly simpler. For details see~\cite{Dencker-2004-57} or \cite[Thm.1, Sec.~VII]{Krejcirik-2015-56}; cf.~also \cite{Davies-1999-200, Krejcirik-2019-276}.
	\qedhere
\end{enumerate}
\end{proof}

Corollary~\ref{cor:T.Pn} yields the following result on the spectral projections of $H$, cf.~Theorem~\ref{thm:L.Pn} for $L$.

\begin{theorem}\label{thm:H.even.Pn}
Let $H$ be as in \eqref{H.even.def} with $b \in (2/3,\infty)$ and let $P_n$, $n \in \N$, be the one-dimensional Riesz projections related to the eigenvalues $\la_n$ of $H$, see \eqref{la.n.H.even}. Then there exists $\gamma>0$ such that
\begin{equation}\label{Pn.H.even.limsup}
\limsup_{n \to \infty} \frac{\|P_n\| }{\exp (\gamma n)} = \infty.
\end{equation}
\end{theorem}
\begin{proof}
We apply Corollary~\ref{cor:T.Pn} with $p=4$ to the rotated operator and shifted $H$, namely,
\begin{equation}\label{T.H.even.def}
T:= e^{-\ii \frac\pi4} \left(H + 1 \right).
\end{equation}
To this end, we first verify that Assumption~\ref{asm:T} is satisfied. 

By the enclosure of $\Num(H)$ in \eqref{Num.H.even} and the resolvent bound \eqref{res.Num.T}, we get 
\begin{equation}\label{H.even.res}
	\|(e^{ \pm \frac{\ii}{4} \pi} r -T)^{-1}\|  = \BigO (1), 
		\quad  
		\|(e^{\pm \ii \frac{3}{4} \pi} r -T)^{-1}\|  = \BigO (r^{-1}), \quad r \to + \infty.
\end{equation}
Thus \eqref{T.res.pol.asm} is satisfied. The spectrum of $T$ reads
\begin{equation}\label{H.T.even.sp}
\sigma(T) = \{ e^{\ii \frac{2-b}{2+b} \frac \pi 4} \mu_n + e^{-\ii \frac\pi4}\}_{n\in \N},
\end{equation}
thus \eqref{cor.T.thet0} holds with 
\begin{equation}
	\vartheta_0 = \frac{2-b}{2+b} \frac\pi4 \in (-\pi/4,\pi/8).
\end{equation} 
Moreover, \eqref{la.n.asm} holds with $\kappa$ as in \eqref{kap.H.even} and as all eigenvalues of $T$ are simple also Assumption~\ref{asm:T} \ref{asm:T.EV.simple} is satisfied. 
Finally, the system of spectral projections $P_n \equiv Q_{\la_n}$, $n \in \N$, is complete by Proposition~\ref{prop:compl} since $T^{-1} \in \cS_{\frac 1 \kappa+\eps}$ with any $\eps>0$, $1/\kappa <2$, see \eqref{kap.b.23}, and we have the resolvent bounds \eqref{H.even.res}. In summary, Assumption~\ref{asm:T} with $p=4$ is indeed satisfied.

Using \eqref{H.res.lb.even} and the rotation in \eqref{T.H.even.def}, we obtain that \eqref{res.lb.limsup} holds with $ \tau := 1/\kappa<2 = p/2$ for $z>R$, i.e., with $\Gamma(t) =t$, $t>0$. Hence Corollary~\ref{cor:T.Pn} yields
\begin{equation}
	\limsup_{n \to \infty} \frac{\|P_n\|}{\exp(\gamma |\la_n|^\tau)} = \infty
\end{equation}
if $\gamma \in (0,\gamma_0)$ where $\gamma_0 = \eta \cos \left( \tau \vartheta_0 \right)$. The claim \eqref{Pn.H.even.limsup} follows since by \eqref{H.T.even.sp}
\begin{equation}
|\la_n| = |\mu_n|(1+o(1)) = (dn)^\kappa(1+o(1)), \quad n \to \infty,	
\end{equation}
where $d>0$ is as in \eqref{kap.G_0.even} with $l=b$ and $\kappa \tau = \kappa/\kappa=1$.
\end{proof}

We remark that in the proof of Theorem~\ref{thm:H.even.Pn}, we can take $p=2$ for the case $b > 2$ and $p=3$ for $b > 1$. Moreover, if $b>1$, we can use $T:= \exp(- \ii \pi /(b+2))(H+1)$, cf.~\eqref{T.H.even.def}, so that the spectrum of $T$ is real and \eqref{T.res.pol.asm} is satisfied; then $\vartheta_0=0$ and  $\gamma_0$ is improved. We omit details.

\begin{remark}\label{rem:H.even}
Theorem~\ref{thm:H.even.Pn} improves \cite[Thm.~6]{Davies-2000-32}, where it was established that the spectral projection norms cannot be bounded by any power of $n$, and confirms the conjectured exponential growth based on the numerical evidence. The obtained exponential order is sharp for $b \in 2 \N$, where a direct asymptotic analysis of the differential operators shows that for $b=2$ 
\begin{equation}\label{L.shift.Pn.norm}
	\lim_{n \to \infty} \frac{\log \|P_n\|}{n} = \log (1+\sqrt 2),
\end{equation} 
see \cite{Davies-2004-70} and also \cite{Viola-2013-4,Krejcirik-2015-56}, and for $b \in 2 \N$, 
\begin{equation}
 c(b) \equiv	\lim_{n \to \infty} \frac{\log \|P_n\|}{n} >0,
\end{equation} 
see \cite{Henry-2012-350,Henry-2014-4} for more details and a more precise asymptotic expansion. 
For $b \notin 2\N$, our result appears to be new; moreover, the exponential order $1$ is sharp as well, see Subsection~\ref{ssec:JC}. 

\emph{We expect that \eqref{Pn.H.even.limsup} holds with $\lim$ instead of $\limsup$.}
\end{remark}

\subsection{Imaginary odd anharmonic oscillators}
\label{ssec:H.odd}

Let $b \in \N$ and consider the imaginary odd oscillators
\begin{equation}\label{H.odd.def}
	H = \Dt + \ii x^{2b+1}, \quad \Dom(H) = \left\{ f \in H^2(\R) \, : \, |x|^{2b+1} f(x) \in L^2(\R) \right\}.
\end{equation}
The basis spectral properties of $H$ are summarized in the following.

\begin{proposition}\label{prop:H.odd}
Let $H$ be as in \eqref{H.odd.def} with $b \in \N$. Then
	\begin{enumerate}[\upshape (i)]
		\item the numerical range of $H$ satisfies
		\begin{equation}\label{Num.H.odd}
			\Num (H) \subset \{ z \in \C \, : \, \Re z \geq 0\}; 
		\end{equation}
		\item $H$ is m-accretive with compact resolvent, moreover
		\begin{equation}\label{sk.H.odd}
			\frs_k((1+H)^{-1}) = \BigO(k^{- \kappa}), \quad k \to \infty,
		\end{equation}
		where
		\begin{equation}\label{kap.H.odd}
		\kappa=	\frac{4b+2}{2b+3};
		\end{equation}
		\item For any $b=4j-3$, $j \in \N$, i.e., $2b+1 = 8j-5$, there exists $N_0 \in \N$ such that for all $n\geq N_0$ the eigenvalues $\la_n$ of $H$ are simple and real, moreover, they obey the asymptotics 
		\begin{equation}\label{la.n.H.odd}
		\begin{aligned}
		\la_n &= (d n)^\kappa(1+ \BigO(n^{-\frac{1}{2b+1}})), 
		\\
		\la_{n+1}-\la_n & = \kappa d (dn)^{\kappa-1}(1+o(1)), \qquad n \to \infty,
		\end{aligned}
		\end{equation}
		where
		\begin{equation}
		d = \frac{\pi}{B\left(\tfrac12,1+ \tfrac1{2b+1}\right) \cos \left(\tfrac{1}{2b+1} \tfrac \pi 2\right)} >0
		\end{equation}
		(we recall the convention that the eigenvalues are repeated with the algebraic multiplicity and ordered with a non-decreasing modulus);
		\item the system of Riesz projections $Q_{\la}$, $\la \in \sigma(H)$,
		is complete in $L^2(\R)$;
		\item we have as $r \to + \infty$
		\begin{equation}\label{H.odd.res.iR}
		\|(\pm \ii r - H)^{-1}\|  = \BigO \left( r^{-\frac{4b}{3(2b+1)}} \right), 
		\end{equation}
		\item for any $\eps \in (0,\pi/2)$ there exist $\eta, R >0$ such that for all $z \in \rho(H)$ with 
		\begin{equation}
			|z| >R \quad \text{and} \quad |\arg z| \leq \frac{\pi}{2} - \eps,
		\end{equation}
		we have (with $\kappa$ as in \eqref{kap.H.odd})
		\begin{equation}\label{H.res.lb.odd}
			\|(z-H)^{-1}\| \geq \exp(\eta |z|^{\frac 1 \kappa}).
		\end{equation}
	\end{enumerate}
\end{proposition}
\begin{proof}
\begin{enumerate}[\upshape (i), wide]
	\item The enclosure of the numerical range \eqref{Num.H.odd} follows by a straightforward estimate.
	\item For m-accretivity and compactness of resolvent see Subsection~\ref{ssec:1D.Schr}.  The estimate \eqref{sk.H.odd} follows by \eqref{sn.GH} in Subsection~\ref{ssec:1D.Schr} by comparing $H$ with $G_0$ for $l=2b+1$.
	\item The claims follow from \cite[Cor.~2.1, 2.2]{Shin-2010-35} with $m = 2b+1 = 8j -5$, $l = 2j-1$ and \cite[Thm.~1.2]{Shin-2010-35} since $H$ is $\cP\cT$-symmetric.
	\item By \eqref{sk.H.odd} and \eqref{kap.H.odd}, the resolvent of $H$ belongs to $\cS_{\frac 1 \kappa + \eps}$ for any $\eps>0$. Since $b \geq 1$, we have
	\begin{equation}\label{kap.b.odd}
		\frac{1}{\kappa} = \frac12 + \frac 1{2b+1} < 1. 
	\end{equation}
	By the enclosure of the numerical range \eqref{Num.H.odd} and the resolvent bound \eqref{res.Num.T}, we obtain that for every $\omega \in [\pi/2, 3\pi/2]$
	\begin{equation}
		\|(r e^{\ii \omega}-(1+H)^{-1} \| = \BigO(1), \qquad r \to + \infty.
	\end{equation}
	The claim on the completeness follows by Proposition~\ref{prop:compl}, see Subsection~\ref{ssec:ot}.
	\item See \cite[Thm.~3.2, Sec.~7.1]{Arnal-2023-284}; cf.~also \cite{BordeauxMontrieux-2013}.
	\item See~\cite{Dencker-2004-57} or \cite[Thm.1, Sec.~VII]{Krejcirik-2015-56}; the proof follows the arguments in Appendix~\ref{app:res.L.lb} for $L$ from \eqref{L.def}.
	\qedhere
\end{enumerate}
\end{proof}

Corollary~\ref{cor:T.Pn} yields the following claim on the spectral projections of $H$. Notice that the restriction $b \in 4\N -3$ is used to express $\la_n$ as in \eqref{la.n.H.odd}. In fact, Corollary~\ref{cor:T.Pn} can be used with any $b \in \N$, however, $\|Q_{\la_n}\|$ is then compared with $\exp(\gamma|\la_n|^{1/\kappa})$ only and the behavior of $\la_n$ is not known for a general $b \in \N$.

\begin{theorem}\label{thm:H.odd.Pn}
Let $H$ be as in \eqref{H.odd.def} with $b=4j-3$, $j \in \N$, and let $P_n$ be the one-dimensional Riesz projections related to the eigenvalues $\la_n$, $n > N_0$. Then there exists $\gamma>0$ such that
\begin{equation}\label{Pn.H.odd.limsup}
\limsup_{n \to \infty} \frac{\|P_n\|}{\exp (\gamma n)} = \infty.
\end{equation}
\end{theorem}
\begin{proof}
The proof is analogous to the one of Theorem~\ref{thm:H.even.Pn}. Here we take $p=2$ and $T=H+1$ since $1/\kappa<1$, see \eqref{kap.b.odd}, and for all $r>0$
\begin{equation}
\|(\ii r - T)^{-1}\| \leq 1 
\end{equation}
by \eqref{Num.H.odd} and \eqref{res.Num.T}. We omit further details.
\end{proof}

\begin{remark}\label{rem:H.odd}
In the special case $b=1$, it was showed  by the direct analysis that
\begin{equation}
	\lim_{n \to \infty} \frac{\log \|P_n\|}{n} = \frac{\pi}{\sqrt3},
\end{equation}
see in \cite{Henry-2014-15}, thus the exponential order $1$ in \eqref{Pn.H.odd.limsup} is sharp in this case and the claim \eqref{Pn.H.odd.limsup} holds in fact with $\lim$ instead of $\limsup$. Moreover, like in Remark~\ref{rem:H.even}, the exponential order $1$ in \eqref{Pn.H.odd.limsup} is sharp in all cases, see Subsection~\ref{ssec:JC}.

\emph{We expect that \eqref{Pn.H.odd.limsup} holds with $\lim$ in all cases in Theorem~\ref{thm:H.odd.Pn} and also for other $b \in \N$.} 

\end{remark}

\subsection{Conjugated real anharmonic oscillators}
\label{ssec:MSV.op}

Let $G_0$ be as in \eqref{G0.def} with $l = b \geq 2$ and let $v \in C^2(\R)$ be a non-negative odd function which is increasing on $(0,\infty)$ and for some $s \in (0,1)$
\begin{equation}\label{v.conj.def}
v(x) = \frac12 x^{\frac{2+b}{2} s}, \quad x>1.
\end{equation}
With this $v$, let $H$ be the following conjugation of $G_0$
\begin{equation}\label{H.conj.def}
	H = e^v G_0 e^{-v} = \left(-\ii \Do + \ii v'\right)^2 + |x|^b , \qquad \Dom(H) = \Dom(G_0)
\end{equation}
in $L^2(\R)$. Notice that $e^v$ is not a bounded operator in $L^2(\R)$ and so it does not induce a (bounded) similarity transform of $H$ and $G_0$. Nonetheless, some properties of $G_0$ are preserved by this unbounded conjugation. In particular, the following was established in \cite{Mityagin-2017-272}.

\begin{proposition}\label{prop:H.conj}
Let $H$ be as in \eqref{H.conj.def} with fixed $b \geq 2$ and $s \in (0,1)$. Let $\mu_n$ and $\kappa$ be as in \eqref{mu.n.asym}, \eqref{kap.G_0.even} with $l=b$, in particular
\begin{equation}\label{kap.H.conj}
	\kappa = \frac{2b}{b+2}. 
\end{equation}
Then
\begin{enumerate}[\upshape (i)]
	\item there exist $C, \alpha_0 > 0$ such that
	\begin{equation}
	\Num(H+\alpha_0) \subset \{ z \in \C \, : \, \Re z \geq 1, |\Im z| \leq C (\Re z)^\beta\}
	\end{equation}
	with 
\begin{equation}
	\beta = \frac 12 + \max \left\{ 0, \frac{s}{2}+ \frac{s-1}{b} \right \} < 1;
\end{equation}
\item $H$ is m-sectorial with compact resolvent, $0 \in \rho(H)$ and 
\begin{equation}
\frs_k(H^{-1}) = \BigO(k^{- \kappa}), \quad k \to \infty;
\end{equation}
\item the spectrum of $H$ and $G_0$ coincide, i.e., it comprises infinitely many simple eigenvalues $\la_n= \mu_n$, $n \in \N$;
\item the system of one-dimensional Riesz projections $P_n \equiv Q_{\la_n}$, related to the eigenvalues $\la_n=\mu_n$, $n \in \N$, is complete in $L^2(\R)$;
\item the norms of $P_n$ satisfy (with $\kappa$ as in \eqref{kap.H.conj})
\begin{equation}\label{Pn.conj.lim}
\|P_n\| = \exp \left(\la_n^\frac{s}{\kappa} (1+o(1)) \right), \quad n \to \infty.
\end{equation}
\end{enumerate}
\end{proposition}

Since the behavior of $\|P_n\|$ is known, we can apply Theorem~\ref{thm:T} or Corollary~\ref{cor:T.par} to obtain an upper bound on the resolvent norm of $H$ (cf.~Theorems~\ref{thm:L.Pn}, \ref{thm:H.even.Pn}, \ref{thm:H.odd.Pn} where a lower bound on the resolvent norm was used to obtain a result on the norms of $P_n$).

\begin{theorem}\label{thm:H.conj}
Let $H$ be as in \eqref{H.conj.def} with fixed $b \geq 2$ and $s \in (0,1)$. Then for every $\delta>0$, there exist $C \equiv C(\delta)>0$ and $r \equiv r(\delta)>0$ such that for all $z \in \rho(H)$ satisfying $|z|> r$ and $|\arg z| \leq \pi/4$,
\begin{equation}\label{H.conj.res.est}
\|(z-H)^{-1}\| \leq \frac{C}{\dist(z,\sigma(H))} \exp \left( (1+\delta) \cos \left(\frac s \kappa \arg z\right) |z|^\frac s \kappa \right),
\end{equation}
where $\kappa$ is as in \eqref{kap.H.conj}.
\end{theorem}
\begin{proof}
Since $\Num(H)$ is contained in a ``parabolic'' region, see Subsection~\ref{ssec:ot}, the claim follows by Corollary~\ref{cor:T.par} with $\epsilon_1 = 1/2$, $p=2$ and 
\begin{equation}
\varrho = \frac{s}{2 \kappa} = \frac{s}{2} \left( \frac 12 + \frac 1b  \right) \leq \frac{s}{2}< \frac 12,
\end{equation}
see Proposition~\ref{prop:H.conj}.
\end{proof}

\begin{remark}
In the special case $b=2$ and $v(x)=x/2$ (so $\kappa=1$ and $s =1/2$),  the operator $H$ reads 
\begin{equation}\label{H.conj.shift}
H = \left(-\ii \Do + \frac\ii 2 \right)^2 + x^2. 
\end{equation}
It is unitarily equivalent (via Fourier transform) to the shifted harmonic oscillator
\begin{equation}
- \frac{\dd^2 }{\dd \xi^2}+ \left(\xi + \frac \ii 2 \right)^2,
\end{equation}
i.e., to~$L$ in \eqref{L.shift.HO.def} with $\alpha_1 = \ii$ and $\alpha_0 = -1/4$. Hence \eqref{L.res.lb.odd} shows that the exponential order $s/\kappa=1/2$ in \eqref{H.conj.res.est} is sharp in this case (see also Remark~\ref{rem:L.shifted.HO}). 

\emph{We expect that the exponential order $s/\kappa$ in \eqref{H.conj.res.est} is sharp also in all other cases in Theorem~\ref{thm:H.conj}.}
\end{remark}

\section{Further remarks}
\label{sec:misc}

We collect several remarks on our results and related conjectures.

\subsection{Remarks on Theorem~\ref{thm:T}}
\label{ssec:extensions}

As in \cite{Davies-2000-32}, one can formulate Theorem~\ref{thm:T} in a Banach space; no changes in the proof are needed.

The condition on the equality of algebraic and geometric multiplicities is satisfied e.g.~for (multi-dimensional) differential operators obtained from self-adjoint operators by transformations which preserve the eigenvalues and their multiplicities like complex scaling and complex shift of the variable or a conjugation as in \eqref{H.conj.def}.

Remark \ref{rem:ak.modification} indicates a flexibility in the choice of the gauge function $F$. Besides possible modifications of zeros of $F$ by additive perturbations and scaling, we also used the substitution $w = y^p$ in Section~\ref{sec:T}, see~e.g.~\eqref{pf.decomp.p}. 

\subsection{Sharpness of the exponential order the estimate of $\|P_n\|$}
\label{ssec:JC}

The exponential orders in \eqref{Pn.L.limsup}, \eqref{Pn.H.even.limsup} and \eqref{Pn.H.odd.limsup} are sharp for the special cases of $a$ and $b$, which were analyzed directly in earlier works, see Remarks~\ref{rem:L.shifted.HO}, \ref{rem:H.even} and \ref{rem:H.odd}. Here we discuss the sharpness of the exponential order for $H$ in \eqref{H.even.def} and \eqref{H.odd.def} and for $L$ in \eqref{L.def} for other $a,b$.

First, for $H$ in \eqref{H.even.def} or \eqref{H.odd.def}, it follows from \cite{CueSie-2026} that 
\begin{equation}\label{CuSi.bound}
	\|(z-H)^{-1}\| \leq \frac{\exp(\gamma (1+|z|^\frac1\kappa)
		)}{\dist(z,\sigma(H))}, \quad z \in \rho(H), \ |z| \geq R,
\end{equation}
where $R>0$ is sufficiently large, $\gamma>0$ depends only on $b$ and $\kappa$ is as in \eqref{kap.H.even} or \eqref{kap.H.odd}. Thus with $\delta_n := \eps n^{\kappa-1}$ and $\eps>0$ sufficiently small, see Proposition~\ref{prop:H.even} or \ref{prop:H.odd}, we get that there exist $\gamma' >\gamma$ and $\gamma'' > \gamma$ such that for all sufficiently large $n$ (so that \eqref{CuSi.bound} can be used in the $\delta_{n}$-disc around $\la_n$)
\begin{equation}\label{H.Pn.sharp}
	\begin{aligned}
		\|P_n\| &\leq \frac1{2\pi} \int_{|z-\la_n| = \delta_{n}} \|(z-H)^{-1}\| \, |\dd z| 
		 \leq \exp \left(\gamma \left(1+ \sup_{|z-\la_n|=\delta_{n}}|z|^\frac1\kappa \right) \right) 
		\\
		& \leq \exp (\gamma'|\la_n|^\frac1\kappa ) \leq \exp (\gamma'' n).
	\end{aligned}
\end{equation}
This justifies that exponential order $1$ is sharp in \eqref{Pn.H.even.limsup} and \eqref{Pn.H.odd.limsup}.

For $L$, the bound like \eqref{CuSi.bound} is valid with $\kappa$ as in \eqref{kap.L.def}. As above, an analogue of \eqref{H.Pn.sharp} together with \eqref{la.L.intro} show that $\|P_n\| \leq \exp( \gamma n)$, $n \in \N$, for some $ \gamma >0$. However, this does not yield the sharpness of the exponential order $\sigma<1$ in \eqref{Pn.L.limsup}. Notice that $\kappa$ is determined by the behavior of the singular values of $L^{-1}$ which in turn stems from the leading term $x^{2a}$ of the potential, irrespective of its imaginary part.

\subsection{A subsequence in the estimate of $\|P_n\|$}
\emph{We expect that \eqref{Pn.L.limsup}, \eqref{Pn.H.even.limsup} and \eqref{Pn.H.odd.limsup} hold with $\lim$ instead of $\limsup$.} It happens in the special cases discussed in Remarks~\ref{rem:L.shifted.HO}, \ref{rem:H.even} and \ref{rem:H.odd}.

\subsection{Absence of basis properties for $L$ with non-integer $a,b$}
\label{ssec:no.basis}

In Theorem~\ref{thm:L.Pn} we explained that the eigensystem of $L$ does not contain a basis as the spectral projections are not uniformly bounded. The essential condition on the real and imaginary part of the potential is $b>a-1$; for $b<a-1$, the eigensystem of $L$ contains a Riesz basis, see \cite[Thm.~6.4]{Mityagin-2019-139}. In the borderline case $b=a-1$ it is known that the eigensystem of $L$ contains a Riesz basis of finite dimensional subspaces, see \cite[Theorem, case (b)]{Agranovich-1995-28}. In the case of the operator $L =  -(\dd/\dd x)^2 + x^2 + q(x)$, $q \in L^\infty (\R)$, in $L^2(\R)$, 
A.~Shkalikov's results and constructions \cite{Shkalikov-2016-71,Shkalikov-1983-268,Shkalikov-1983} could explain the possibility of choosing Riesz basis system of spectral block-projections $\{P_j\}$ with uniformly bounded ranks; moreover, $\Rank (P_j) \leq 2 \|q\|_{L^\infty}$; \cite{Shkalikov-2019-private}.

We expect that $b=a-1$ is the borderline case also for non-integer $a$, $b$. As a concrete example, consider 
$L = -(\dd/\dd x)^2 + V$ in $L^2(\R)$ with a potential $V \in C^\infty(\R)$ that satisfies
\begin{equation}
	V(x) = |x|^{2a} + \ii \sgn(x) |x|^b \quad \text{for } |x|>1
\end{equation}
and assume that $a \geq 1$, $b \in \R$.

First, if $ b < a-1$, then the eigensystem of $L$ contains a Riesz basis, see \cite[Thm.~6.4]{Mityagin-2019-139}. We note that for $W \in L^p(\R)$, $p \in (0,\infty)$, also the eigensystem of $L+W$ has the Riesz basis property; this holds even for some distributional $W$, see \cite{Mityagin-2019-139} for details. 

Second, assume that $b>a-1$. It follows from \cite[Thm.~3.7, Ex.~3.8]{Krejcirik-2019-276} that for any fixed $\beta \in \R$ and $N \in \N$, there exists $C>0$ such that for all $\alpha + \ii \beta \in \rho(L)$ satisfying $\alpha \geq 1$, we have
\begin{equation}\label{res.lb.gen}
\|(\alpha+\ii \beta - L)^{-1}\| \geq C \alpha^N. 
\end{equation}
We note that an analogous claim holds also for less regular $V$ and consequently possibly slower rates of the resolvent growth, see \cite[Thm.~3.7, Ex.3.8]{Krejcirik-2019-276}; for estimates on general curves in $\C$ see \cite[Sec.~5]{Krejcirik-2019-276}. Recall also that for polynomial potentials, Theorem~\ref{thm:L.res.lb} yields exponential lower bounds for the resolvent norm.
 
Since $\Num(L + t) \subset \{ z \in \C \, : \, \Re z \geq 0\}$ for a sufficiently large $t$ and $\frs_k((L+t+1)^{-1}) = \BigO(k^{-\kappa})$ as $k \to +\infty$ with 
\begin{equation}
\kappa = \frac{2m}{m+2}, \quad m:=\max \{2a,b\} \geq 2,
\end{equation}
see Subsection~\ref{ssec:1D.Schr}, it follows from Proposition~\ref{prop:compl} that the spectrum of $L$ comprises infinitely many eigenvalues $\{\la_n\}_{n \in \N}$ and the system of spectral projections $Q_\la$, $\la \in \sigma(L)$, cf.~\eqref{Q.la.T.def} with $T=L$, is complete. We note also that by \eqref{lak.sk.ineq} we get $|\la_n| \geq c n^\kappa$ for some $c>0$ and all $n>n_0$ with a sufficiently large $n_0 \in \N$. Moreover,  \cite[Lemma~3]{Agranovich-1995-28} shows that there exists $c'>0$ and a strictly increasing sequence $\{j_k\}_{k \in \N} \subset \N$ such that
\begin{equation}\label{Agr.gaps}
	|\la_{j_k+1}| - |\la_{j_{k}}| \geq c' j_k^{\kappa-1}.
\end{equation}

We show below that the eigensystem of $L$ does not contain a basis \emph{if the eigenvalues $\la_n$ are eventually simple.} (The latter is not known for non-polynomial potentials.) We argue by contradiction and assume that the eigensystem of $L$ contains a basis. Following the steps in the proof of Proposition~\ref{prop:L:no.basis}, we get that
\begin{equation}\label{L.bas.is.gen}
\|(z-L)^{-1} \| 
\leq 
\frac{C|z|}{\dist(z, \sigma(L))}. 
\end{equation}
However, for the points $z_k =(|\la_{j_k}| + |\la_{j_k+1}|)/2$, the bounds \eqref{res.lb.gen} and \eqref{L.bas.is.gen} lead to a contradiction as $k \to \infty$. 

Finally, the eigensystem of $L$ contains a Riesz basis of finite dimensional subspaces due to \cite[Theorem, case (b)]{Agranovich-1995-28} if $b=a-1$ and an Abel basis of finite dimensional invariant subspaces of the order $\tau + \eps$ with $\tau$ from \eqref{rho.L.lb} and an arbitrarily small $\eps>0$ if $a-1<b<2a$, see \cite[Theorem, case (c)]{Agranovich-1995-28} where $p=\kappa$ and $q=b/(2a)$.

\subsection{Critical case in Theorem~\ref{thm:H.even.Pn}}
\label{ssec:b23}

For $b=2/3$, the completeness of the eigensystem of $H$ is non-trivial (a borderline case of the application of Proposition~\ref{prop:compl}). It is established for the operators $H = -(\dd/\dd x)^2 + \ii |x|^\frac23$ in $L^2(\R_+)$ with Dirichlet boundary condition at $0$, see \cite{Tumanov-2021-280,Tumanov-2022-319}, and the given proof can be extended to Neumann boundary conditions and also for $H$ in $L^2(\R)$. 

For $b=2/3$, we have $1/\kappa = 2$, cf.~\eqref{kap.b.23}. Notice that $\tau=1/\kappa=2$ is expected in \eqref{cor.T.Pn} due to \eqref{CuSi.bound} and \eqref{H.Pn.sharp}. Thus one can attempt to apply Corollary~\ref{cor:T.Pn} with $p=4$ and $\tau =p/2=2$. In this case, the additional condition $N < -1 +1/(2p) =-7/8$ in \eqref{T.res.pol.asm} needs to be verified (the remaining conditions in Assumption~\ref{asm:T} are satisfied). However, neither the rotation and shift like in \eqref{T.H.even.def} nor the (sharp) rates in \eqref{H.even.res.iR} can be used to show that \eqref{T.res.pol.asm} holds with such $N$. 

Nonetheless, Corollary~\ref{cor:T.Pn} is applicable for any $\tau <2$ and hence also an analogue of the statement of Theorem~\ref{thm:H.even.Pn} holds for $b=2/3$, namely 
\begin{equation}\label{H.23.PN.eps}
	\limsup_{n \to \infty} \frac{\|P_n\|}{\exp(n^{1-\eps})} = \infty
\end{equation}
for any $\eps>0$.

\subsection{The case $b=0$ in Theorem~\ref{thm:H.odd.Pn}}

The case $b=0$ in Theorem~\ref{thm:H.odd.Pn}, i.e., ~$H = -(\dd/\dd x)^2 + \ii x$ is not included in Theorem~\ref{thm:H.odd.Pn} (nor in Proposition~\ref{prop:H.odd}). The spectrum of this $H$ is empty, so there are no non-trivial Riesz projections for this $H$.

\appendix

\section{Partial fraction decomposition} 
\label{app:pfc}

Decomposition of meromorphic functions into simple fractions starting from Mittag-Leffler \cite{Mittag-Leffler-1884-4}, \cite[Chap.~7]{Markushevich-1977}, has a long history. Papers of M.~Krein \cite{Krein-1944-44,Krein-1947-11} and H.~Hamburger \cite{Hamburger-1944-66} linked to and inspired by the moment theory were the initial milestone in
the after-WW2 development of this topic, cf.~\cite{Akhiezer-1965,Bakan-1998,Bakan-2009,Berg-1995-171,Borichev-2001-45,DeBranges-1959-10,Koosis-1980,Loya-2017,Maergoiz-2002-9,Maergoiz-2016-8,Ostrovskii-1976-229,Sherstyukov-2011-202,Sherstyukov-2016,Sherstyukov-2021-257}. 
We refer the reader who would be interested in general theory to monographs and articles mentioned above. But in our paper and this appendix we focus our attention on the narrow question about the behavior of functions $F$ and $1/F$, where $F(w)$ is defined in (2.28), and the partial fraction decomposition of $1/F$. For convenience of a reader (and our own -- 
in the proofs in Sections \ref{sec:L} and \ref{sec:T}) we give in this appendix a few facts and statements which are well known in literature or in folklore.

\subsection{Proposition~\ref{prop:F.Titch} and its refinement}

We follow \cite[Thm.~I, p.~199]{Titchmarsh-1927-26} or \cite[Chap.~8.64]{Titchmarsh-1958}.
\begin{proof}[Sketch of the proof of Proposition~\ref{prop:F.Titch}]
First let $x \in \R$. Using the partial summation in the second step, we obtain
\begin{equation}
\begin{aligned}
\log F(x) &= \sum_{n=1}^\infty \log \left(1+\frac{x}{a_n}\right)
=
\sum_{n=1}^\infty n \left[\log \left(1+\frac{x}{a_n}\right)-\log \left(1+\frac{x}{a_{n+1}}\right)\right]
\\
& = 
\sum_{n=1}^\infty n \int_{a_n}^{a_{n+1}} \frac{x \, \dd r}{r(x+r)}.
\end{aligned}
\end{equation}
Next, by the definition of $n(r)$ in \eqref{n.r.asm}, we have $n(r) = 0$ if $0 \leq r < a_1$ and we arrive at
\begin{equation}\label{F.int}
\log F(x) = x \sum_{n=1}^\infty \int_{a_n}^{a_{n+1}} \frac{ n(r) \, \dd r}{r(x+r)} 
=
x \int_{0}^{\infty} \frac{ n(r) \, \dd r}{r(x+r)}.
\end{equation}
Thus for the main term (see \eqref{n.r.asm}), we have
\begin{equation}
\nu x \int_{0}^{\infty} \frac{ r^\varrho \, \dd r}{r(x+r)}
= 
\nu x^\varrho \int_{0}^{\infty} \frac{ t^{\varrho-1} \, \dd t}{1+t}
= 
\nu \frac{\pi}{\sin(\pi \varrho)} x^\varrho
\end{equation}
and with the integrand being positive we get
\begin{equation}
	\log F(x) = \nu \frac{\pi}{\sin(\pi \varrho)} x^\varrho(1+o(1)), \quad x \to \infty.
\end{equation}

The integral representation of $\log F$ in \eqref{F.int} guarantees that
\begin{equation}\label{F.int.C}
\log F(w) = w \int_{0}^{\infty} \frac{ n(r) \, \dd r}{r(w+r)}, 
\quad w \in G = \C \setminus (-\infty,-a_1]
\end{equation}
because both functions of $w$ are analytic in $G$ and coincide on the positive semi axis $[ 0, + \infty)$.
Using Cauchy formula, we turn the line of integration to $r = u \exp(\ii \arg w)$ and obtain
\begin{equation}
\nu w \int_{0}^{\infty} \frac{ r^\varrho \, \dd r}{r(w+r)}
= 
\nu |w| \exp(\ii \varrho \arg w) \int_{0}^{\infty} \frac{ u^\varrho \, \dd u}{u(|w|+u)} 
= 
\nu \frac{\pi}{\sin(\pi \varrho)} w^\varrho,
\end{equation}
which leads to \eqref{F.asym.Titch}.
\end{proof}

\begin{remark}
For the special $a_k = k^\frac1\varrho$, $k \in \N$, one can refine \eqref{F.asym.Titch} to
\begin{equation}\label{F.refined}
	F(x) = \exp \left(
	\frac{\pi}{\sin(\pi \varrho)} x^\varrho
	\right) \psi(x), \quad x>0,
\end{equation}
where 
\begin{equation}
\frac{\exp\left(-\frac 1\varrho \right)}{1+x}	\leq \psi(x) \leq \exp \left(-\frac 1\varrho \frac{x}{1+x} \right), \quad x>0.
\end{equation}
Indeed, in this case $n(r) = 0$ for $r<1$ and $n(r) = \lfloor r^\varrho \rfloor$, so 
\begin{equation}
r^\varrho -1 \leq n(r) \leq r^\varrho, \quad r \geq 1,
\end{equation}
and 
\begin{equation}
- x \int_{1}^{\infty} \frac{ \dd r}{r(x+r)} - x \int_{0}^{1} \frac{r^\varrho \, \dd r}{r(x+r)} 	
\leq  
x \int_{0}^{\infty} \frac{ n(r) \, \dd r}{r(x+r)} - \frac{\pi x^\varrho}{\sin(\pi \varrho)}  
\leq 
- x \int_{0}^{1} \frac{r^\varrho \, \dd r}{r(x+r)}. 
\end{equation}
Finally, we use
\begin{equation}
\begin{aligned}
x \int_{1}^{\infty} \frac{ \dd r}{r(x+r)} &= \log(1+x),
\\	
\frac{x}{x+1} \frac1 \varrho \leq  x \int_{0}^{1} \frac{r^\varrho \, \dd r}{r(x+r)} & \leq \frac1 \varrho,
\end{aligned}	
\end{equation}
which leads to \eqref{F.refined}; cf.~\eqref{F.int}.
\end{remark}

\subsection{Partial fraction decompositions}

First, we give a series of elementary identities and inequalities which make Proposition~\ref{prop:B-P} below applicable in our setting.

For any $w \in \C \setminus \Z$ define a periodic function on $\R$
\begin{equation}
\begin{aligned}
h(t) & = \cos w t, \quad t \in [-\pi,\pi],
\\
h(t+2k\pi) & = h(t), \quad k \in \Z.	
\end{aligned}
\end{equation}
Its Fourier series
\begin{equation}
h(t) = \frac{\sin\pi w}{\pi w} + \sum_{k=1}^\infty (-1)^k \frac{2 w \sin \pi w}{\pi(w^2-k^2)} \cos kt
\end{equation}
converges uniformly on $\R$ by Dini-Lipschitz criterion, so $t=0$ leads to the identity
\begin{equation}\label{w/sin.exp}
\frac{\pi w}{\sin \pi w} = 1 + 2 \sum_{k=1}^\infty (-1)^k \frac{w^2}{w^2-k^2}, \qquad w \in \C \setminus \Z,
\end{equation}
see also \cite[Eq.~4.22.5]{DLMF}, and $t=\pi$ gives
\begin{equation}\label{cot.exp}
\cot \pi w = \frac{1}{\pi w} + \sum_{k=1}^\infty \frac{2 w}{\pi} \frac{1}{w^2-k^2},  \qquad w \in \C \setminus \Z,
\end{equation}
see also \cite[Eq.~4.22.3]{DLMF}.

Later, we need to know estimates or asymptotics for the derivatives $F'(-a_n)$ of the products
\begin{equation}
	F(w) = \prod_{k=1}^\infty \left( 1 + \frac{w}{a_k} \right).
\end{equation}
In advance, let us analyze the sequence
\begin{equation}\label{An.b.def}
A(n ; b) = (-1)^{n-1} \prod_{\substack{k=1 \\ k\neq n}}^\infty \left( 1 - \frac{n^b}{k^b} \right), \quad 1 < b < \infty.
\end{equation}

\subsubsection{Special case $b=2$}
If we use Euler product, 
\begin{equation}\label{Euler.prod}
\frac{\sin \pi w}{\pi w} = \prod_{k=1}^\infty \left(1 - \frac{w^2}{k^2}  \right),
\end{equation}
we find that
\begin{equation}
\begin{aligned}
A(n;2) & = (-1)^{n-1} \lim_{w \to n} \frac{\sin \pi w - \sin \pi n}{\pi w \left(1 - \frac{w^2}{n^2}\right) } 
= (-1)^{n} \lim_{w \to n} \frac{\sin \pi w - \sin \pi n}{\pi (w-n)} \frac{n^2}{w(n+w)}
\\
& 
= (-1)^{n} \cos(\pi n) \cdot \frac 12 = \frac 12. 
\end{aligned}	
\end{equation}
But let us do an elementary straightforward evaluation avoiding Euler product \eqref{Euler.prod}
\begin{equation}\label{An.2.id}
\begin{aligned}
A(n;2) & = (-1)^{n-1} \prod_{\substack{k=1 \\ k\neq n}}^\infty \left( 1 - \frac{n^2}{k^2}  \right)
= 
(-1)^{n-1} \lim_{m \to \infty} \prod_{\substack{k=1 \\ k\neq n}}^m \frac{(k-n)(k+n)}{k^2}
\\
& =
\lim_{m \to \infty}  
\left(
(n-1)! (m-n)!
\right)
\left(
\frac{(m+n)!}{n! 2n}
\right)
\left(
\frac{n.n}{m!m!}
\right)
\\
& =
\lim_{m \to \infty}  
\frac12 \frac{(m+1)\dots (m+n)}{m(m-1)\dots(m-n+1)} = \frac12.
\end{aligned}
\end{equation}

\subsubsection{Estimates for case $b\neq 2$}
We basically follow \cite[Ex.~5.5]{Berg-1995-171}. By \eqref{An.2.id},
\begin{equation}\label{1.An.2.fact}
1 = 2 A(n;2) = \frac 1{2 A(n;2)} = \frac{1}{2} (-1)^{n-1} \prod_{\substack{k=1 \\ k\neq n}}^\infty \frac{k^2}{k^2-n^2}.
\end{equation}
To analyze $A(n;b)$, see \eqref{An.b.def}, we ``improve''
\footnote{This trick has been suggested by Leif Mejlbro and used in \cite{Berg-1995-171} to get estimates \eqref{An.b.est} below.}
the sequence of factors in the r.h.s.~of \eqref{An.b.def} by multiplying them by corresponding factors in factorization of $1$ in \eqref{1.An.2.fact}. We get
\begin{equation}\label{An.b.GH}
	A(n;b) = \frac12 \prod_{\substack{k=1 \\ k\neq n}}^\infty \frac{k^b-n^b}{k^2-n^2} \frac{k^2}{k^b} = \frac12 G \cdot H,
\end{equation}
where
\begin{equation}\label{GH.def}
G = \prod_{k=1}^{n-1} r \left(\frac kn\right), \quad H = \prod_{k=n+1}^\infty  r \left(\frac kn\right),
\end{equation}
and
\begin{equation}\label{r.def}
r(t) = \frac{1-v^\frac b2}{1-v}, \quad 0 < t < \infty; \quad v = \frac{1}{t^2}.
\end{equation}
The mean value theorem shows that if $b>2$, then
\begin{equation}
r(t) > 
\begin{cases}
1, & t>1,
\\
\frac{b}{2}, & 0<t<1,
\end{cases}
\end{equation}
and if $b<2$, then
\begin{equation}
	r(t) < 
	\begin{cases}
		1, & t>1,
		\\
		\frac{b}{2}, & 0<t<1.
	\end{cases}
\end{equation}

Therefore, by \eqref{An.b.GH} and \eqref{GH.def}, 

\begin{equation}\label{An.b.est}
A(n;b) > \frac 12 \left( \frac{b}{2}\right)^{n-1} = \frac{1}{b} \left( \frac{b}{2}\right)^n \quad \text{if} \ b>2,
\end{equation}
and
\begin{equation}\label{An.b.est.2}
A(n;b) < \frac 1b \left( \frac{b}{2}\right)^{n} \quad \text{if} \ b<2.
\end{equation}

\subsubsection{Asymptotics for case $b\neq 2$}

Estimates \eqref{An.b.est} are good enough to guarantee convergence of the partial fraction decomposition when we need it. But to complete the picture let us show the validity of asymptotics, for $1<b<\infty$,
\begin{equation}\label{An.b.asym}
A(n;b) = \exp \left(
n \pi \cot \frac \pi b + o(n)
\right), \quad n \to \infty.
\end{equation}
The proof has been given by A.~Sedletski in \cite[Lemma~5]{Sedletskii-1983-34}; this asymptotics was used by E.~Titchmarsh in \cite{Titchmarsh-1927-26}.

\begin{proof}[Sketch of a proof of \eqref{An.b.asym}]
By \eqref{An.b.def},
\begin{equation}
\log A(n;b) = \log G_0 + \log G_1,
\end{equation}
where
\begin{equation}
G_0 = \prod_{k=1}^{n-1} \left( \frac{n^b}{k^b} -1 \right), \qquad G_1 = \prod_{k=n+1}^{\infty} \left(1- \frac{n^b}{k^b} \right).
\end{equation}
With $t_k = k/n$
\begin{equation}\label{H01.def}
H_0 = \frac{1}{n} G_0 = \sum_{k=1}^{n-1} \log \left(  t_k^{-b} -1 \right) \cdot \frac1n,
\quad
H_1 = \frac{1}{n} G_1 = \sum_{k=n+1}^{\infty} \log \left(1-  t_k^{-b} \right) \cdot \frac1n.
\end{equation}
Sums \eqref{H01.def} are Riemann sums of the integrals
\begin{equation}
\mathcal J_0 = \int_0^1 g(t) \, \dd t, \qquad 	\mathcal J_1 = \int_1^\infty g(t) \, \dd t,
\end{equation}
where 
\begin{equation}
	g(t) = \log \left| t^{-b}-1  \right|, \quad 0<t<\infty, \quad t \neq 1.
\end{equation}
We evaluate $\cJ_0$ and $\cJ_1$ by using N.~Mercator decomposition 
\begin{equation}
\log(1-v) = - \sum_{k=1}^\infty \frac{v^k}{k}, \quad 0 \leq v <1,
\end{equation}
and Lebesgue dominated convergence theorem.
We have
\begin{equation}\label{J0.def}
\cJ_0 = \int_0^1 \log \left( \frac{1}{t^b}-1\right) \, \dd t = b \int_0^1 \log \frac{1}{t} \, \dd t + \int_0^1 \log \left( 1- t^b\right) \, \dd t = b - \sum_{k=1}^\infty \frac 1k \frac{1}{kb+1},
\end{equation}
and
\begin{equation}\label{J1.def}
\cJ_1 = \int_1^\infty \log \left(1-t^{-b} \right) \, \dd t = - \sum_{k=1}^\infty \frac 1k \frac{1}{kb-1}.
\end{equation}
Therefore, by \eqref{cot.exp} with $w = 1/b$
\begin{equation}
\cJ_0 + \cJ_1 = b - \sum_{k=1}^\infty \frac 1k \left(
\frac{1}{kb+1} + \frac{1}{kb-1} 
\right)
= 
b + \sum_{k=1}^\infty \frac{\frac 2b}{-k^2+ \left(\frac1b\right)^2} = \pi \cot \frac{\pi}{b},
\end{equation}
which proves \eqref{An.b.asym}.
\end{proof}
Notice that a more accurate analysis of an error when Riemann sums \eqref{H01.def} are changed to the integrals \eqref{J0.def}, \eqref{J1.def} shows: if  $1<b<\infty$ and $b+1 \leq n$
\begin{equation}\label{An.b.asym.improved}
A(n;b) = \varphi(n) \exp \left(
n \pi \cot \frac \pi b
\right),
\end{equation}
where
\begin{equation}
	(en)^{-b} \leq \varphi(n) \leq e^3 b^{-2} n^2.
\end{equation}

\subsubsection{Convergence of partial fraction decomposition}

In the context of the necessary and sufficient conditions for construction of all undetermined moment problems, H.~Hamburger defined in \cite{Hamburger-1944-66} the class $\frU$ of entire functions $F(w)$ of finite order which fulfill the following conditions:
\begin{enumerate}[\upshape (i)]
	\item all roots $\{\xi_n\}_{n \in \N}$ of $F$ are simple and real;
	\item $F(x)$ is real for all real $x$;
	\item for all $w \neq \xi_n$, $n \in \N$,
	\begin{equation}
	\frac{1}{F(w)} = \sum_{k=1}^\infty \frac{1}{F'(\xi_k)} \frac{1}{w-\xi_k}; 
	\end{equation}
	\item for all $m=0,1,2, \dots$
	\begin{equation}
		\sum_{k=1}^\infty \xi_k^m \frac{1}{|F'(\xi_k)|} < \infty.
	\end{equation}
\end{enumerate}
He proved in \cite[Thm.~5a]{Hamburger-1944-66} the following:
Let $0 < a_1 < a_2 < \dots < a_n \to \infty$, and for $\nu >0$ and $\varrho \in (0,1/2)$
\begin{equation}
	n(r) = \# \{ k \, : \, a_k \leq r \} = \nu r^\varrho(1+o(1)), \quad r \to + \infty.
\end{equation}
Then
\begin{equation}
F(w) = \prod_{k=1}^\infty \left( 1 + \frac{w}{a_k} \right) \in \frU.
\end{equation}

There was a series of improvements and refinements of such a statement, see \cite[p.~167]{Akhiezer-1965}, \cite{DeBranges-1959-10}, \cite{Koosis-1980}. We choose as our basic reference the statement in \cite[Thm.~6.6]{Berg-1995-171}.

\begin{proposition}\label{prop:B-P}
Let $F$ be a non-constant entire function of minimal exponential type and assume that the zeros $\{\xi_n\}_{n \in \N}$ of $F$ are all real and simple. If
\begin{equation}\label{pfd.asm}
	\sum_{n=1}^\infty \frac{1}{|F'(\xi_n)|} < \infty,
\end{equation}
then
\begin{eqnarray}
	\frac 1{F(w)} = \sum_{n=1}^\infty \frac{1}{F'(\xi_n)} \frac{1}{w - \xi_n},
\end{eqnarray}
for all $w \neq \xi_n$, $n \in \N$.
\end{proposition}

\subsection{Proof of Propositions~\ref{prop:F.pfc.new} and \ref{prop:F.pfc.new.12}}
\label{ssec:App.A.3}

In the following, we have
\begin{equation}\label{F.def.app}
F(w) = \prod_{k=1}^\infty \left( 1 + \frac{w}{a_k} \right), \quad \text{where} \quad a_k = \left( \frac{k}{\nu} \right)^{\frac 1 \varrho}, \quad k \in \N,
\end{equation}
with $\nu>0$ and $\varrho \in (0,1)$ and we set $b=1/\varrho$.

\begin{proof}[Proof of Proposition~\ref{prop:F.pfc.new}]
By Proposition~\ref{prop:F.Titch}, the function $F$ in \eqref{F.def.app} is of exponential order $\varrho<1$ and therefore of minimal exponential type. 

Notice that	
\begin{equation}\label{F'.scal}
F'(-a_n) = \frac{1}{a_n} \prod_{\substack{k=1 \\k \neq n}}^\infty \left(1 - \frac{a_n}{a_k} \right) =
\frac{\nu^b}{n^b} \prod_{\substack{k=1 \\k \neq n}}^\infty \left(1 - \frac{n^b}{k^b} \right) 
%\\ & = \nu^b G_b' (-n^b) 
= (-1)^{n-1}\nu^b A(n;b), \quad n \in \N.
\end{equation}
Thus \eqref{F'.an.asym} follows from \eqref{An.b.asym.improved}; cf.~also \eqref{An.b.asym}.

If $\varrho <1/2$, i.e., $b>2$, then even a weaker estimate \eqref{An.b.est} justifies the
assumption~\eqref{pfd.asm}. Hence Proposition~\ref{prop:B-P} leads to \eqref{pfc.F.b>2.new}. Finally, \eqref{pf.decomp} follows by
	\begin{equation}\label{pfd.zw<12}
		\begin{aligned}
			\frac{1}{z-w} \left(
			\frac{1}{F(w)} - \frac{1}{F(z)}
			\right) &= \frac{1}{z-w} \sum_{n=1}^\infty \frac{1}{F'(-a_n)} \left( \frac{1}{w+a_n} - \frac{1}{z+a_n}   \right)
			\\
			& = \sum_{n=1}^\infty \frac{1}{F'(-a_n)} \frac{1}{(w+a_n)(z+a_n)} 
		\end{aligned}
	\end{equation}
for all $w \neq z$ and $w \neq -a_n$, $n \in \N$. 
\end{proof}

\begin{proof}[Proof of Proposition~\ref{prop:F.pfc.new.12}]
Recall that we have $\varrho =1/2$, i.e., $b=2$. To show \eqref{pfc.F.12}, we use \eqref{w/sin.exp} and \eqref{Euler.prod}, namely
\begin{equation}
\begin{aligned}
\frac{1}{F(w)} &= \frac{1}{\prod_{n=1}^\infty \left(1 + \frac{w^2}{\left(\frac{n}{\nu}\right)^2}\right)}
= 
\frac{1}{\prod_{n=1}^\infty \left(1 - \frac{(\ii \nu w)^2}{n^2}\right)}
= 
\frac{\pi \ii \nu w}{\sin(\ii \nu w)}
\\
&=
1 + 2\sum_{n=1}^\infty (-1)^n \frac{(\ii \nu w)^2}{(\ii \nu w)^2-n^2}
=
1 + 2\sum_{n=1}^\infty (-1)^n \frac{w^2}{w^2 + \left(\frac n\nu\right)^2}. 
\end{aligned}
\end{equation}
Thus \eqref{pf.decomp.12} is obtained by 
\begin{equation}
\begin{aligned}
\frac{1}{z-w} \left(
\frac{1}{F(w)} - \frac{1}{F(z)}
\right) &= \frac{2}{z-w} \sum_{n=1}^\infty (-1)^n \left( \frac{w}{w+a_n} - \frac{z}{z+a_n}   \right)
\\
& =  2 \sum_{n=1}^\infty (-1)^{n+1}\frac{a_n}{(z+a_n)(w+a_n)},
\end{aligned}
\end{equation}
for all $w$ which satisfy $w\neq z$ and $w \neq -a_n$, $n \in \N$. 
\end{proof}

\newpage

\section{Proof of Theorem~\ref{thm:L.res.lb}}
\label{app:res.L.lb}

\subsection{Setting and strategy}

Recall that $a, b \in \N$ with
\begin{equation}\label{a.b.asm.app}
	a -1 < b < 2a
\end{equation}
and the polynomial potential reads
\begin{equation}\label{V.L.app.def}
V(x) = x^{2a} + V_1(x), \quad x \in \R,
\end{equation}
with $\deg \Re V_1 \leq 2a-1$, $\deg \Im V_1 = b$ and $\Im V_1(x) = c_b x^b + \cdots$ with $c_b >0$,  cf.~\eqref{L.def}--\eqref{V.deg}. The lower order terms in $\Re V$ and $\Im V$ will be treated as a perturbation, using a constant $M>0$ such that
\begin{equation}\label{M.def}
\begin{aligned}
|V(z) - z^{2a}| &\leq M \left( |z|^{2a-1} +1 \right) , \quad z \in \C,
\\
|\Im V(x) - c_b x^b| &\leq M \left(|x|^{b-1} + 1\right), \quad x \in \R.
\end{aligned}
\end{equation}

We write $\la = \alpha + \ii \beta$ and assume that 
\begin{equation}\label{alp.bet.0}
\alpha \geq 1 + c_b^{2a}, \quad  \beta \geq 0,
\end{equation}
and, with fixed $\eps \in (0,c_b)$ and $\omega = \gamma \omega_0$ with $\gamma \in (0,1)$, see~\eqref{var.kap.def} and \eqref{rho.L.lb}, that
\begin{equation}\label{alp.bet.odd}
0 \leq \beta \leq (c_b-\eps) \alpha^\frac{b}{2a} \leq c_b \alpha^\frac{b}{2a}   \quad \text{if $b$ is odd}
\end{equation}
and 
\begin{equation}\label{alp.bet.even.0}
\alpha^{\frac{b}{2(b+1)} + \omega} \leq 	\beta \leq (c_b-\eps) \alpha^\frac{b}{2a}
\quad \text{if $b$ is even.}
\end{equation}
So for $b$ even, we have 
\begin{equation}\label{alp.bet.even.new}
\alpha^{\frac b{b+1} \left(
\frac12 + \gamma \tau	
	\right)} \leq \beta \leq  c_b \alpha^{\frac b{2a}}.
\end{equation}
The proof for the remaining case with $b$ odd and $\beta<0$ in the claim of Theorem~\ref{thm:L.res.lb} is analogous and we omit it.

In the following, $\la \in \C$ is assumed to always satisfy \eqref{alp.bet.0}--\eqref{alp.bet.even.new}. 
Notice that 
\begin{equation}\label{la.alp.ineq.0}
\alpha	\leq |\la| 
\leq 
\alpha +  c_b\alpha^{\frac{b}{2a}}
\leq 
\alpha + c_b \alpha^{1- \frac1{2a}} 
= \alpha \left(1 + c_b\alpha^{- \frac1{2a}}\right), 
\end{equation}
so by \eqref{alp.bet.0}
\begin{equation}\label{la.alp.ineq}
\alpha	\leq |\la| \leq 2 \alpha.
\end{equation}
Notice also that if $b$ is even, then $\beta \to + \infty$; if $b$ is odd, then $\beta$ can stay bounded.

We define the turning points $x_\alpha, y_\beta \geq 0$ of the leading terms of the real and imaginary parts of $V$ by
\begin{equation}\label{xa.yb.def}
	x_\alpha^{2a} = \alpha, \qquad c_b y_\beta^b = \beta.
\end{equation}
Notice that due to \eqref{alp.bet.odd}
\begin{equation}\label{yb.xa.est}
y_\beta = \left(\frac{\beta}{c_b}\right)^\frac1b \leq \left(1 - \frac{\eps}{c_b}\right)^\frac1b \alpha^\frac{1}{2a} 
=  \left(1 - \frac{\eps}{c_b}\right)^\frac1b x_\alpha < x_\alpha.
\end{equation}
We define
\begin{equation}\label{zeta.delta.def}
	\zeta \equiv \zeta(\eps) : = \left(1 - \frac{\eps}{c_b}\right)^\frac1b \in (0,1), \qquad \delta \equiv \delta(\eps) := \frac {1-\zeta}{32}  \in (0,\tfrac1{32})
\end{equation}
and
\begin{equation}\label{D.la.def}
	\Delta_\la := 
	\begin{cases}
		\delta x_\alpha & \text{if $b$ is odd},
		\\
		\delta y_\beta & \text{if $b$ is even}.
	\end{cases}	
\end{equation}
Notice that from \eqref{yb.xa.est},
\begin{equation}\label{yb.Del.est}
	y_\beta \leq \frac{\Delta_\la} \delta,
\end{equation}
and so in either case in \eqref{D.la.def}, using \eqref{alp.bet.even.new}, we obtain
\begin{equation}\label{Del.alp.ineq}
K^\frac{1}{b+1} \alpha^{\frac{1}{b+1}\left(\frac12 + \gamma \tau\right)}	\leq  \Delta_\la \leq \delta \alpha^\frac{1}{2a} \leq \alpha^\frac{1}{2a},
\end{equation}
where
\begin{equation}\label{K.def}
K:= \delta^{b+1} c_b^{-\frac{b+1}b}.
\end{equation}
Hence, if $\alpha \to +\infty$, then $\Delta_\la \to + \infty$. Moreover, we define
\begin{equation}\label{mu.la.def}
\mu_\la := \frac{\Delta_\la^{b+1}}{\alpha^\frac12}
\end{equation}
and note that by \eqref{Del.alp.ineq} we get
\begin{equation}\label{mu.la.ineq}
K  \alpha^{\gamma \tau}	\leq \mu_\la \leq \alpha^\tau.
\end{equation}
So in either case in \eqref{D.la.def}, if $\alpha \to +\infty$, then $\mu_\la \to +\infty$.

The inequalities \eqref{Del.alp.ineq} and \eqref{mu.la.def} show that there exists $\Lambda_1 \geq 1+c_b^{2a}$ such that for all $\alpha \geq \Lambda_1$, we have
\begin{equation}\label{Del.mu.alpha.ineq}
\Delta_\la \geq 4, \quad \mu_\la \geq 4, \quad \log \alpha \leq k_1 \log \Delta_\la \leq k_2 \log \mu_\la \leq k_3 \log \alpha,
\end{equation}
where the constants $k_1, k_2, k_3>0$ are independent of $\alpha$, $\Delta_\mu$, $\mu_\la$; $\Lambda_1$ can be take as

\begin{equation}\label{La.1.def}
\Lambda_1 := \max \left\{ 1+c_b^{2a}, \left(\frac{4^{b+1}}K\right)^\frac{1}{\gamma\tau} \right\}.
\end{equation}

Our goal is to show that there exists $\wt \eta>0$ such that for any $\la$ satisfying \eqref{alp.bet.0}--\eqref{alp.bet.even.new} and with $\alpha$ large enough we can find a function $u_\la \in C_0^\infty(\R) \setminus \{0\}$ which satisfies the inequality
\begin{equation}\label{pseudo.est}
	\|(\la-L)u_\la \|_{L^2} \leq \exp \left(-\wt \eta \mu_\la \right)  \|u_\la\|_{L^2}.
\end{equation}
Such functions $u_\la$ are often called pseudomodes. The resolvent estimates in the claims \eqref{L.res.lb.odd} and \eqref{L.res.lb.even} then follow from \eqref{pseudo.est}, \eqref{mu.la.def}, \eqref{D.la.def} and \eqref{xa.yb.def}. Indeed, for $\la \in \rho(L)$ satisfying \eqref{alp.bet.0}--\eqref{alp.bet.even.new} and with $\alpha$ large enough and $v_\la:=(\la-L)u_\la$, we obtain
\begin{equation}\label{L.res.lb.proof}
\|(\la-L)^{-1}\| \geq \frac{\|(\la-L)^{-1} v_\la\|_{L^2}}{\|v_\la\|_{L^2}} = \frac{\|u_\la\|_{L^2}}{\|(\la-L)u_\la\|_{L^2}} \geq 
\exp \left(\wt \eta \mu_\la \right).
\end{equation}

We find pseudomodes $u_\la$ in the form 
\begin{equation}\label{u.la.def}
u_\la = \chi_\la e^{\ii \varphi_\la} \sum_{j=0}^{N_\la} a_{j,\la}
\end{equation}
with the functions $\chi_\la$, $\varphi_\la$, $a_{j,\la}$ and the integer $N_\la \in \N$ constructed in the next subsections. 

\subsection{Notation}
To simplify notation, we suppress the subscripts $\la$ in the following. 
It is however essential to keep the $\la$-(in)dependence in mind. In particular, we have
\begin{equation}
\text{$\la$-dependent}: \quad u,\chi, a_j, \varphi, N, \Delta, \mu, r_k, I_j, 
\end{equation}
and
\begin{equation}
\text{$\la$-independent}: \quad \delta, \zeta, M, K, K_j, \Lambda_j, k_j.
\end{equation}

For $r>0$ and an analytic function $f$ in $\D_r$ we write
\begin{equation}\label{disc.norm.def}
\D_r := B_r(y_\beta)=\{z \in \C \, : \, |z-y_\beta|<r \}, \qquad \|f\|_{r} := \|f\|_{L^\infty(\D_r)}.
\end{equation}

\subsection{Smooth cut-offs $\chi$}
Let $h \in C_0^\infty(\R)$ be such that $0 \leq h \leq 1$ and
\begin{equation}
h(x) = 
\begin{cases}
1, \qquad x \in [-1,1],
\\
0, \qquad x \notin (-2,2);
\end{cases}
\end{equation}
the construction of such $h$ via mollification is standard, see e.g.~\cite[Lemma~V.1.9]{EE}. In our special case, we can also take 
\begin{equation}
h(x):= \frac1p \int_{-\infty}^x E(\tfrac32 +t) \, \dd t, \quad x \leq 0, \quad h(x):= h(-x), \quad x > 0,
\end{equation}
where 
\begin{equation}
p = \int_{-\infty}^\infty E(t) \, \dd t  > 0.22
\end{equation}
and $E$ is the bump function
\begin{equation}
E(t) = 
\begin{cases}
e^{- \frac{1}{1-4t^2}}, & |t| < \frac 12,
\\
0, & |t| \geq \frac12.		
\end{cases}
\end{equation}

The cut-off functions $\chi$ used below are defined via shift and scaling
\begin{equation}
	\chi(x) := h \left(\frac1{\Delta}(x-y_\beta)\right), \quad x \in \R,
\end{equation}
and so they satisfy $\chi \in C_0^\infty(\R)$, $0 \leq \chi \leq 1$, 
\begin{equation}\label{chi.supp}
\chi(x) = 
\begin{cases}
1, & x \in [y_\beta - \Delta, y_\beta + \Delta], 	
\\
0, & x \notin (y_\beta - 2  \Delta, y_\beta + 2 \Delta),
\end{cases}
\end{equation}
and 
\begin{equation}\label{chi.der}
\|\chi'\|_{L^\infty} = \frac{1}{\Delta} \|h'\|_{L^\infty} \leq \frac{5}{\Delta}, \quad \|\chi''\|_{L^\infty} = \frac{1}{\Delta^2} \|h''\|_{L^\infty} \leq \frac{10}{\Delta^2}.
\end{equation}

\subsection{Phase functions $\varphi$}

By \eqref{M.def} and triangle inequality, we get 
\begin{equation}
\begin{aligned}
|V(z)| 
&\leq |z|^{2a} + M \left( |z|^{2a-1} + 1\right)
\\
& 
\leq \left(|z-y_\beta| + y_\beta\right)^{2a} + M \left( \left(|z-y_\beta| + y_\beta\right)^{2a-1} + 1\right). 
\end{aligned}
\end{equation}
Thus, recalling \eqref{yb.xa.est}, \eqref{zeta.delta.def} and \eqref{D.la.def}, we obtain for all $z \in \D_{8 \Delta}$ that
\begin{equation}
\begin{aligned}
|V(z)|  & \leq  
(8 \Delta + y_\beta)^{2a} + M \left( \left( 8 \Delta + y_\beta \right)^{2a-1} + 1 \right) 
\\ 
& \leq 
\left(8 \delta x_\alpha  + \zeta x_\alpha \right)^{2a} \alpha + M \left( \left( 8 \delta x_\alpha  + \zeta x_\alpha  \right)^{2a-1}  + 1 \right)
\\
& 
= \left(
\frac14 \left(1 + 3 \zeta\right)
\right)^{2a}\alpha
+ M \left( \left(
\frac14 \left(1 + 3 \zeta\right)
\right)^{2a-1}\alpha^{1-\frac1{2a}} + 1 \right)
\\
& \leq 
\left[
\left(
\frac14 \left(1 + 3 \zeta\right)
\right)^{2a} 
+ M \left( \alpha^{-\frac1{2a}} + \alpha^{-1} \right)
\right] \alpha.
\end{aligned}
\end{equation}
Since $\zeta \in (0,1)$ and therefore $(1 + 3 \zeta)/4 < (1 + \zeta)/2<1$, there exists $\Lambda_2 \geq \Lambda_1$ such that for all $\alpha \geq \Lambda_2$ and all $z \in \D_{8 \Delta}$
\begin{equation}\label{V.ext}
\begin{aligned}
|V(z)| 
\leq
\left(
\frac12 \left(1 +  \zeta\right)
\right)^{2a}\alpha
\leq (1-K_1^2)\alpha < \alpha,
\end{aligned}
\end{equation}
where
\begin{equation}\label{K1.def}
	K_{1}:= 
	\left( 1 - 
	\left(
	\frac12 \left(1 +  \zeta\right)
	\right)^{2a}  
	\right)^\frac12 \in (0,1);
\end{equation}
a straightforward estimate shows that $\Lambda_2$ can be taken as
\begin{equation}\label{La.2.def}
\Lambda_2 = \max \left\{\Lambda_1,\left[\frac{2 M}{\left(
	\frac12 \left(1 +  \zeta\right)
	\right)^{2a}
 - \left(
 \frac14 \left(1 +  3\zeta\right)
 \right)^{2a}
 }\right]^{2a} \right\}.
\end{equation}
Hence, for all $\alpha \geq \Lambda_2$ and all  $z \in \D_{8 \Delta}$ 
\begin{equation}\label{alp.V.lb}
\alpha - |V(z)| \geq K_1^2 \alpha.
\end{equation}

The phase functions $\varphi$ are defined by
\begin{equation}
\varphi(z) := - \int_{y_\beta}^z (\la - V(w))^\frac12 \, \dd w, \qquad z \in \D_{8 \Delta},
\end{equation}
where the principal part of the complex square root is used and the integral is taken along the straight path connecting $y_\beta$ and $z$. Because of \eqref{V.ext} and $\alpha>0$, the functions $\varphi$ are analytic on $\D_{8 \Delta}$ and
\begin{equation}\label{phi.def}
\varphi'(z) = -(\la - V(z))^\frac12, \qquad z \in \D_{8 \Delta}.
\end{equation}
Moreover, by \eqref{V.ext}, \eqref{la.alp.ineq} and $\zeta<1$, we get for all $\alpha \geq \Lambda_2$ and all $z \in \ov \D_{8 \Delta}$
\begin{equation}\label{phi.der}
\begin{aligned}
|\varphi'(z)|
& \leq 	
\left( 1 + 
\left(
\frac12 \left(1 +  \zeta\right)
\right)^{2a}  
\right)^\frac12 |\la|^\frac12
\leq 4 \alpha^\frac12,
\\
|\varphi'(z)| &\geq 
\left( 1 - 
\left(
\frac12 \left(1 +  \zeta\right)
\right)^{2a}  
\right)^\frac12  |\la|^\frac12
\geq K_1 \alpha^\frac12,
\end{aligned}
\end{equation}
and therefore
\begin{equation}\label{1.phi'.norm}
\frac{1}{4\alpha^\frac12} 
\leq 
\left\| \frac{1}{\varphi'} \right\|_{4\Delta} 
\leq 
\frac{1}{K_1 \alpha^\frac12}.
\end{equation}

\subsection{Estimates of $\Re (\ii \varphi(x))$.}
For $t \in [y_\beta-2\Delta,y_\beta+ 2\Delta]$, we have from \eqref{phi.def} that
\begin{equation}\label{phi.sqrt}
	\begin{aligned}
		- \Re \left(\ii \varphi'(t)\right) &= \Im (\varphi'(t)) 
		\\
		&=  \frac{\Im V_1(t) -\beta }{2^\frac12\left(
			\left[
			(\alpha-\Re V(t))^2 + (\beta -\Im V_1(t))^2
			\right]^\frac12		   
			+ \alpha-\Re V(t) 		   
			\right)^\frac12}
		\\
		&= \frac{\Im V_1(t) -\beta }{D(t)},
	\end{aligned}
\end{equation}
where
\begin{equation}\label{Dt.def}
D(t) = 2^\frac12\left(
\left[
(\alpha-\Re V(t))^2 + (\beta -\Im V_1(t))^2
\right]^\frac12		   
+ \alpha-\Re V(t) 		   
\right)^\frac12;
\end{equation}
recall that $\alpha-\Re V(t)>0$ due to \eqref{alp.V.lb}. 

First we estimate $D(t)$ for $t \in [y_\beta-2 \Delta,y_\beta+2 \Delta]$. By \eqref{alp.V.lb}, we have for all $\alpha \geq \Lambda_2$ that
\begin{equation}\label{den.est.1}
D(t) \geq 
2(\alpha-\Re V(t))^\frac12 
\geq  2 K_1 \alpha^\frac12.
	\end{equation}
On the other hand, $\beta \leq c_b \alpha^{b/(2a)}$ and $b \leq 2a-1$, cf.~\eqref{alp.bet.odd}, \eqref{alp.bet.even.new},  \eqref{a.b.asm.app}, yield
\begin{equation}\label{den.est.2a}
D(t) \leq 
2
\big(
\alpha + |\Re V(t)| + \beta + |\Im V(t)|
\big)^\frac12
\leq
2  
\left(\alpha \big(1+ c_b \alpha^{-\frac{1}{2a}} \big) + 2 |V(t)|
\right)^\frac12.
\end{equation}
Thus using \eqref{V.ext} and \eqref{alp.bet.0}, we get for all $\alpha \geq \Lambda_2$ that
\begin{equation}\label{den.est.2}
D(t)		\leq
	2  
	\left( 2 \alpha + 2  \alpha \right)^\frac12 
	= 	4 \alpha^\frac12.
\end{equation}

Next we consider $\Re (\ii \varphi)$, namely,
\begin{equation}\label{I12.def.0}
- \Re (\ii \varphi(x))  = \int_{y_\beta}^x \Im(\varphi'(t)) \, \dd t 
		 = I_1(x) +  I_2(x), \quad x \in  [y_\beta-2\Delta,y_\beta+ 2\Delta],
\end{equation}
where
\begin{equation}\label{I12.def}
I_1(x):=c_b \int_{y_\beta}^x \frac{t^b-y_\beta^b }{D(t)}  \, \dd t, \qquad 
I_2(x):= \int_{y_\beta}^x \frac{\Im V_1(t) - c_b t^b }{D(t)} \, \dd t.
\end{equation}

First notice that 
\begin{equation}\label{I1.pos}
	I_1(x) \geq 0, \quad x \in [y_\beta - 2 \Delta, y_\beta + 2 \Delta].
\end{equation}
Indeed, $D(t)>0$ for all $t \in [y_\beta - 2 \Delta, y_\beta + 2 \Delta]$ and for $x \geq y_\beta$, we get
\begin{equation}
I_1(x) = c_b \int_{y_\beta}^x \frac{t^b-y_\beta^b }{D(t)}  \, \dd t \geq 0.
\end{equation}
On the other hand, for $x \in [y_\beta-2 \Delta,y_\beta]$, we have 
\begin{equation}\label{I1.neg}
I_1(x) = c_b \int_{y_\beta}^x \frac{t^b-y_\beta^b }{D(t)}  \, \dd t = c_b \int_x^{y_\beta} \frac{y_\beta^b-t^b }{D(t)}  \, \dd t.
\end{equation}
If $b$ is odd, then $y_\beta^b-t^b \geq 0$ for all $t \leq y_\beta$. If $b$ is even, we selected $\Delta \leq y_\beta$, see \eqref{D.la.def}, thus for $t \in [y_\beta-2\Delta,y_\beta]$, we have $|t| \leq y_\beta$ and hence
$
y_\beta^b-t^b \geq y_\beta^b-y_\beta^b 	=0.
$

Second we estimate $I_1(x)$ from below for all $x \in \supp \chi'$; it follows from \eqref{chi.supp} that
\begin{equation}\label{x.supp.chi'}
\supp \chi' \subset [y_\beta-2 \Delta,y_\beta- \Delta] \cup [y_\beta + \Delta,y_\beta + 2 \Delta].
\end{equation}
For all $x \in [y_\beta + \Delta,y_\beta + 2\Delta]$ and all $\alpha \geq \Lambda_2$, we obtain (recall \eqref{den.est.2}, \eqref{mu.la.def})
\begin{equation}\label{I1.lb.0}
	\begin{aligned}
		I_1(x) & \geq \frac{c_b}{4} \frac{1}{\alpha^\frac12} \int_{y_\beta}^x (t^b-y_\beta^b) \, \dd t  
		\geq
		\frac{c_b}{4} \frac{1}{\alpha^\frac12} \int_{y_\beta+\frac \Delta 2 }^{y_\beta+ \Delta} (t^b-y_\beta^b) \, \dd t   	
		\\
		& 
		\geq
		\frac{c_b}{4} \frac{1}{\alpha^\frac12}
		\left( \left(y_\beta+\frac \Delta 2 \right)^b - y_\beta^b \right) \frac{\Delta}{2}
		\geq 
		\frac{c_b}8 \frac{1}{2^b} \frac{\Delta^{b+1}}{\alpha^\frac12} 
		= \frac{c_b}8 \frac{1}{2^b} \mu.
	\end{aligned}
\end{equation}
Similarly, for all $x \in [y_\beta - 2\Delta,y_\beta - \Delta]$ and all $\alpha \geq \Lambda_2$, we have (cf.~\eqref{I1.neg}, \eqref{den.est.2})
\begin{equation}\label{I1.lb.1}
	I_1(x)  \geq \frac{c_b}{4} \frac{1}{\alpha^\frac12} \int_{y_\beta-\Delta}^{y_\beta-\frac\Delta2} (y_\beta^b-t^b) \, \dd t  
		\geq
		\frac{c_b}{8} \frac{\Delta }{\alpha^\frac12}
		\left( y_\beta^b - \left(y_\beta-\frac \Delta 2 \right)^b\right).
\end{equation}
In the next steps we analyze separately the cases $\Delta > 3 y_\beta$ and $\delta y_\beta \leq \Delta \leq 3 y_\beta$; the lower bound on $\Delta$ is due to \eqref{yb.Del.est}. The first case can occur only for $b$ is odd since $\Delta \leq y_\beta$ for $b$ even, see \eqref{D.la.def}. Thus for $b$ odd and $\Delta > 3 y_\beta$,
\begin{equation}\label{I1.neg.1}
 y_\beta^b - \left(y_\beta-\frac \Delta 2\right)^b = y_\beta^b + \left(\frac \Delta 2-y_\beta \right)^b  \geq   \left(\frac\Delta 6 + \frac \Delta3  - y_\beta \right)^b  
 \geq \frac1 {6^b} \Delta^b.
\end{equation}
In the second case $\delta y_\beta \leq \Delta \leq 3 y_\beta$, we get
\begin{equation}
 -\frac12 y_\beta \leq 	y_\beta- \frac{\Delta}{2} \leq \left(1-\frac{\delta}{2}\right) y_\beta.
\end{equation}
Since $\delta<1$, we have 
\begin{equation}
\left|y_\beta- \frac{\Delta}{2}\right| \leq \left(1-\frac{\delta}{2}\right) y_\beta
\end{equation}
and therefore
\begin{equation}\label{I1.neg.2}
y_\beta^b - \left(y_\beta-\frac \Delta 2\right)^b 
\geq
y_\beta^b - \left(1-\frac{\delta}{2}\right)^b y_\beta^b
\geq 
\frac{1}{3^b}
\left(
1- \left(1-\frac{\delta}{2}\right)^b
\right) \Delta^b.
\end{equation}

Summarizing \eqref{I1.lb.0} and \eqref{I1.lb.1} together with \eqref{I1.neg.1} and \eqref{I1.neg.2}, we obtain that for all $ x \in \supp \chi'$, cf.~\eqref{x.supp.chi'}, and all $\alpha \geq \Lambda_2$

\begin{equation}\label{I1.lb}
I_1(x) \geq K_2 \mu,
\end{equation}
where 
\begin{equation}\label{K2.def}
K_2 := \frac{c_b}{8} \min \left\{\frac{1}{6^{b}}, \frac{1}{3^b } \left(
1- \left(1-\frac{\delta}{2}\right)^b
\right)  \right\}.
\end{equation}

Next, to estimate $I_1(x)$ from above for $x \in [y_\beta,y_\beta + 2 \Delta]$, we use \eqref{den.est.1}, the mean value theorem and $y_\beta \leq \delta^{-1} \Delta$, cf.~\eqref{yb.Del.est}, and obtain for all $\alpha \geq \Lambda_2$
\begin{equation}\label{I1.ub}
\begin{aligned}
I_1(x) & \leq \frac{c_b}{2 K_1} \frac{1}{\alpha^\frac12} \int_{y_\beta}^x  (t^b-y_\beta^b) \, \dd t 
\leq 
\frac{b c_b}{2 K_1} \frac{x^{b-1}}{\alpha^\frac12}  \int_{y_\beta}^x  (t-y_\beta) \, \dd t 
\\
&
\leq 
\frac{b c_b}{4 K_{1}} 
\frac{(y_\beta+2\Delta)^{b-1}}{\alpha^\frac12} (x-y_\beta)^2
\\
&\leq
\frac{b c_b (\delta^{-1} +2)^{b-1}}{4 K_{1}} 
\frac{\Delta^{b-1}}{\alpha^\frac12} (x-y_\beta)^2
= \frac{K_3}{2} \frac{\mu}{\Delta^2} (x-y_\beta)^2,
\end{aligned}
\end{equation}
where 
\begin{equation}\label{K3.def}
K_3 := 	\frac{b c_b (\delta^{-1} +2)^{b-1}}{2 K_{1}}. 
\end{equation}

Finally, we estimate $I_2$ from above and also show that it is a perturbation of $I_1$ on $\supp \chi'$, cf.~\eqref{x.supp.chi'}. To this end, recalling \eqref{M.def}, \eqref{den.est.1} and $y_\beta  \leq \delta^{-1} \Delta$, cf.~\eqref{yb.Del.est}, we obtain for all $x \in [y_\beta-2 \Delta,y_\beta+2 \Delta]$ and all $\alpha \geq \Lambda_2$ that 
\begin{equation}\label{I2.ub.0}
\begin{aligned}
|I_2(x)| 
&\leq 
\frac{M}{2 K_1} \frac{1}{\alpha^\frac12} \int_{y_\beta-2 \Delta}^{y_\beta+2\Delta}  (|t|^{b-1} +1) \, \dd t 
\leq 
\frac{M}{2 K_1} \frac{1}{\alpha^\frac12} 4 \Delta ((y_\beta+2\Delta)^{b-1}+1) 
\\
&\leq 
\frac{2 M}{ K_1}
\frac{\Delta}{\alpha^\frac12} \left( (\delta^{-1}+2)^{b-1}\Delta^{b-1} +1 \right)
=
\frac{2 M}{ K_1}
\mu \left( \frac{(\delta^{-1}+2)^{b-1}}{\Delta} +\frac{1}{\Delta^{b}} \right).
\end{aligned}
\end{equation}
Hence, there exists $\Lambda_3 \geq \Lambda_2$ such that for all $\alpha \geq \Lambda_3$ and all $x \in [y_\beta-2 \Delta,y_\beta+2 \Delta]$
\begin{equation}\label{I2.ub.1}
|I_2(x)|
\leq 
\frac{4 M}{ K_1} \frac{\mu}{\Delta} 
\leq 
\frac{K_2}{2} \mu;
\end{equation}
using \eqref{Del.alp.ineq}, $\Lambda_3$ can be taken as
	\begin{equation}
\Lambda_3 = \max \left\{ \Lambda_2, \frac{(\delta^{-1}+2)^{2(b-1)(b+1)}}{K^2}, \left(\frac{8M}{K_1 K_2}\right)^{2(b+1)} \frac{1}{K^2} \right\} .
	\end{equation}

In summary, using \eqref{I1.lb} and \eqref{I2.ub.1}, we have for all $\alpha \geq \Lambda_3$ and all $x \in \supp \chi'$, cf.~\eqref{x.supp.chi'}, that
\begin{equation}\label{I12.est.0}
- \Re (\ii \varphi(x)) = I_1(x) + I_2(x) \geq K_2\mu - \frac{K_2}2 \mu = \frac {K_2}2  \mu.
\end{equation}

\subsection{Functions $a_{j}$}

The application of the operator $\la-L$ on $u$ from \eqref{u.la.def} yields
\begin{equation}\label{Lu.1}
(\la-L)u = \chi'' e^{\ii \varphi} \sum_{j=0}^{N} a_{j} + 2 \chi' e^{\ii \varphi} \sum_{j=0}^{N} \left(\ii \varphi'  a_{j} + a_{j}'   \right) 
+ \chi (\la-L) e^{\ii \varphi} \sum_{j=0}^{N} a_{j}
\end{equation}
and, using \eqref{phi.def},
\begin{equation}
\begin{aligned}
(\la-L) \left(e^{\ii \varphi} \sum_{j=0}^{N} a_{j}\right) & =
e^{\ii \varphi} \left(
\Dtp + 2 \ii \varphi' \Do + \ii \varphi'' - \varphi'^2 - V +\la 
\right)\sum_{j=0}^{N} a_{j}
\\
& = 
e^{\ii \varphi}  
\sum_{j=0}^{N} \left(
 a_{j}''
+ 2 \ii \varphi' a_{j}'
+ \ii \varphi'' a_{j}
\right).
\end{aligned}
\end{equation}

The functions $a_{j}$ are selected to satisfy the system of equations
\begin{equation}
\begin{aligned}
2\varphi' a_{0}' +  \varphi'' a_{0} &= 0,
\\
2\ii \varphi' a_{j}' +  \ii \varphi'' a_{j} &= - a_{j-1}'', \quad j \in \N, 
\end{aligned}
\end{equation}
which leads to
\begin{equation}\label{Lu.N}
(\la-L) e^{\ii \varphi} \sum_{j=0}^{N} a_{j}  =  e^{\ii \varphi} a_{N}''.
\end{equation}
For $z \in \D_{8 \Delta}$, we set
\begin{equation}\label{a0.int}
a_{0}(z) = \frac{\varphi'(y_\beta)^\frac12}{\varphi'(z)^\frac12} =
\frac{(\la-V(y_\beta))^\frac14}{(\la-V(z))^\frac14},	
\end{equation}
and by induction
\begin{equation}\label{aj.int.1}
a_{j}(z) = \frac{1}{\varphi'(z)^\frac12}\int_{y_\beta}^z \frac{ \ii a_{j-1}''(w)}{2 \varphi'(w)^\frac12} \, \dd w 
= 
\frac{1}{\varphi'(z)^\frac12}\int_{0}^{z-y_\beta} \frac{ \ii a_{j-1}''(y_\beta+w)}{2 \varphi'(y_\beta+w)^\frac12} \, \dd w
, \quad j \in \N;
\end{equation}
the integrals are taken along the straight paths connecting the end-points.
The functions $\varphi'^{-1/2}$ and $a_j$, $j = 0, 1,\dots $, are analytic on  $\D_{8 \Delta}$.

\subsection{Estimates of $a_{j}$}
First, it follows from \eqref{phi.der} and \eqref{1.phi'.norm} that for all $z \in \ov \D_{8 \Delta}$ and all $\alpha \geq \Lambda_3$ 
\begin{equation}\label{a0.est}
\frac{K_1^\frac12}{2}  \leq |a_{0}(z)| \leq \frac2{K_1^\frac12}.
\end{equation}
Next, by Cauchy estimates for an analytic and bounded function $g$ on $B_r(y)$ with $y \in \C$ and $r>0$, we have
\begin{equation}\label{Cauchy.est}
|g^{(l)}(y)|\leq \frac{l!}{r^l} \|g\|_{L^\infty(B_r(y))}, \quad l = 1,2. 
\end{equation}

Fix $j \in \N$ and define
\begin{equation}
r_k \equiv r_k^{(j)} = \left(4 - \frac kj \right) \Delta, \quad k =0,\dots, j;
\end{equation}
notice that
\begin{equation}\label{rk.k-1}
r_{k-1} - r_k = \frac{\Delta}{j}, \quad k = 1, \dots, j,
\end{equation}
and
\begin{equation}\label{rj.ineq}
	4 \Delta  = r_0 > r_1 > \dots > r_j =  3 \Delta.
\end{equation}

Let $z \in \D_{r_k}$ for $k \in \{1,\dots,j\}$. From \eqref{aj.int.1} we get
\begin{equation}\label{ak.est}
|a_{k}(z)| \leq \frac{1}{2} \left\| \frac{1}{\varphi'} \right\|_{4\Delta} \int_{0}^{|z-y_\beta|}  |a_{k-1}''(y_\beta+w)| \, |\dd w|.
\end{equation}
By \eqref{Cauchy.est} for $a_{k-1}''$ with $y=y_\beta+w$ and $r=r_{k-1}-|w|$ in the disc
\begin{equation}
B_{r_{k-1}-|w|}(y_\beta+w) \subset \D_{r_{k-1}} =B_{r_{k-1}}(y_\beta), 
\end{equation}
we have
\begin{equation}\label{ak.Cauchy}
|a_{k-1}''(y_\beta+w)| \leq \frac{2 \|a_{k-1}\|_{r_{k-1}}}{( r_{k-1}-|w| )^2}.
\end{equation}
Thus returning to \eqref{ak.est}, we get for all $z \in \D_{r_k}$ that
\begin{equation}\label{ak.est.2}
\begin{aligned}
|a_{k}(z)| & \leq \|a_{k-1}\|_{r_{k-1}} \left\| \frac{1}{\varphi'} \right\|_{4\Delta}  \int_{0}^{|z-y_\beta|}  \frac{|\dd w|}{( r_{k-1}-|w| )^2}
\\
&= \|a_{k-1}\|_{r_{k-1}} \left\| \frac{1}{\varphi'} \right\|_{4\Delta}  
\left(
\frac{1}{r_{k-1}}  
- \frac{1}{r_{k-1}-|z-y_\beta|}
\right)
\\
& = \|a_{k-1}\|_{r_{k-1}} \left\| \frac{1}{\varphi'} \right\|_{4\Delta}  \frac{|z-y_\beta|}{r_{k-1} (r_{k-1}-|z-y_\beta|)}.
\end{aligned}
\end{equation}
Moreover, using $r_k< r_{k-1}$ and $r_{k-1} - r_k = \Delta/j$, see \eqref{rj.ineq} and \eqref{rk.k-1}, we arrive at

\begin{equation}\label{ak.est.3}
	\begin{aligned}
		\|a_{k}\|_{r_k}  & \leq  \|a_{k-1}\|_{r_{k-1}} \left\| \frac{1}{\varphi'} \right\|_{4\Delta}  \frac{r_{k-1}}{r_{k-1} (r_{k-1}-r_k)}
		%\\
		& 
		= \|a_{k-1}\|_{r_{k-1}} \left\| \frac{1}{\varphi'} \right\|_{4\Delta}   \frac{j}{ \Delta},
	\end{aligned}
\end{equation}
for $k=1,\dots,j$. 

In particular for $k=j$ we get from \eqref{ak.est.2} that for all $z \in \D_{r_j}$ (recall $r_{j-1}>r_j = 3 \Delta$)

\begin{equation}
\begin{aligned}
|a_{j}(z)| & \leq 
\|a_{j-1}\|_{r_{j-1} } \left\| \frac{1}{\varphi'} \right\|_{4\Delta}
\frac{1}{r_{j-1}-|z-y_\beta|} 
\frac{|z-y_\beta|}{r_{j-1} }
\\
& \leq 
\|a_{j-1}\|_{r_{j-1} } \left\| \frac{1}{\varphi'} \right\|_{4\Delta}
\frac{1}{r_{j-1} - r_{j} } 
\frac{|z-y_\beta|}{r_{j} }
\\
& = 
\frac{|z-y_\beta|}{3 \Delta}  \left\| \frac{1}{\varphi'} \right\|_{4\Delta} \frac{j}{\Delta}
\|a_{j-1}\|_{r_{j-1} }.
\end{aligned}
\end{equation}
Moreover, the iterated application of \eqref{ak.est.3} yields that for all $z\in \D_{3 \Delta}$
\begin{equation}\label{aj.est.z}
|a_{j}(z)| \leq \frac{|z-y_\beta|}{3 \Delta}
\left(
\left\| \frac{1}{\varphi'} \right\|_{4\Delta} \frac{j}{\Delta}
\right)^j
\|a_{0}\|_{4 \Delta}, \qquad j = 1, 2, \dots,
\end{equation}
and hence also
\begin{equation}\label{aj.est.norm}
	\|a_{j}\|_{3 \Delta} \leq 
	\left(
	\left\| \frac{1}{\varphi'} \right\|_{4\Delta} 
	\frac{j}{\Delta}
	\right)^j
	\|a_{0}\|_{4 \Delta}, \qquad j = 1, 2, \dots.
\end{equation}

Finally, by the Cauchy estimates in the disc $B_\Delta(z)$, we obtain for all $z \in \D_{2 \Delta}$  
\begin{equation}\label{aj.Cauchy}
|a_{j}^{(l)}(z)| \leq \frac{l!}{\Delta^l} \|a_{j}\|_{3 \Delta}, \quad l = 1,2, \quad j = 0, 1, \dots.
\end{equation}

\subsection{The choice of $N$ and related estimates}

Let 
\begin{equation}
N := \left \lfloor  \left\| \frac{1}{\varphi'} \right\|_{4\Delta}^{-1}  \Delta  \frac 1e    \right \rfloor,
\end{equation}
thus
\begin{equation}\label{N.e.ineq}
\left\|\frac{1}{\varphi'} \right\|_{4\Delta} \frac{N}{\Delta}	\leq \frac{1}{e}
\end{equation}
and by \eqref{1.phi'.norm}, for all $\alpha \geq \Lambda_3$
\begin{equation}\label{N.la.est}
\frac{K_1}e \alpha^\frac12 \Delta  -1 \leq N \leq  4 \alpha^\frac12 \Delta \frac 1 e.
\end{equation}
In particular, $N \to + \infty$ if $\alpha \to +\infty$ due to $\Delta \to +\infty$, see \eqref{Del.alp.ineq}.

Taking $j=N$ in the estimates above, namely \eqref{aj.Cauchy} with $l=2$, \eqref{aj.est.norm}, we obtain by \eqref{N.e.ineq} and the estimate on $\|a_0\|_{4 \Delta}$ from \eqref{a0.est} that
\begin{equation}\label{aN.est.0}
\|a_{N}''\|_{2 \Delta}  \leq \frac{2}{\Delta^2} \|a_{N}\|_{3 \Delta}
\leq  \frac{2}{\Delta^2} \left(
\left\| \frac{1}{\varphi'} \right\|_{4\Delta} \frac{N}{\Delta}
\right)^{N}
 \|a_{0}\|_{4 \Delta} 
\leq 
\frac{4}{K_1^\frac12} \frac{1}{\Delta^2}  
\exp(-N).
\end{equation}
Since by \eqref{D.la.def}, \eqref{mu.la.def}, \eqref{yb.xa.est} and $b\leq 2a-1$
\begin{equation}
	\frac{\mu}{\alpha^\frac12 \Delta}
	=
	\frac{\Delta^{b}}{\alpha} \leq \frac{x_\alpha^b}{\alpha}  
	= \alpha^{\frac{b}{2a}-1} \leq  \alpha^{-\frac{1}{2a}},
\end{equation}
we have
\begin{equation}
\frac{N}{\mu} 
\geq 
\frac{K_1}{e} \frac{\alpha^\frac12 \Delta}{\mu} - \frac1\mu
\geq 
\frac{K_1}{e} \alpha^\frac{1}{2a} - \frac{1}{\mu} \geq K_2;
\end{equation}
the last step is valid if $\alpha \geq \La_4$ where
\begin{equation}\label{La4.def}
\La_4: =\La_3+ \left( \frac3 {K_1}(K_2+1) \right)^{2a}.
\end{equation}
Hence for all $\alpha \geq \Lambda_4$, we get
\begin{equation}\label{aN.est.1}
\|a_{N}''\|_{2 \Delta} 
\leq  
\frac{4}{K_1^\frac12} \frac{1}{\Delta^2}  
\exp \left(
-\mu \frac N \mu
\right)
\leq   
\exp \left(
-K_2 \mu
\right).
\end{equation}

By \eqref{aj.est.norm} and \eqref{N.e.ineq}, we get for all $\alpha \geq \Lambda_4$ that

\begin{equation}\label{aj.est.fin}
\begin{aligned}
\sum_{j=0}^{N} \| a_{j} \|_{2 \Delta} 
&
\leq 
\|a_{0}\|_{4 \Delta}
\left(1+ 
\sum_{j=1}^{N}
\left(
\left\| \frac{1}{\varphi'} \right\|_{4\Delta} 
\frac{j}{\Delta}
\right)^j
\right)
\\
& 
\leq 
\|a_{0}\|_{4 \Delta}
\left(1+ 
\sum_{j=1}^{N}
\left(
\left\| \frac{1}{\varphi'} \right\|_{4\Delta} 
\frac{N}{\Delta}
\right)^j
\right)
\\
& \leq 
\|a_{0}\|_{4 \Delta} \sum_{j=0}^{N} e^{-j} 
\leq  
\frac2{K_1^\frac12}  \frac{e}{e-1} \leq \frac4{K_1^\frac12}.
\end{aligned}
\end{equation}
Similarly by \eqref{aj.Cauchy} with $l=1$ and \eqref{aj.est.norm},  we have for all $\alpha \geq \Lambda_4$ that (see \eqref{Del.mu.alpha.ineq} for the last step)
\begin{equation}\label{aj'.est.fin}
\sum_{j=0}^{N} \| a_{j}' \|_{2 \Delta} 
\leq 
\frac{\|a_{0}\|_{4 \Delta}}{\Delta}
 \sum_{j=0}^{N} e^{-j} 
 \leq 
\frac2{K_1^\frac12} \frac{1}{\Delta} \frac{e}{e-1} \leq \frac4{K_1^\frac12} \frac1\Delta \leq \frac1{K_1^\frac12}.
\end{equation}
Finally, for $s \in (0,1)$, by \eqref{aj.est.z}, we get for all $\alpha \geq \Lambda_4$ that
\begin{equation}\label{aj.s.est}
\sum_{j=1}^{N} \| a_{j} \|_{s \Delta} 
\leq 
\frac{s}{3} \| a_{0} \|_{4 \delta \Delta} \sum_{j=1}^{N} e^{-j} 
\leq 
\frac s3  \frac2{K_1^\frac12} \frac{1}{e-1}
\leq 
\frac13 \frac2{K_1^\frac12} s.
\end{equation}
Hence for 
\begin{equation}\label{s0.def}
s_0 := \frac{K_1}2  \in (0, \tfrac12)
\end{equation}
we get by \eqref{a0.est} and \eqref{aj.s.est} that for all $x \in [y_\beta-s_0 \Delta,y_\beta+s_0 \Delta]$ and all $\alpha \geq \Lambda_4$
\begin{equation}\label{a0.aj.est}
|a_{0}(x)| - \sum_{j=1}^{N} | a_{j}(x)| \geq
\frac{K_1^\frac12}2 - \frac13 \frac2{K_1^\frac12} \frac{K_1}2
= 
\frac{K_1^\frac12}6 .
\end{equation}

\subsection{The upper estimate of $\|(\la-L)u\|$}
We estimate the terms in \eqref{Lu.1} and \eqref{Lu.N}. Recall \eqref{chi.supp}, \eqref{chi.der} and
\begin{equation}
\supp \chi'' \subset \supp \chi' \subset [y_\beta-2 \Delta,y_\beta- \Delta] \cup [y_\beta + \Delta,y_\beta + 2 \Delta].
\end{equation}
Moreover, on $\supp \chi'$ we have the estimates of $\Re (\ii \varphi)$ in \eqref{I12.est.0} and on the sum of $\|a_{j}\|_{2 \Delta}$ in \eqref{aj.est.fin}. Hence, for all $\alpha \geq \Lambda_4$, 
\begin{equation}\label{chi''.aj.est}
\begin{aligned}
\left\|\chi'' e^{\ii \varphi} \sum_{j=0}^{N} a_{j} \right\|^2_{L^2}
&\leq 
\int \limits_{\supp \chi''} \|\chi''\|_{L^\infty}^2  \exp \left(2 \Re (\ii \varphi(x) ) \right) \left(\sum_{j=0}^{N} \|a_{j}\|_{2 \Delta}\right)^2 \, \dd x  
\\
& 
\leq 
\frac{100}{\Delta^4}  
\int \limits_{\supp \chi''}  \exp \left(- K_2 \mu \right) \frac{16}{K_1} \, \dd x
\\
& \leq 
\frac{1600}{K_1 \Delta^4}  \exp \left(- K_2 \mu \right)  2 \Delta 
\leq 
\frac{50}{K_1}
\exp\left(- K_2 \mu \right),
\end{aligned}
\end{equation}
for the last step recall that $\Delta \geq 4$ by \eqref{Del.mu.alpha.ineq}.

Similarly, using \eqref{aj.est.fin}, \eqref{aj'.est.fin}, \eqref{phi.der} and \eqref{chi.der}, we obtain for all $\alpha \geq \Lambda_4$
\begin{equation}\label{chi'.aj.est}
\begin{aligned}
 \left\|2 \chi' e^{\ii \varphi} \sum_{j=0}^{N} \left(\ii \varphi'  a_{j} + a_{j}'   \right) \right\|^2_{L^2}
& \leq 4 \|\chi'\|_{L^\infty}^2 2 \left( 16 \alpha  \frac{16}{K_1} + \frac{1}{K_1} \right) \int \limits_{\supp \chi'}  \exp \left(- K_2 \mu \right) \, \dd x
\\
& 
\leq 
\frac{16^5}{K_1}  \alpha  \exp \left(- K_2 \mu \right).
\end{aligned}
\end{equation}

Finally, we estimate the resulting term from \eqref{Lu.N}, i.e., $\chi  e^{\ii \varphi} a_{N}''$.
To this end, we first employ \eqref{aN.est.1} and  \eqref{I12.def} and obtain for all $\alpha \geq \Lambda_4$
\begin{equation}
\begin{aligned}
\left\| \chi  e^{\ii \varphi} a_{N}'' \right\|^2_{L^2}  & \leq \int_{y_\beta-2\Delta}^{y_\beta+2\Delta} \exp(2 \Re(\ii \varphi(x))) 
\exp \left(
-2 K_2 \mu 
\right) \, \dd x
\\
&\leq
\int_{y_\beta-2\Delta}^{y_\beta+2\Delta} \exp(- 2 I_1(x)) \exp( 2 |I_2(x)|)
\exp \left(
- 2 K_2 \mu
\right) \, \dd x.
\end{aligned}
\end{equation}
Next by \eqref{I1.pos} and \eqref{I2.ub.1}, we have for all $\alpha \geq \Lambda_4$
\begin{equation}\label{chi.aN.est}
\left\| \chi  e^{\ii \varphi} a_{N}'' \right\|^2_{L^2}   
\leq 
\int_{y_\beta-2\Delta}^{y_\beta+2\Delta} 
\exp \left(
-2 K_2 \mu
+ 
K_2 \mu
\right) \, \dd x
= 4 \Delta \exp \left(
- K_2 \mu 
\right).
\end{equation}
Putting together \eqref{chi''.aj.est}, \eqref{chi'.aj.est}, \eqref{chi.aN.est} and \eqref{Lu.1}, \eqref{Lu.N}, we arrive at
\begin{equation}\label{Lu.ub}
\|(\la-L)u\|^2_{L^2}    
 \leq 
4\left(  \frac{50}{K_1}
+  \frac{16^5}{K_1} \alpha 
+ 4 \Delta 
\right)
\exp \left(
- K_2 \mu 
\right)
\leq \frac{16^6}{K_1} \alpha \exp \left(
- K_2 \mu 
\right).
\end{equation}
for all $\alpha \geq \Lambda_4$; in the second step we used \eqref{alp.bet.0}, \eqref{Del.alp.ineq} and \eqref{K1.def}. 
\subsection{The lower estimate of $\|u\|$}
 
Let $s_0 \in (0,1/2)$ be as in \eqref{s0.def}. By \eqref{chi.supp} and \eqref{I12.def} we get 

\begin{equation}
\begin{aligned}
\|u\|^2_{L^2} & \geq 
\int_{y_\beta}^{y_\beta + s_0 \Delta} |\chi(x)|^2 \exp \left(2 \Re(\ii \varphi(x)) \right)  \left|
\sum_{j=0}^{N} a_{j}(x) 
\right|^2
\, \dd x
\\ & \geq 
\int_{y_\beta}^{y_\beta + s_0 \Delta} \exp \left(- 2 I_1(x) - 2|I_2(x)|\right) \left(
|a_{0}(x)| - \sum_{j=1}^{N} |a_{j}(x)| 
\right)^2
\, \dd x.
\end{aligned}
\end{equation}
Next, employing \eqref{I1.ub}, \eqref{I2.ub.1}, \eqref{a0.aj.est} and changing the integration variable in the second step, we obtain that for all $\alpha \geq \Lambda_4$
\begin{equation}
\begin{aligned}
\|u\|^2_{L^2} & \geq 
\exp
\left(-\frac{8 M}{K_1}
\frac{\mu}{\Delta} 
\right) 
\int_{y_\beta}^{y_\beta + s_0 \Delta} 
\exp \left( - K_3  \frac{\mu}{\Delta^2} (x-y_\beta)^2\right)  \frac{K_1}{36} \, \dd x
\\
& =
\frac{K_1}{36}
\exp
\left(-\frac{8 M}{K_1}
\frac{\mu}{\Delta} 
\right) 
\frac{\Delta}{\mu^\frac12}
\int_{0}^{s_0 \mu^\frac12} \exp \left( - K_3 y^2\right) \, \dd y.
\end{aligned}
\end{equation}
Thus if $\mu \geq 1/ s_0^2$, we have for all $\alpha \geq \Lambda_4$ (recall that $\Delta\geq 4$ by \eqref{Del.mu.alpha.ineq})
\begin{equation}\label{u.est}
\|u\|^2_{L^2} \geq 
\frac{K_1}{9} 
\exp
\left(-\frac{8 M}{K_1}
\frac{\mu}{\Delta} 
\right) 
\frac{1}{\mu^\frac12}
\exp(-K_3). 
\end{equation}

\subsection{The lower estimates of the resolvent norm}
Finally, combing \eqref{Lu.ub} and \eqref{u.est}, we have for all $\alpha \geq \Lambda_4$ that (if $\mu \geq 1/s_0^2$)
\begin{equation}
\frac{\|(\la-L)u\|_{L^2}^2}{\|u\|_{L^2}^2} 
\leq 
\exp \left(-K_2 \mu + \log \alpha+  \frac{8 M}{K_1}
\frac{\mu}{\Delta} + \frac12 \log \mu  + \log \frac{9.16^6}{K_1^2}  
+   K_3
\right).
\end{equation}
Recalling \eqref{Del.mu.alpha.ineq}, we obtain that there exists $\Lambda_5 \geq \Lambda_4$ such that for all $\alpha \geq \Lambda_5$, we have $\mu \geq 1/s_0^2$ and 
\begin{equation}
\frac{\|(\la-L)u\|_{L^2}^2}{\|u\|_{L^2}^2} 
\leq 
\exp \left(-\frac {K_2}2 \mu \right);
\end{equation}
thus \eqref{pseudo.est} and consequently \eqref{L.res.lb.proof} follow with $\wt \eta = K_2/4$.
\hfill \qedsymbol

%\newpage

\printbibliography

\end{document}